\documentclass[10pt,reqno]{amsart}
\usepackage[backend=bibtex,isbn=false,style=alphabetic]{biblatex}
\bibliography{literature}

\usepackage{a4wide}

\usepackage{amsmath,amsfonts,amsthm,amssymb,amsxtra,dsfont,mathtools,bm}
\usepackage{amsxtra, amssymb, mathrsfs, pstricks}
\usepackage[english]{babel}
\usepackage{amsmath,amsfonts,amstext,amsbsy,amsopn,amssymb,amsthm,amscd}
\usepackage{amsxtra,latexsym}
\usepackage{exscale}
\usepackage{graphicx,mathrsfs}
\usepackage{psfrag}
\usepackage{paralist,eucal,enumerate}
\usepackage{color,graphics}

\usepackage{subcaption}

\newtheorem{theorem}{Theorem}[section]
\newtheorem{proposition}[theorem]{Proposition}

\newtheorem{lemma}[theorem]{Lemma}
\newtheorem{corollary}[theorem]{Corollary}

\theoremstyle{definition}

\newtheorem{definition}[theorem]{Definition}

\theoremstyle{remark}
\newtheorem{remark}[theorem]{\bf Remark}

\numberwithin{equation}{section}

\usepackage{graphicx,mathrsfs}
\newcommand{\dd}{\mathrm{d}}

\DeclareMathOperator{\sym}{sym}

\newcommand{\field}[1]{\mathbb{#1}}
\newcommand{\N}{\field{N}}
\newcommand{\R}{\field{R}}

 \def\bbN{{\mathbb N}} 
  \def\bbR{{\mathbb R}}

 \def\sfT{{\mathsf T}}

 \def\calE{{\mathcal E}} \def\calF{{\mathcal F}}
\def\calG{{\mathcal G}}

 \def\calQ{{\mathcal Q}} \def\calR{{\mathcal R}}
  \def\calU{{\mathcal U}}
\def\calV{{\mathcal V}}  \def\calX{{\mathcal X}}
 \def\calZ{{\mathcal Z}}

\DeclareMathOperator{\tr}{tr}

\DeclareMathOperator{\dom}{dom}
\DeclareMathOperator{\dist}{dist}

\DeclareMathOperator{\argmin}{argmin}

\DeclareMathOperator{\Lin}{Lin}

\makeatletter%
\newcommand\xrupharpoonast[2][]{%
	\ext@arrow 0099{\rupharpoonfill@}{#1}{#2}\hspace{-4pt}{}^\ast\hspace{3pt}}
\def\rupharpoonfill@{%
	\arrowfill@\relbar\relbar\rightharpoonup}
\makeatother
\makeatletter%
\newcommand\xrupharpoon[2][]{%
	\ext@arrow 0099{\rupharpoonfill@}{#1}{#2}}
\def\rupharpoonfill@{%
	\arrowfill@\relbar\relbar\rightharpoonup}
\makeatother
\makeatletter%
\newcommand\xrarrow[2][]{%
	\ext@arrow 0099{\rarrowfill@}{#1}{#2}}
\def\rarrowfill@{%
	\arrowfill@\relbar\relbar\longrightarrow}
\makeatother
\newcommand{\rtok}{\xrarrow{k\to\infty}}
\newcommand{\rton}{\xrarrow{n\to\infty}}

\newcommand{\rtom}{\xrarrow{m\to\infty}}
\newcommand{\rtol}{\xrarrow{l\to\infty}}
\newcommand{\rhupk}{\xrupharpoon{k\to\infty}}
\newcommand{\rhupn}{\xrupharpoon{n\to\infty}}

\newcommand{\rhupm}{\xrupharpoon{m\to\infty}}
\newcommand{\rhupastk}{\xrupharpoonast{k\to\infty}}
\newcommand{\rhupastn}{\xrupharpoonast{n\to\infty}}

\newcommand{\abs}[1]{\left\lvert#1\right\rvert}      

\newcommand{\babs}[1]{\bigl\vert#1\bigr\vert} 

\newcommand{\norm}[1]{\left\lVert#1\right\rVert}

\newcommand{\bnorm}[1]{\bigl\Vert#1\bigr\Vert} 

\newcommand{\Set}[2]{\left\{\,#1\,:\,#2\,\right\}}

\definecolor{darkgreen}{rgb}{0.0, 0.5, 0.0}

\usepackage[ruled,linesnumbered]{algorithm2e}
\usepackage{eucal}

\renewcommand{\dd}[1]{\ensuremath{\operatorname{d}\!{#1}}}
\newcommand{\D}{\operatorname{D}}

\newcommand{\strain}{\bm{\varepsilon}}
\newcommand{\boldC}{\textbf{C}}

\begin{document}
	\DeclareRobustCommand{\subtitle}[1]{\\#1}
	\title[Approximation of BV solutions of a rate-independent damage model]{Approximation of balanced viscosity solutions of a rate-independent damage model by combining alternate minimization with a local minimization algorithm}
	
	\author{Samira Boddin}
	\address{Samira Boddin, Institute of Mathematics, University of Kassel,
		Heinrich-Plett Str.~40, 34132 Kassel, Germany. Phone:
		+49 0561 8044376
	}
	\email{sboddin@mathematik.uni-kassel.de}
	
	\author{Felix Rörentrop}
	\address{Felix Rörentrop, Institute of Mechanics, TU Dortmund,
		Leonhard-Euler-Str. 5, 44227 Dortmund, Germany. Phone:
		+49 231 7556843 }
	\email{felix.roerentrop@tu-dortmund.de}
	
	\author{Dorothee Knees}
	\address{Dorothee Knees, Institute of Mathematics, University of Kassel,
		Heinrich-Plett Str.~40, 34132 Kassel, Germany. Phone:
		+49 0561 8044355}
	\email{dknees@mathematik.uni-kassel.de}
	
	\author{Jörn Mosler}
	\address{Jörn Mosler, Institute of Mechanics, TU Dortmund,
		Leonhard-Euler-Str.~5, 44227 Dortmund, Germany. Phone:
		+49 231 7555744}
	\email{joern.mosler@tu-dortmund.de}

\date{}

\maketitle

\begin{abstract}
	The modeling of cracks has been an intensely researched topic for decades -- both from the mechanical as well as from the mathematics point of view. As far as the modeling of sharp cracks/interfaces is concerned, the resulting free boundary problem is numerically very challenging. For this reason, diffuse approximations in the sense of phase-field theories have become very popular. Within this paper, the focus is on rate-independent damage models. Since the  resulting phase-field energies in general are non-convex, we are faced with a  discontinuous evolution of the phase-field variable. Solution concepts have to be carefully chosen in order to predict discontinuities that are physically reasonable. We focus here on the concept of balanced viscosity solutions and develop a convergence scheme that combines alternate minimization with a local minimization ansatz due to Mielke/Efendiev, \cite{EfendievMielke2006}. We proof the convergence of the incremental solutions to balanced viscosity solutions. Moreover, the discretization concept is implemented  and several carefully selected  examples show the performance of this combined approach. Particularly, the effect of different norms and arc-length parameters in the local minimization scheme is  investigated.
\end{abstract}

{\small



}



\setcounter{tocdepth}{3}
{\small \tableofcontents}


	\section{Introduction}
\label{s:introduction}
The modeling of cracks has been an intensely researched topic for decades -- both from the mechanical as well as from the mathematics point of view. As far as the modeling of sharp cracks/interfaces is concerned, the resulting free boundary problem is numerically very challenging. For this reason, diffuse approximations in the sense of phase-field theories have become very popular, cf. \cite{BFM2008}. Within this paper, the focus is on a rate-independent damage model of Ambrosio-Tortorelli type. The underlying energy is not simultaneously convex in the displacement field and the phase-field variable representing the damage state.   
Due to this non-convexity solutions typically are discontinuous in time even if the applied loads depend smoothly on the time parameter. 

Starting with the seminal paper \cite{MielkeTheilLevitas2002} different solution concepts for rate-independent systems with non-convex energies were proposed and in general they are not equivalent. The concept of (global) energetic solutions (GES) based on a global energy minimization principle is now well established. For damage models we mention the papers \cite{ThomasMielke2010,MielkeRoubicek2006,Giacomini2005} 
where the existence of GES is shown. In addition in \cite{Giacomini2005} the $\Gamma$-convergence to the Francfort/Marigo fracture model is proved. 
However, the evolution predicted by the GES is questionable from a physical point of view. Indeed, GES may develop a discontinuity even though  a local  criterion of Griffith type  predicts a stable evolution. In contrast to this, the class of balanced viscosity solutions (BV solutions) seems to predict discontinuities that better reflect the physical behavior of the underlying system.   
For this reason, we focus here on balanced viscosity solutions. BV solutions for instance are obtained as inviscid limit of viscously regularized versions of the original rate-independent system and we refer to \cite{Jumps2008,MRS16}
for abstract results about BV solutions and to \cite{KRZ,KRZ15}, 
where they are discussed for damage problems. 
An illustrative example showing the different predictions of GES and BV solutions is presented in  \cite{KneesSchroeder2013}, where the evolution of a single crack is analyzed.  Numerical algorithms are needed that reliably approximate that type of solution one is interested in.

In this paper, we analyze an algorithm  that combines a time-adaptive, incremental local minimization scheme  (originally proposed in \cite{EfendievMielke2006}) 
with an alternate minimization scheme. To be more specific, we consider a   phase-field energy of Ambrosio-Tortorelli type \cite{AmbrosioTortorelli1992}
\begin{align}\label{energy functional F}
	\calF(t,u,z) = 
	\frac{1}{2}\int_\Omega \left(z^2+\eta\right)\boldC(x) \strain(u)\colon\strain(u)\dd{x}+
	\underbrace{g_c\int_\Omega \frac{\left( 1-z\right)^2}{4\theta}+\theta\abs{\nabla z}^2
		\dd{x}}_{=: \calG_f(z)}
	-\left\langle \ell(t),u\right\rangle
\end{align}
where $u\in\calU$ is the \textit{displacement field} with \textit{linearized strain} $\strain(u)$ and \textit{elasticity tensor} $\boldC$, $z\in\calZ$ is the scalar \textit{phase-field damage variable} with $0\leq z\leq 1$, $z=1$ corresponding to no damage and $z=0$ corresponding to maximum damage, $\ell(t)$ is a linear functional representing an \textit{external loading}, $g_c$ is \textit{fracture toughness} while $\theta>0$ is the phase-field \textit{length scale parameter} and $\eta>0$ is a lower bound for the elasticity coefficient often used to model incomplete damage. The latter two parameters are typically very small and they will be kept fixed throughout the paper.  
For notational details about our function spaces $\calU$ and $\calZ$ and the precise set-up we refer to section~\ref{sec:set-up and notation}.  
Observe that $\calF(t,\cdot,\cdot)$ in general is not convex in the pair $(u,z)$ but it is separately convex (and even quadratic) in $u$ if we keep $z$ fixed and vice versa. 

When analyzing rate-independent phase-field models for damage evolution in elastic bodies we basically face two difficulties. First, as already emphasized, the energy functional we consider is not convex but  only separately convex in the two variables $u$ and $z$. Therefore it is numerically challenging to minimize with respect to both variables at once. The second difficulty arises from the rate-independent nature of the model leading to discontinuous solutions as explained above. 

For each of the two problems there can be found some approaches to solve them. Regarding the two field problem two convenient approaches would be either to switch (numerically) to the reduced energy functional, which then depends for instance  on the damage variable only and not on the displacement field, or as a second approach to minimize alternately with respect to damage and displacement field with the aim to calculate at least critical points. With respect to rate-independent problems there exist basically two approaches. Either one could perform the vanishing viscosity procedure or one could use the discretization scheme with adaptive time-step size proposed in \cite{EfendievMielke2006}.

In summary this leads to four combined approaches for solving both problems at once and we give here a short overview on the literature. The vanishing viscosity approach is followed in \cite{KRZ} where the reduced energy is considered and in \cite{Almi2020} where it is combined with an alternate minimization scheme. The example in \cite{KneesSchroeder2013} 
shows that although convergence is guaranteed analytically, from a computational point of view it is difficult to choose the discretization parameters (in particular to scale the viscosity parameter correctly with the time increment) in a way such that the correct behavior is visible already for rather coarse discretizations. 
It seems that the vanishing viscosity procedure is rather appropriate to provide an abstract existence result than for the numerical approximation. 

The time-adaptive local minimization algorithm proposed by \cite{EfendievMielke2006}  provides an alternative to the vanishing viscosity approach. We refer to \cite{Knees2018}, 
where a convergence proof for a semilinear problem was given based on a reduced energy. 
In \cite{Sievers2021} this approach was tested numerically, again based on the reduced energy: 
the energy is first discretized and then the discretized energy is reduced in the discrete state space. Unfortunately, a rigorous convergence analysis of this method is missing in \cite{Sievers2021}. 
Numerical experiments suggest that the time-adaptive local minimization algorithm is better suited for the practical approximation of BV solutions.

Let us finally comment on schemes that are based on a time-discretization (with fixed time-increment) combined with alternating minimization: Motivated by the separately convex structure of the energy $\calF(t,\cdot,\cdot)$, time incremental alternate minimization algorithms or staggered schemes are commonly used in practice, see for instance \cite{MarigoMaurini2016,Miehe2010}. Here, at a given point $t_k$ in time one alternately minimizes $\calF(t_k,\cdot,\cdot)$ in $u$ and $z$ until a stopping criterion is reached (the aim is to find a critical point of $\calF(t_k,\cdot,\cdot)$) and this procedure is then repeated at the next time increment. It was shown in \cite{KN2017} 
and later also in \cite{AlmiNegri2020}
that discrete solutions obtained with this procedure converge to BV-like solutions. However there is a subtle difference between the limits of this version of the  alternate minimization procedure and BV solutions in their original definition: In case of BV solutions, the transition curve connecting the starting point of a jump with  the end point of the jump   is given by a viscous equation for the damage variable, where the viscosity corresponds to the one from the vanishing viscosity approach or is induced by the norm chosen for the local minimization algorithm. Jump trajectories created by   pure alternate minimization schemes (without additional viscosity or locality) are characterized by viscous equations which include both, visco-damage as well as visco-elastic effects.  This might finally  lead to completely different predictions of the different numerical schemes.

As already mentioned, in this paper we combine the time-adaptive local minimization scheme from \cite{EfendievMielke2006} with an alternate minimization ansatz in order to calculate suitable critical points (fixed points of a certain minimization problem) in each discretization step. We prove the convergence of suitable interpolants to BV solutions and complete the paper with numerical examples that illustrate the behavior of this combined approach.  We also show the influence of different choices for the norms as well as for the arc-length parameter in the local minimization scheme. From the analytic point of view we combine and extend the arguments 
from \cite{KN2017,KRZ} 
and \cite{Knees2018,Sievers2020} 
in order to obtain our convergence result.

	\section{Set-up and notation}\label{sec:set-up and notation}

	Hereafter $\Omega\subset\bbR^d$ with $d=2$ is a bounded domain with Lipschitz boundary and $\partial\Omega=\overline{\Gamma_D}\cup\overline{\Gamma_N}$ with $\Gamma_D\cap\Gamma_N=\emptyset$, where $\Gamma_D$ is an open Dirichlet boundary with Hausdorff measure $\CMcal{H}^{d-1}\big(\Gamma_D\big)>0$ and $\Gamma_N$ is a Neumann boundary.
	Additionally we shall assume that $\Gamma_D$ and $\Gamma_N$ are regular in the sense of Gröger (comp. \cite{HMW2011}).
	
	\subsection{Function spaces}\label{sec:spaces} We denote by $u:\Omega\to\R^d$ the displacement field and by $z:\Omega\to\bbR$ the (scalar) damage/phase-field variable. 
	The corresponding state spaces are
	\begin{align*}
		\calU=\Set{u\in H^1\big(\Omega,\bbR^d\big)}{u|_{\Gamma_D}=0}=W_{\Gamma_D}^{1,2}\big(\Omega,\bbR^d\big)\qquad
		\text{and}\qquad\calZ=H^1(\Omega).
	\end{align*}
	For $z\in\calZ$ we use the notation
	\begin{align*}
		\norm{z}_\calZ^2=\norm{z}_{L^2(\Omega)}^2+\norm{\nabla z}_{L^2(\Omega,\R^d)}^2.
	\end{align*}
	Notice the compact embedding $\calZ\hookrightarrow\hookrightarrow L^r (\Omega)$ for every $r\in \left[1,\infty\right)$ given for $d=2$ by Rellich-Kondrachov's theorem \cite[Theorem~12.44]{Leoni2017}. We will further set $\calV=L^\alpha(\Omega)$ with $\alpha>3$ specified in \eqref{alpha} and $\calX=L^1(\Omega)$. We thus have the following dense and compact respectively continuous embeddings
	\begin{align}\label{assume:embed}
		\calZ\hookrightarrow\hookrightarrow\calV\hookrightarrow \calX.
	\end{align}
	As mentioned in the introduction the physically reasonable range for the damage/phase-field variable $z$ is $\left[0,1\right]$, where $z(x)=1$ corresponds to no damage and $z(x)=0$ corresponds to maximum damage at $x\in\Omega$. Therefore we will often use the notation
	\begin{align*}
		\calZ_{\left[0,1\right]}\coloneq\Set{z\in\calZ}{z\in\left[0,1\right]\text{ a.e. in }\Omega}.
	\end{align*}
	With regard to the improved integrability of the displacement and the stresses stated in the paragraph below, we will additionally use for $p\in\left(1,\infty\right)$ the notation 
	\begin{align*}
		\calU^p&\coloneq W_{\Gamma_D}^{1,p}\bigl(\Omega,\bbR^d\bigr)=\Set{u\in W^{1,p}\big(\Omega,\bbR^d\big)}{u|_{\Gamma_D}=0}\\
		\text{with dual space}\quad 
		\left(\calU^{p}\right)^\ast&=\left(W_{\Gamma_D}^{1,p}\bigl(\Omega,\bbR^d\bigr)\right)^\ast=W_{\Gamma_D}^{-1,p'}\big(\Omega,\bbR^d\big)\quad\text{for }p'>1\text{ with }\frac{1}{p}+\frac{1}{p'}=1.
	\end{align*}

	\subsection{Energy functional}\label{sec:energy} The energy is given by $\calE:\left[0,T\right]\times \calU\times\calZ\to\R_{\infty}$ with
	\begin{align}\label{energy functional}
		\calE(t,u,z) = \underbrace{\frac{\kappa}{2}\norm{z}_{\calZ}^2}_{\eqcolon \calE_1(z)}+\underbrace{\frac{1}{2}\int_\Omega \left(z^2+\eta\right)\bigl(\boldC \strain(u)\bigr):\strain(u)\dd{x}-\left\langle \ell(t),u\right\rangle_{\calU^\ast,\calU}}_{\eqcolon \calE_2(t,u,z)}.
	\end{align}
	Note that for $u\in\calU$ and $z\in\calZ$ the energy $\calE(t,u,z)$ does not have to be finite. But if we choose $u\in\calU^{\tilde{p}}$ or $z\in\calZ_{\left[0,1\right]}$, it is always finite.
	
	Here $\boldC\in L^\infty \big(\Omega, \Lin\big(\bbR_{\sym}^{d\times d},\bbR_{\sym}^{d\times d}\big)\big)$ is the \textit{elastic tensor} and is supposed to fulfill the coercivity condition
	\begin{align}\label{assume:elastictensor}
		\exists \gamma>0,\,\forall \xi\in \bbR_{\sym}^{d\times d}, \text{ f.a.a. } x\in\Omega: \quad (\boldC(x)\xi):\xi\geq \gamma\abs{\xi}^2 
	\end{align}
	and the symmetry condition
	\begin{align}\label{assume:symmetrictensor}
		\forall \xi_1,\xi_2\in \bbR_{\sym}^{d\times d}, \text{ f.a.a. } x\in\Omega: \quad (\boldC(x)\xi_1)\colon\xi_2=(\boldC(x)\xi_2)\colon\xi_1.
	\end{align}
	By $A \colon B=\tr\bigl(B^\mathsf{T}A\bigr)$ we denote the Frobenius inner product of two matrices $A, B\in\R^{d\times d}$.
	Furthermore, $\strain(u)=\frac{1}{2}\left(\nabla u+\left(\nabla u\right)^\sfT \right)$ is the \textit{symmetrized strain tensor} and
	$\ell\in C^{1,1}(\left[0,T\right],\calU^\ast)$ is a given \textit{external loading}. 
	The small parameter $\eta>0$ is on the one hand often used to model incomplete damage of elastic materials \cite{KRZ}, since the lower bound on the coefficient of the elastic energy density means, that even if the material is completely damaged, it is still capable to store some elastic energy. On the other hand we need this lower bound to have coercivity of the energy functional with respect to the displacement field.
	For modeling complete damage we would need to send $\eta$ to zero, which we do not in this paper. Moreover, a constant $\kappa>0$ is given.

	Now we will show some coercivity estimates for $\calE$, which especially imply that $\calE$ is bounded from below.
	\begin{lemma}[Coercivity of the energy functional]\label{lemma:coercivity}
		The energy functional $\calE$ given by \eqref{energy functional} is bounded from below and satisfies the following coercivity estimate
		\begin{align}\label{coercivity}
			\exists c_0,c_1,c_2>0:\quad \calE(t,u,z)\geq c_1\norm{u}_\calU^2+c_2\norm{z}_\calZ^2-c_0.
		\end{align}
	\end{lemma}
	\begin{proof}
		Choose $c_2\coloneq \frac{\kappa}{2}$, then $\calE_1(z)\geq c_2\norm{z}_\calZ^2$ for all $z\in\calZ$. By \eqref{assume:elastictensor} and since $\ell\in C^{1,1}(\left[0,T\right],\calU^\ast)$ we can find $c_0,c_1,c_\ell>0$ such that for all $t\in\left[0,T\right]$, $u\in\calU$ and $z\in\calZ$
		\begin{align}\label{proof:prop_be4}
			\calE_2(t,u,z)\geq \frac{\eta\gamma}{2} \norm{\strain(u)}_{L^2(\Omega)}^2-\underbrace{\norm{\ell(t)}_{\calU^*}}_{\leq c_\ell}\norm{u}_{\calU}
			\geq  c\norm{u}_\calU^2-c_\ell\norm{u}_\calU\geq c_1\norm{u}_\calU^2-c_0,
		\end{align}
		where we have used in the second inequality that by Korn's first inequality and the Poincaré inequality there exists $c>0$ such that
		\begin{align*}
			\norm{\strain(u)}_{L^2(\Omega)}\geq c\norm{\nabla u}_{L^2(\Omega)}\geq c\norm{u}_\calU.
		\end{align*}
	\end{proof}

	\subsection{Dissipation potential}\label{sec:dissipation} For a given constant $\kappa>0$ and $v\in\calZ$ the \textit{unidirectional} dissipation potential $\calR:\calX\to[0,\infty]$ is defined by
	\begin{align*}
		\calR\left(v\right)\coloneq
		\begin{cases}
			\int_\Omega \kappa\abs{v(x)}\dd{x}, & \text{if }v\leq 0\text{ a.e. in }\Omega\\
			\infty, & \text{ otherwise.}
		\end{cases}
	\end{align*}
	This is obviously a proper convex function. Furthermore, $\calR$ is positively homogeneous of degree one, what implies the rate-independence in the damage variable $z$ of the evolutionary system driven by the functionals $\calE$ and $\calR$ \cite[section~1.3.1]{MielkeRoubicek2015}. Obviously it holds
		\begin{align}
			\forall v\in\calX:\quad \kappa\norm{v}_\calX\leq \calR(v).\label{R2}
		\end{align}
	
	\subsection{Improved integrability of displacement and stresses}
	All results obtained by the approaches mentioned in the introduction are based on the following result from \cite[Theorem~1.1 and Proposition~1.2]{HMW2011}. It provides an improved integrability for the solution of the elliptic Euler-Lagrange equation corresponding to the minimization of the considered energy functional with respect to the displacement $u$ and fixed phase-field $z$.
	\begin{theorem}[Improved integrability of displacement and stresses]\label{thm:improved integrability}
		Let $M>0$ be given. Under the above assumptions and for $z\in \calZ$ with $\norm{z}_{L^\infty(\Omega)}\leq M$ let $L_z:\calU\to\calU^\ast$ be the linear elliptic operator defined by
		\begin{align*}
			\forall v,w\in\calU:\quad \left\langle L_z(v),w\right\rangle_{\calU^\ast,\calU} \coloneq \int_\Omega \left(z^2+\eta\right) \left(\normalfont\boldC\strain (v)\right):\strain(w)\dd{x}.
		\end{align*}
		Then
		\begin{align}
			\begin{split}
				\exists p>2,\,\forall \tilde{p}\in\left[2,p\right],\,\forall z\in \calZ\text{ with } \norm{z}_{L^\infty(\Omega)}&\leq M:\\\label{thm:ii1}
				&L_z:\calU^{\tilde{p}}\to\bigl(\calU^{\tilde{p}'}\bigr)^\ast\quad \text{is an isomorphism}
			\end{split}
		\end{align}
		and
		\begin{align}
			\begin{split}
				\exists C_M>0,\,\forall z\in \calZ \text{ with } \norm{z}_{L^\infty(\Omega)}\leq M,\,\forall\tilde{p}\in\left[2,p\right],\,\forall b\in&\bigl(\calU^{\tilde{p}'}\bigr)^\ast:\\\label{thm:ii2}
				&\norm{L_z^{-1}b}_{\calU^{\tilde{p}}}\leq C_M \norm{b}_{\left(\calU^{\tilde{p}'}\right)^\ast}.
			\end{split}
		\end{align}
	\end{theorem}
	Let $p>2$ be given according to Theorem~\ref{thm:improved integrability} with $M=1$. This choice is related to the fact that the phase-field variable $z$ will stay in the physically reasonable range $\calZ_{\left[0,1\right]}$. For $\calV=L^\alpha(\Omega)$ we assume
	\begin{align}\label{alpha}
		\alpha\geq\frac{3p}{p-2}.
	\end{align}
	Observe that $\frac{3p}{p-2}>3$. Let us give a short comment on this assumption. Given a $\tilde{p}\in\left(2,p\right)$ we obtain the continuous dependencies from Lemma~\ref{lemmaA1} with $r_1=\frac{\tilde{p}p}{p-\tilde{p}}$ and from Lemma~\ref{lemmaA2EM} with $r_2=\frac{\tilde{p}}{\tilde{p}-2}$. Within the convergence analysis we need $\alpha\geq \max\{r_1,r_2\}$. Thus the minimal possible $\alpha$ is $\frac{3p}{p-2}$, which is available for 
	\begin{align}\label{tildep}
		\tilde{p}=\frac{3p}{p+1}\in\left(2,p\right).
	\end{align}
	
	\subsection{Refined assumptions on the external loading}
	In view of the above improved integrability result hereafter we shall require that
	\begin{align}\label{assume:external}
		\ell \in C^{1,1}\Big([0,T];\bigl(\calU^{p'}\bigr)^\ast\Big)\quad\text{with } p>2 \text{ from \eqref{thm:ii1} and }\frac{1}{p}+\frac{1}{p'}=1.
	\end{align}
	Note that this is a stronger assumption than $\ell \in C^{1,1}\bigl([0,T];\calU^\ast\bigr)$, since we get from $p>2$ that $p'<2$ and thus $\calU\hookrightarrow\calU^{p'}$ respectively $\bigl(\calU^{p'}\bigr)^\ast\hookrightarrow\calU^\ast$.
	
	\subsection{Continuity properties of the energy functional and the dissipation potential}
	
	\begin{lemma}[Continuity properties of energy and dissipation]\label{lemma:continuity}
		Let $\tilde{p}\in\left[2,p\right]$ with $p>2$ from \eqref{thm:ii1}. The energy functional $\calE$ given by \eqref{energy functional} is weakly lower semicontinuous on $\left[0,T\right]\times\calU^{\tilde{p}}\times\calZ$ and the dissipation potential $\calR$ is weakly lower semicontinuous on $\calZ$. 
		Furthermore, the energy functional $\calE$ is continuous on $\left[0,T\right]\times\calU^{\tilde{p}}\times\calZ$.
	\end{lemma}
	\begin{proof}
		First we know that $\calE_1$ is weakly lower semicontinuous on $\calZ$ as convex and continuous function. Now assume for $\left(t_n\right)_{n\in\N}\subset\left[0,T\right]$, $\left(u_n\right)_{n\in\N}\subset \calU^{\tilde{p}}$ and $\left(z_n\right)_{n\in\N}\subset \calZ$ that
		\begin{align*}
			t_n\rton t\quad\text{in }\left[0,T\right],\quad u_n\rhupn u\quad\text{in } \calU^{\tilde{p}} \quad\text{and}\quad z_n\rhupn z\quad\text{in } \calZ.
		\end{align*}
		Then secondly by \eqref{assume:external} we have that $\ell(t_n)\rton\ell(t)$ in $\bigl(\calU^{p'}\bigr)^\ast$. Notice that $p'<2<\tilde{p}$ such that $\bigl(\calU^{p'}\bigr)^\ast\hookrightarrow\bigl(\calU^{\tilde{p}}\bigr)^\ast$ and thus  $\ell(t_n)\rton\ell(t)$ in $\bigl(\calU^{\tilde{p}}\bigr)^\ast$. Together with $u_n\rhupn u$ in $\calU^{\tilde{p}}$ and a well known convergence principle (see for example \cite[p.~58, Lemma~0.3(ii)]{Ruzicka2004}) this implies
		\begin{align*}
			\langle \ell(t_n),u_n\rangle\rton \langle \ell(t),u\rangle.
		\end{align*}
		Thirdly we obtain by adding a zero
		\begin{multline}\label{term}
			\int_{\Omega}\left(z_n^2+\eta\right)\left(\boldC\strain(u_n)\right)\colon\strain(u_n)\dd{x}\\
			=
			\int_{\Omega}\left(z_n^2-z^2\right)\left(\boldC\strain(u_n)\right)\colon\strain(u_n)\dd{x}+\int_{\Omega}\left(z^2+\eta\right)\left(\boldC\strain(u_n)\right)\colon\strain(u_n)\dd{x}.
		\end{multline}
		Hölder's inequality with $\frac{2}{2r_2}+\frac{2}{\tilde{p}}=1$ and the compact embedding $\calZ\hookrightarrow\hookrightarrow L^{2r_2}(\Omega)$ leads together with the assumed convergences to
		\begin{align*}
			\abs{\int_{\Omega}\left(z_n^2-z^2\right)\left(\boldC\strain(u_n)\right)\colon\strain(u_n)\dd{x}}&=\abs{\int_{\Omega}\left(z_n-z\right)\left(z_n+z\right)\left(\boldC\strain(u_n)\right)\colon\strain(u_n)\dd{x}}\\
			&\leq C\norm{z_n-z}_{L^{2r_2}(\Omega)}\norm{z_n+z}_{L^{2r_2}(\Omega)}\norm{\strain(u_n)}_{L^{\tilde{p}}(\Omega)}^2\rton 0.
		\end{align*}
		For the second term on the right hand side of \eqref{term} we notice that the mapping
		\begin{align*}
			\calU^{\tilde{p}}\ni w\mapsto \int_{\Omega}\underbrace{\left(z^2+\eta\right)}_{\in L^{r_2}(\Omega)}\left(\boldC\strain(w)\right)\colon\strain(w)\dd{x}
		\end{align*}
		is convex (by \eqref{assume:elastictensor} and Korn's inequality, its second Fréchet derivative is positive definite) and that it is continuous on $\calU^{\tilde{p}}$. Therefore it is weakly lower semicontinuous on $\calU^{\tilde{p}}$. This handles the remaining term of $\calE$, thus finally the weak lower semicontinuity of $\calE$ on $\left[0,T\right]\times\calU^{\tilde{p}}\times\calZ$ is proved.

		By the compact embedding $\calZ\hookrightarrow\hookrightarrow \calX$ the weak lower semicontinuity of $\calR$ on $\calZ$ follows, if we can show lower semicontinuity of $\calR$ on $\calX=L^1(\Omega)$, i.e. that for $\left(v_n\right)_{n\in\N}\subset\calX$ and $v\in\calX$
		\begin{align}\label{proof:continuity1}
			v_n\rton v\quad\text{in } \calX \quad\text{implies}\quad \liminf_{n\to\infty}\calR(v_n)\geq \calR(v).
		\end{align}		
		Let us recap that $\dom(\calR)=\Set{v\in\calX}{v\leq 0\text{ a.e. in }\Omega}$ and $\calR(v)=\kappa\norm{v}_\calX$ for all $v\in\dom(\calR)$, so obviously $\calR$ is on $\dom(\calR)$ continuous with respect to $\calX$.
		If $\liminf_{n\to\infty}\calR(v_n)=\infty$ then implication \eqref{proof:continuity1} is obviously satisfied. Otherwise there exists a subsequence $\left(v_{n_k}\right)_{k\in\N}\subset\dom(\calR)$ such that
		\begin{align*}
			\lim_{k\to\infty}\calR\big(v_{n_k}\big)=\liminf_{n\to\infty}\calR(v_n)<\infty.
		\end{align*}
		Now in view of the continuity of $\calR$ on $\dom(\calR)$ with respect to $\calX$ it suffices to prove that $\dom(\calR)$ is closed in $\calX$. Therefore let $\left(v_n\right)_{n\in\N}\subset\dom(\calR)$ and $v\in\calX$ with
		$v_n\rton v$ in $\calX$. Suppose that $v\notin\dom(\calR)$. Then there exists $\tilde{\Omega}\subset\Omega$ with $\big\vert\tilde{\Omega}\big\vert>0$ such that $v>0$ on $\tilde{\Omega}$, so
		\begin{align*}
			0<\int_{\tilde{\Omega}}v(x)\dd{x}
			\leq \int_{\tilde{\Omega}}v(x)-\underbrace{v_{n}(x)}_{\leq 0}\dd{x}
			\leq
			\int_{\Omega}\abs{v(x)-v_{n}(x)}\dd{x}=\norm{v-v_{n}}_{L^1(\Omega)}\rton 0
		\end{align*}
		a contradiction. Thus $v\in\dom(\calR)$ and therefore $\dom(\calR)$ is closed in $\calX$.
		Finally the continuity of $\calE$ on $\left[0,T\right]\times\calU^{\tilde{p}}\times\calZ$ can be shown using Hölder's inequality with $\frac{1}{r_2}+\frac{2}{\tilde{p}}=1$ and the continuous embeddings $\calZ\hookrightarrow L^{2r_2}(\Omega)\hookrightarrow L^{r_2}(\Omega)$.
	\end{proof}

\subsection{Partial derivatives of the energy functional}
\begin{lemma}[Partial derivatives of the energy]\label{lemma:partialderivatives}
	Let $\calE$ be given by \eqref{energy functional} and let $t\in\left[0,T\right]$ be fixed.
	\begin{enumerate}[(a)]
		\item For fixed $u\in\calU^{\tilde{p}}$ the function $\calE(t,u,\cdot):\calZ\to\R$ is Fréchet differentiable in $z\in\calZ$ with
		\begin{align}\label{frechetz}
			\langle \D_z\calE(t,u,z),v\rangle =\kappa \int_{\Omega} zv+\nabla z\cdot\nabla v\dd{x} +\int_{\Omega} zv\left(\normalfont\boldC\strain(u)\right)\colon\strain(u)\dd{x}\quad \forall v\in\calZ.
		\end{align}
		\item For fixed $z\in\calZ$ the function $\calE(t,\cdot,z):\calU^{\tilde{p}}\to\R$ is Frèchet differentiable in $u\in\calU^{\tilde{p}}$ with
		\begin{align}\label{frechetu}
			\langle \D_u\calE(t,u,z),w\rangle=\int_{\Omega}\bigl(z^2+\eta\bigr)\left(\normalfont\boldC\strain(u)\right)\colon\strain(w)\dd{x}-\langle \ell(t),w\rangle\quad\forall w\in\calU^{\tilde{p}}.
		\end{align}
	\end{enumerate}
	Moreover, for every $u\in\calU^{\tilde{p}}$ the Fréchet derivative $\D_z\calE(t,u,\cdot)$ is continuous on $\calZ$ and for every $z\in\calZ$ the Fréchet derivative $\D_u\calE(t,\cdot,z)$ is continuous on $\calU^{\tilde{p}}$. Together with the assumptions on the external loading $\ell$ from \eqref{assume:external} we eventually obtain $\calE\in C^1\left(\left[0,T\right]\times\calU^{\tilde{p}}\times\calZ,\R\right)$.
\end{lemma}
\begin{proof}
	\begin{enumerate}[(a)]
		\item For $z,v\in\calZ$ it holds
		\begin{align*}
			\calE(&t,u,z+v)-\calE(t,u,z)=\frac{\kappa}{2}\left(\norm{z+v}_\calZ^2-\norm{z}_\calZ^2\right)+\frac12 \int_{\Omega} \left(\left(z+v\right)^2-z^2\right)\left(\boldC\strain(u)\right)\colon\strain(u)\dd{x}\\
			&=\kappa \int_{\Omega} zv+\nabla z\cdot\nabla v\dd{x} +\int_{\Omega} zv\left(\boldC\strain(u)\right)\colon\strain(u)\dd{x}+\frac{\kappa}{2}\norm{v}_\calZ^2+\frac12\int_{\Omega}v^2\left(\boldC\strain(u)\right)\colon\strain(u)\dd{x}\\
			&=\kappa \int_{\Omega} zv+\nabla z\cdot\nabla v\dd{x} +\int_{\Omega} zv\left(\boldC\strain(u)\right)\colon\strain(u)\dd{x}+o(\norm{v}_\calZ)\quad\text{for } v\to 0 \text{ in }\calZ,
		\end{align*}
		where the last equation follows, since by Hölder's inequality with $\frac{1}{r_2}+\frac{2}{\tilde{p}}=1$ and by the continuous embedding $\calZ\hookrightarrow L^{2r_2}(\Omega)$
		\begin{align*}
			0&\leq\frac{\kappa}{2}\norm{v}_\calZ^2+\frac12\int_{\Omega}v^2\left(\boldC\strain(u)\right)\colon\strain(u)\dd{x}\\
			&\leq \frac{\kappa}{2}\norm{v}_\calZ^2+C\norm{v}_{L^{2r_2}(\Omega)}^2\norm{\strain(u)}_{L^{\tilde{p}}}^2
			\leq 
			\left(\frac{\kappa}{2}+C\norm{u}_{\calU^{\tilde{p}}}^2\right)\norm{v}_\calZ^2.
		\end{align*}
		\item For $u,w\in\calU^{\tilde{p}}$ taking into account symmetry condition \eqref{assume:symmetrictensor} we have
		\begin{align*}
			\calE(t,u+w,&z)-\calE(t,u,z)\\
			&=\frac12 \int_{\Omega} \left(z^2+\eta\right)\Bigl(\left(\boldC\strain(u+w)\right)\colon\strain(u+w)-\left(\boldC\strain(u)\right)\colon\strain(u)\Bigr)\dd{x}-\langle \ell(t),w\rangle\\
			&\hspace{-4.5pt}\overset{\eqref{assume:symmetrictensor}}{=}\int_{\Omega} \left(z^2+\eta\right)\left(\boldC\strain(u)\right)\colon\strain(w)\dd{x}-\langle \ell(t),w\rangle+\frac12 \int_{\Omega} \left(z^2+\eta\right)\left(\boldC\strain(w)\right)\colon\strain(w)\dd{x}\\
			&=\int_{\Omega} \left(z^2+\eta\right)\left(\boldC\strain(u)\right)\colon\strain(w)\dd{x}-\langle \ell(t),w\rangle+o\left(\norm{w}_{\calU^{\tilde{p}}}\right)\quad\text{for } w\to 0 \text{ in }\calU^{\tilde{p}},
		\end{align*}
		where the last equation follows, since by Hölder's inequality with $\frac{1}{r_2}+\frac{2}{\tilde{p}}=1$ and by the continuous embedding $\calZ\hookrightarrow L^{2r_2}(\Omega)$
		\begin{align*}
			\hspace{13pt}0\leq\frac12\int_{\Omega}\left(z^2+\eta\right)\left(\boldC\strain(w)\right)\colon\strain(w)\dd{x}\leq C\left(\norm{z}_{L^{2r_2}(\Omega)}^2+\eta\right)\norm{\strain(w)}_{L^{\tilde{p}}}^2\leq C\left(\norm{z}_{\calZ}^2+\eta\right)\norm{w}_{\calU^{\tilde{p}}}^2.
		\end{align*}
	\end{enumerate}
	The continuity properties of the Fréchet derivatives can immediately be shown with similar arguments.
\end{proof}
	
	\section{Approximation scheme combining Alternate Minimization with the Efendiev\&Mielke scheme}
	
	Let initial values $t_0=0$ and $z_0\in\calZ_{\left[0,1\right]}$ and an arc-length parameter $\rho>0$ be given. Moreover assume that
	\begin{align}
		\tilde{p}\in\left(2,p\right) \text{ with } p>2 \text{ from \eqref{thm:ii1}}.\label{assumption:AMp}
	\end{align}
	We consider the following scheme combining alternate minimization (AM) with the scheme proposed by Efendiev and Mielke (E\&M).
	
	\SetKwComment{Comment}{/* }{ */}
	
	\begin{algorithm}
		\caption{Approximation scheme combining AM with E\&M}\label{alg1}
		\KwData{initial values $t_0=0$ and $z_0\in \calZ_{\left[0,1\right]}$}
		$z_{0,0}^\rho\coloneq z_{-1}^\rho\coloneq z_0$\Comment*[r]{$k=0$: check for jump at $t_0$}
		$t_0^\rho\coloneq t_0$\;
		\ForAll{$k\in\bbN_0$}{
			\ForAll{$i\in\bbN$}{
				$u_{k,i}^\rho\coloneq \argmin\Set{\calE\big(t_k^\rho,u,z_{k,i-1}^\rho\big)}{u\in \calU}$\;
				$z_{k,i}^\rho\coloneq \argmin\Set{\calE\big(t_k^\rho,u_{k,i}^\rho,z\big)+\calR\big(z-z_{k-1}^\rho\big)}{z\in \calZ,\, \norm{z-z_{k-1}^\rho}_\calV\leq \rho}$\;
			}
			$\exists$ subsequence $\left(i_m\right)_{m\in\bbN}$ such that the following limits exist\Comment*[r]{see Prop.~\ref{lemma_convergentsubsequence}}\PrintSemicolon
			$u_k^\rho\coloneq \lim\limits_{m\to\infty} u_{k,i_m}^\rho\quad \text{strongly in } \calU^{\tilde{p}}$\;
			$z_{k+1,0}^\rho\coloneq z_k^\rho\coloneq \lim\limits_{m\to\infty} z_{k,i_m}^\rho\quad\text{strongly in } \calZ$\;
			$t_{k+1}^\rho\coloneq \min\left\{\, t_k^\rho+\rho-\norm{z_k^\rho-z_{k-1}^\rho}_\calV,\,\, T \,\right\}$\;
		}
	\end{algorithm}
	Here, the inner for-loop realizes the AM step, whereas the E\&M scheme and its adaptive time-incrementation is implemented within the outer for-loop. From now on we will therefore refer to the iterates of the inner loop $u_{k,i}^\rho$ and $z_{k,i}^\rho$ (for $k\in\N_0$ fixed and $i\in\N$) as AM iterates and to the iterates of the outer loop $u_k^\rho$ and $z_k^\rho$ as well as the discrete time-steps $t_k^\rho$ (for $k\in\N_0$) as E\&M iterates.
	
	\subsection{Convergence analysis for the Alternate Minimization iterates}\label{section:AM}
	Within this section we will often need the following assumptions\\
	\noindent
	\textbf{Assumptions. } 
	\begin{align}\label{assumption:AMkrhoz0}
			\text{Let } k\in\bbN_0 \text{ and } \rho>0 \text{ be fixed and suppose that } t_k\in\left[0,T\right]\text{ and } z_{k-1}^\rho\in\calZ_{\left[0,1\right]}.
	\end{align}
	\begin{align}
		\begin{split}
			&\text{Further let }\bigl(u_{k,i}^\rho\bigr)_{i\in\bbN}\subset \calU \text{ and } \bigl(z_{k,i}^\rho\bigr)_{i\in\bbN}\subset\calZ \text{ be the sequences generated by the inner}\\\label{assumption:AMiterates}
			&\text{for-loop of Algorithm~\ref{alg1}}.
		\end{split}
	\end{align}
  	\textbf{Notation.}
	Moreover, to ease notation we set for the rest of section~\ref{section:AM}
	\begin{align}\label{notation:AMz}
		\begin{split}
			z_{min}(u)&\coloneq \argmin\Set{\calE\bigl(t_k^\rho,u,z\bigr)+\calR\big(z-z_{k-1}^\rho\big)}{z\in\calZ, \norm{z-z_{k-1}^\rho}_\calV\leq \rho}\quad \text{for }u\in \calU\\
			\text{and}\hspace{30pt}&\\
			\quad u_{min}(z)&\coloneq\argmin\Set{\calE\big(t_k^\rho,u,z\big)}{u\in\calU}\quad \text{for }z\in \calZ.
		\end{split}
	\end{align}
	We remark that the minimization problems in \eqref{notation:AMz} are uniquely solvable. In particular the AM iterates $u_{k,i}^\rho$ and $z_{k,i}^\rho$ are uniquely determined. Indeed, the energy functional $\calE$ given in \eqref{energy functional} is strictly convex separately in $u$ and in $z$. Moreover, the necessary coercivity and weak lower semicontinuity assumptions of the abstract existence theorem of Tonelli from the direct method in the calculus of variations are satisfied. In particular, the set $\Set{z\in\calZ}{\norm{z-z_{k-1}^\rho}_\calV\leq\rho}$ is weakly sequentially closed in $\calZ$. Thus the existence and uniqueness of the above minimizers is guaranteed.

	The following basic estimates for the AM iterates $z_{k,i}^\rho$ are proven analogously to \cite[Prop.~5.5]{KRZ}.
	
	\begin{lemma}[Basic estimates for the AM iterates $z_{k,i}^\rho$]\label{lemma:be_AM}
		Under assumption \eqref{assumption:AMkrhoz0} we obtain
		\begin{align}\label{prop_be_AM}
			\forall i\in\bbN_0:\quad z_{k,i}^\rho\in\calZ_{\left[0,1\right]}.
		\end{align}
	\end{lemma}
	\begin{proof}
		Observe first that we have $\calR\big(z_{k,i}^\rho-z_{k-1}^\rho\big)<\infty$ for all $i\in\N$, hence
		\begin{align*}
			\forall i\in\bbN:\quad
			z_{k,i}^\rho\leq z_{k-1}^\rho\quad\text{a.e. in }\Omega.
		\end{align*}
		Since
		$z_{k-1}^\rho\leq 1$ a.e. in $\Omega$ by \eqref{assumption:AMkrhoz0}, we find
		\begin{align*}
			\forall i\in\bbN:\quad
			z_{k,i}^\rho\leq 1\quad\text{a.e. in }\Omega.
		\end{align*}
		On the other hand due to $z_{k-1}^\rho\geq 0$ we obtain
		\begin{align*}
			\norm{\big(z_{k,i}^\rho\big)^+-z_{k-1}^\rho}_\calV\leq \bnorm{z_{k,i}^\rho-z_{k-1}^\rho}_\calV\leq \rho,
		\end{align*}
		where $z^+\coloneq \max\{ z,\, 0 \}$ denotes the positive part of a function $z:\Omega\to\R$.
		Moreover it holds $\big(z_{k,i}^\rho\big)^+\in\calZ$ with 
		\begin{align}\label{proof:beAM1}
			\nabla\bigl(z_{k,i}^\rho\bigr)^+(x)=\begin{cases}
				\nabla z_{k,i}^\rho(x), &\text{if }z_{k,i}^\rho(x)>0\\
				0, & \text{if } z_{k,i}^\rho(x)\leq 0,
			\end{cases}
		\end{align}
		see for instance \cite[Satz~5.20]{Dobro2010}.
		Therefore it is admissible to test with $z=\big(z_{k,i}^\rho\big)^+$ within the AM step and thus it follows for all $i\in\bbN$
		\begin{align}\label{proof:lemma_be_AM}
			\calE\big(t_k^\rho,u_{k,i}^\rho,z_{k,i}^\rho\big)+\calR\big(z_{k,i}^\rho-z_{k-1}^\rho\big)\leq 
			\calE\Big(t_k^\rho,u_{k,i}^\rho,\big(z_{k,i}^\rho\big)^+\Big)+\calR\Big(\big(z_{k,i}^\rho\big)^+-z_{k-1}^\rho\Big).
		\end{align}
		We already know that $\bigl(z_{k,i}^\rho\bigr)^+=\max\{ z_{k,i}^\rho,\, 0\}\leq z_{k-1}^\rho$ and therefore have $z_{k,i}^\rho\leq \bigl(z_{k,i}^\rho\bigr)^+\leq z_{k-1}^\rho$. With this, \eqref{proof:beAM1} and $\babs{\big(z_{k,i}^\rho\big)^+}\leq \babs{z_{k,i}^\rho}$ one easily sees that
		\begin{align*}
			&\calR\Big(\big(z_{k,i}^\rho\big)^+-z_{k-1}^\rho\Big)\leq \calR\big(z_{k,i}^\rho-z_{k-1}^\rho\big),\\
			\norm{\nabla\bigl(z_{k,i}^\rho\bigr)^+}_{L^2(\Omega)}^2&\leq
			\bnorm{\nabla z_{k,i}^\rho}_{L^2(\Omega)}^2\quad\text{and}\quad \norm{\bigl(z_{k,i}^\rho\bigr)^+}_{L^2(\Omega)}^2\leq
			\bnorm{z_{k,i}^\rho}_{L^2(\Omega)}^2.
		\end{align*}
		Furthermore, by $\Bigl(\big(z_{k,i}^\rho\big)^+\Bigr)^2\leq \bigl(z_{k,i}^\rho\bigr)^2$ and noticing also positivity \eqref{assume:elastictensor} we obtain
		\begin{align*}
			\int_\Omega \biggl(\Bigl(\big(z_{k,i}^\rho\big)^+\Bigr)^2+\eta\biggr)\bigl(\boldC \strain\big(u_{k,i}^\rho\big)\bigr):\strain\big(u_{k,i}^\rho\big)\dd{x}\leq \int_\Omega \Bigl(\bigl(z_{k,i}^\rho\bigr)^2+\eta\Bigr)\bigl(\boldC \strain\big(u_{k,i}^\rho\big)\bigr):\strain\big(u_{k,i}^\rho\big)\dd{x}.
		\end{align*}
		Therefore, we actually have an equality in \eqref{proof:lemma_be_AM} and by uniqueness of the minimizer $z_{k,i}^\rho$ in the AM step we find $\big(z_{k,i}^\rho\big)^+=z_{k,i}^\rho$, thus $z_{k,i}^\rho\geq 0$ a.e. in $\Omega$.
	\end{proof} 
	Further we will need the following monotonicity properties.
	\begin{lemma}[Monotonicity of energy dissipation sequences]\label{lemma:monotonicity}
		Assume \eqref{assumption:AMkrhoz0} and \eqref{assumption:AMiterates} are satisfied. 
		Then $\left(\calE\big(t_k^\rho,u_{k,i}^\rho,z_{k,i}^\rho\big)+\calR\big(z_{k,i}^\rho-z_{k-1}^\rho\big)\right)_{i\in\N_0}$ is a monotonically decreasing sequence.
	\end{lemma}
	\begin{proof}
		By construction $z_{k,i}^\rho=z_{min}(u_{k,i}^\rho)$ and $u_{k,i}=u_{min}(z_{k,i-1}^\rho)$ for all $i\in\bbN$, therefore we have 
		\begin{align*}
			\forall i\in\N:\quad\calE\big(t_k^\rho,u_{k,i}^\rho,z_{k,i}^\rho\big)+\calR\big(z_{k,i}^\rho-z_{k-1}^\rho\big)&\leq \calE\big(t_k^\rho,u_{k,i}^\rho,z_{k,i-1}^\rho\big)+\calR\big(z_{k,i-1}^\rho-z_{k-1}^\rho\big)\\
			&\leq \calE\big(t_k^\rho,u_{k,i-1}^\rho,z_{k,i-1}^\rho\big)+\calR\big(z_{k,i-1}^\rho-z_{k-1}^\rho\big).
		\end{align*}
	\end{proof}
	As in \cite[Lemma~3.1]{KN2017} we derive some uniform bounds for the AM iterates.
	\begin{lemma}[Uniform bounds of the AM iterates]\label{lemma:boundsAM}
		Assume \eqref{assumption:AMkrhoz0} and \eqref{assumption:AMiterates} are satisfied and let $p>2$ be as in \eqref{thm:ii1}. Then
		 \begin{align}\label{proof:uniformbounds}
		 	\hspace{-5pt}\exists C>0,\,\forall i\in\bbN,\,\forall \tilde{p}\in\left[2,p\right]:\hspace{6pt} u_{k,i}^\rho\in\calU^{\tilde{p}} \text{ with } \bnorm{u_{k,i}^\rho}_{\calU^{\tilde{p}}}\leq C\text{ and } \bnorm{z_{k,i}^\rho}_\calZ\leq C\big(1+\norm{z_{k-1}^\rho}_\calZ \big),
		 \end{align}
	 	where $C$ can be chosen independent of $\rho>0$ and $k\in\N_0$.
	 	Moreover it holds
	 		\begin{align}\label{upperboundenergyAM}
	 			\forall i\in\N:\quad \calE\big(t_k^\rho,u_{k,i}^\rho,z_{k,i}^\rho\big)\leq \calE\big(t_k^\rho,u_{k,1}^\rho,z_{k-1}^\rho\big).
	 		\end{align}
	\end{lemma}
	\begin{proof}
		Noticing that by Lemma~\ref{lemma:be_AM} $\bnorm{z_{k,i-1}^\rho}_{L^\infty(\Omega)}\leq 1$ for all $i\in\N$, \eqref{thm:ii2} from Theorem~\ref{thm:improved integrability} implies
		\begin{align}\label{proof:ub1}
			\exists C>0,\,\forall i\in\bbN,\,\forall \tilde{p}\in\left[2,p\right]:\quad u_{k,i}^\rho\in\calU^{\tilde{p}} \text{ with } \bnorm{u_{k,i}^\rho}_{\calU^{\tilde{p}}}\leq C\norm{\ell(t_k^\rho)}_{\left(\calU^{\tilde{p}'}\right)^\ast}
		\end{align}
		and $C$ is independent of $\rho$ and $k$.
		Therefore the first estimate in \eqref{proof:uniformbounds} follows by the assumptions on $\ell$ from \eqref{assume:external}.	
		For the second estimate it follows from Lemma~\ref{lemma:monotonicity} via induction and together with $z_{k,1}^\rho=z_{min}(u_{k,1}^\rho)$ that
		\begin{align*}
			\forall i\in\N:\quad\calE\big(t_k^\rho,u_{k,i}^\rho,z_{k,i}^\rho\big)+\calR\big(z_{k,i}^\rho-z_{k-1}^\rho\big)&\leq \calE\big(t_k^\rho,u_{k,1}^\rho,z_{k,1}^\rho\big)+\calR\big(z_{k,1}^\rho-z_{k-1}^\rho\big)\\
			&\leq \calE\big(t_k^\rho,u_{k,1}^\rho,z_{k,0}^\rho\big)+\calR\big(z_{k,0}^\rho-z_{k-1}^\rho\big)
		\end{align*}
		Keeping in mind that $z_{k,0}^\rho=z_{k-1}^\rho$, $\calR(0)=0$ and $\calR(z_{k,i}^\rho-z_{k-1}^\rho)\geq 0$, we find for all $i\in\N$
		\begin{align*}
			\calE\big(t_k^\rho,u_{k,i}^\rho,z_{k,i}^\rho\big)\leq \calE\big(t_k^\rho,u_{k,1}^\rho,z_{k-1}^\rho\big),
		\end{align*}
		which is
		\begin{multline}\label{proof:ub2}
			\frac{\kappa}{2}\norm{z_{k,i}^\rho}_{\calZ}^2+\frac{1}{2}\int_\Omega \left(\bigl(z_{k,i}^\rho\bigr)^2+\eta\right)\bigl(\boldC \strain(u_{k,i}^\rho)\bigr):\strain(u_{k,i}^\rho)\dd{x}-\bigl\langle \ell(t_k^\rho),u_{k,i}^\rho\bigr\rangle\\
			\leq \frac{\kappa}{2}\norm{z_{k-1}^\rho}_{\calZ}^2+\frac{1}{2}\int_\Omega \left(\bigl(z_{k-1}^\rho\bigr)^2+\eta\right)\bigl(\boldC \strain(u_{k,1}^\rho)\bigr):\strain(u_{k,1}^\rho)\dd{x}-\bigl\langle \ell(t_k^\rho),u_{k,1}^\rho\bigr\rangle.
		\end{multline}
		Observe that by \eqref{proof:ub1} and the assumption on $\ell$ from \eqref{assume:external} we have for all $i\in\N$
		\begin{align*}
			\babs{\bigl\langle \ell(t_k^\rho),u_{k,i}^\rho\bigr\rangle}\overset{\text{\eqref{proof:ub1}}}{\leq} C\norm{\ell(t_k^\rho)}_{\left(\calU^{\tilde{p}'}\right)^\ast}^2\overset{\text{\eqref{assume:external}}}{\leq} C\norm{\ell}_{C^{1,1}\left(\left[0,T\right];\left(\calU^{\tilde{p}'}\right)^\ast\right)}^2.
		\end{align*}
		Moreover, by Hölder's inequality and the first estimate in \eqref{proof:uniformbounds} we obtain for all $z\in \calZ_{\left[0,1\right]}$ and all $i\in\N$
		\begin{align*}
			\abs{\int_\Omega \left(z^2+\eta\right)\bigl(\boldC \strain(u_{k,i}^\rho)\bigr):\strain(u_{k,i}^\rho)\dd{x}}\leq C\left( 1+\eta\right)\norm{u_{k,i}^\rho}_{\calU^{\tilde{p}}}\leq C.
		\end{align*}
		Therefore, since $z_{k,i}^\rho,z_{k-1}^\rho\in\calZ_{\left[0,1\right]}$ by Lemma~\ref{lemma:be_AM} and assumption \eqref{assumption:AMkrhoz0}, we can bound the left hand side of \eqref{proof:ub2} from below and its right hand side from above by using the previous two estimates. This yields
%
%
		\begin{align*}
			\exists C>0,\,\forall i\in\bbN,\,\forall \tilde{p}\in\left[2,p\right]:\quad\frac{\kappa}{2}\bnorm{z_{k,i}^\rho}_{\calZ}^2 -C\leq \frac{\kappa}{2}\norm{z_{k-1}^\rho}_{\calZ}^2 +C
		\end{align*}
		with $C$ independent of $\rho$ and $k$.
	\end{proof}

	\begin{proposition}[Convergent subsequence of the AM iterates $u_{k,i}^\rho$ and $z_{k,i}^\rho$]\label{lemma_convergentsubsequence}
	Assume \eqref{assumption:AMp}, \eqref{assumption:AMkrhoz0} and \eqref{assumption:AMiterates} are satisfied. Then there exists an index subsequence $\left(i_m\right)_{m\in\bbN}$ and limits $u_k^\rho\in\calU^{\tilde{p}}$ and $z_k^\rho\in\calZ$ such that
	\begin{align*}
		\forall j\in\N_0:\quad 
		u_{k,i_m+j}^\rho\rtom u_k^\rho\quad\text{in }\calU^{\tilde{p}}\quad\text{and}\quad
		z_{k,i_m+j}^\rho\rtom z_k^\rho\quad\text{in } \calZ.
	\end{align*}
	Moreover, if $\left(i_m\right)_{m\in\bbN}$ is an index sequence such that the above limits exist, then these limits are fixpoints of the above alternate minimization scheme, i.e.
	\begin{align}\label{fixpoint:u}
		u_k^\rho&=\argmin\Set{\calE\bigl(t_k^\rho,u,z_k^\rho\bigr)}{u\in\calU}\\
		\text{and}\quad z_k^\rho&=\argmin\Set{\calE\bigl(t_k^\rho,u_k^\rho,z\bigr)+\calR\bigl(z-z_{k-1}^\rho\bigr)}{z\in\calZ,\norm{z-z_{k-1}^\rho}_\calV\leq \rho}.\label{fixpoint:z}
	\end{align}
	Besides, it is $z_k^\rho\in\calZ_{\left[0,1\right]}$ and $z_k^\rho\leq z_{k-1}^\rho$ almost everywhere in $\Omega$.
\end{proposition}

\begin{proof}
	With the uniform bounds from the previous Lemma~\ref{lemma:boundsAM} we first obtain only weakly convergent subsequences: There exist a index subsequence $\left(i_m\right)_{m\in\N}$ and  $\tilde{\tilde{z}}_k^\rho\in\calZ$ with
	\begin{align*}
		z_{k,i_m-2}^\rho\rhupm \tilde{\tilde{z}}_k^\rho\quad \text{weakly in }\calZ
	\end{align*}
	and by possibly extracting a further subsequence (not relabeled) we may also assume that there exists $\tilde{z}_k^\rho\in\calZ$ with
	\begin{align*}
		z_{k,i_m-1}^\rho\rhupm \tilde{z}_k^\rho\quad \text{weakly in }\calZ.
	\end{align*}
	Then the compact embedding $\calZ\hookrightarrow\hookrightarrow L^r(\Omega)$ for all $r\in\left[1,\infty\right)$ gives us
	\begin{align}\label{proof:css1}
		\forall r\in\left[1,\infty\right):\quad z_{k,i_m-2}^\rho\rtom \tilde{\tilde{z}}_k^\rho\quad\text{and}\quad z_{k,i_m-1}^\rho\rtom \tilde{z}_k^\rho\quad \text{strongly in }L^r(\Omega).
	\end{align} 
	Now we utilize the continuous dependencies from Lemma~\ref{lemmaA1} and Lemma~\ref{lemmaA2EM} to get strong convergence in $\calU^{\tilde{p}}$ and $\calZ$. 
	With $u_{k,i_m-1}^\rho=u_{min}\bigl(z_{k,i_m-2}^\rho\bigr)$, Lemma~\ref{lemmaA1} implies for all $m,l\in\N$ and for $r_1=\frac{\tilde{p}p}{p-\tilde{p}}>\tilde{p}$
	\begin{align*}
		\bnorm{u_{k,i_m-1}^\rho-u_{k,i_l-1}^\rho}_{\calU^{\tilde{p}}}\leq C\bnorm{z_{k,i_m-2}^\rho-z_{k,i_l-2}^\rho}_{L^{r_1}(\Omega)}.
	\end{align*}
	Thus the first part of \eqref{proof:css1} implies that $\bigl(u_{k,i_m-1}^\rho\bigr)_{m\in\N}$ is a Cauchy sequence in $\calU^{\tilde{p}}$. Since $\calU^{\tilde{p}}$ is complete it therefore has a limit $\tilde{u}_k^\rho\in \calU^{\tilde{p}}$, i.e.
	\begin{align}\label{proof:css2}
		u_{k,i_m-1}^\rho\rtom \tilde{u}_k^\rho\quad \text{strongly in }\calU^{\tilde{p}}.
	\end{align} 
	With $z_{k,i_m-1}^\rho=z_{min}\bigl(u_{k,i_m-1}^\rho\bigr)$, Lemma~\ref{lemmaA2EM} implies for all $m,l\in\N$ and for $r_2=\frac{\tilde{p}}{\tilde{p}-2}$
	\begin{multline}
		\bnorm{z_{k,i_m-1}^\rho-z_{k,i_l-1}^\rho}_{\calZ}^2\leq\\ C\bnorm{u_{k,i_m-1}^\rho-u_{k,i_l-1}^\rho}_{\calU^{\tilde{p}}}\left(\bnorm{u_{k,i_m-1}^\rho}_{\calU^{\tilde{p}}}+\bnorm{u_{k,i_l-1}^\rho}_{\calU^{\tilde{p}}}\right)\bnorm{z_{k,i_m-1}^\rho-z_{k,i_l-1}^\rho}_{L^{r_2}(\Omega)}.
	\end{multline}
	Therefore the uniform bound from \eqref{proof:uniformbounds} and the strong convergences \eqref{proof:css1} and \eqref{proof:css2} imply that $\bigl(z_{k,i_m-1}^\rho\bigr)_{m\in\N}$ is a Cauchy sequence in $\calZ$. Since $\calZ$ is complete we thus obtain
	\begin{align*}
		z_{k,i_m-1}^\rho\rtom \tilde{z}_k^\rho\quad \text{even strongly in }\calZ.
	\end{align*}
	By using again the continuous dependencies from Lemma~\ref{lemmaA1} and Lemma~\ref{lemmaA2EM} and repeating the above arguments with $i_m$ instead of $(i_m-1)$ we find $u_k^\rho\in\calU^{\tilde{p}}$ and $z_k^\rho\in\calZ$ such that
	\begin{align*}
		u_{k,i_m}^\rho\rtom u_k^\rho\quad\text{in }\calU^{\tilde{p}}\quad\text{and}\quad
		z_{k,i_m}^\rho\rtom z_k^\rho\quad\text{in } \calZ.
	\end{align*}
	Observe that $z_k^\rho\in \calZ_{\left[0,1\right]}$ and $z_{k}^\rho\leq z_{k-1}^\rho$ almost everywhere in $\Omega$.
	Now we will show that the limits are fixpoints of the alternate minimization scheme.
	By Lemma~\ref{lemma:continuity}  $\calE\big(t_k^\rho,\cdot,\cdot\big)+\calR\big(\cdot-z_{k-1}^\rho\big)$ is lower semicontinuous on $\calU^{\tilde{p}}\times\calZ$ \textcircled{\small 1} and $\calE\bigl(t_k^\rho,\cdot,z_{min}\bigl(u_k^\rho\bigr)\bigr)$ is continuous on $\calU^{\tilde{p}}$ \textcircled{\small 2}. Together with the previous convergences and the minimality of $z_{k,i_m}^\rho=z_{min}\bigl(u_{k,i_m}^\rho\bigr)$ \textcircled{\small 3} this implies
	\begin{align*}
		\calE\big(t_k^\rho,u_k^\rho,z_k^\rho\big)+\calR\big(z_k^\rho-z_{k-1}^\rho\big)
		&\overset{\textcircled{\small 1}}{\leq} \liminf_{m\to\infty} \calE\Big(t_k^\rho,u_{k,i_m}^\rho,z_{k,i_m}^\rho\Big)+\calR\Big(z_{k,i_m}^\rho-z_{k-1}^\rho\Big)\\
		&\overset{\textcircled{\small 3}}{\leq} \liminf_{m\to\infty}\calE\Big(t_k^\rho,u_{k,i_m}^\rho,z_{min}\bigl(u_k^\rho\bigr)\Big)+\calR\big(z_{min}\bigl(u_k^\rho\bigr)-z_{k-1}^\rho\big)\\
		&\overset{\textcircled{\small 2}}{=}\calE\big(t_k^\rho,u_k^\rho,z_{min}\bigl(u_k^\rho\bigr)\big)+\calR\big(z_{min}\bigl(u_k^\rho\bigr)-z_{k-1}^\rho\big).
	\end{align*}
	So by uniqueness of the minimizer $z_{min}\bigl(u_k^\rho\bigr)$ we have $z_k^\rho=z_{min}\bigl(u_k^\rho\bigr)$ and \eqref{fixpoint:z} follows for $j=0$.
	Lemma~\ref{lemma:be_AM}, $z_{k,i_m-1}^\rho\rtom\tilde{z}_k^\rho$ in $\calZ$ and Lemma~\ref{lemma:closedness} imply that $\tilde{z}_k^\rho\in\left[0,1\right]$ a.e. in $\Omega$. Therefore we obtain as in Lemma~\ref{lemma:boundsAM} with Theorem~\ref{thm:improved integrability} that $u_{min}(\tilde{z}_k^\rho)\in\calU^{\tilde{p}}$. Thus we know from Lemma~\ref{lemma:continuity} that $\calE(t_k^\rho,u_{min}(\tilde{z}_k^\rho),\cdot)$ is continuous on $\calZ$ and that $\calE\big(t_k^\rho,\cdot,\cdot\big)$ is lower semicontinuous on $\calU^{\tilde{p}}\times\calZ$. Together with the convergences $u_{k,i_m}^\rho\rtom u_k^\rho$ in $\calU^{\tilde{p}}$ and $z_{k,i_m-1}^\rho\rtom \tilde{z}_k^\rho$ in $\calZ$ as well as the minimality of $u_{k,i_m}^\rho=u_{min}\bigl(z_{k,i_m-1}^\rho\bigr)$ this leads to
	\begin{align*}
		\calE\big(t_k^\rho,u_k^\rho,\tilde{z}_k^\rho\big)&\leq \liminf_{m\to\infty} \calE\Big(t_k^\rho,u_{k,i_m}^\rho,z_{k,i_m-1}^\rho\Big)\\
		&\leq \lim\limits_{m\to\infty}\calE\Big(t_k^\rho,u_{min}\bigl(\tilde{z}_{k}^\rho\bigr),z_{k,i_m-1}^\rho\Big)=\calE\big(t_k^\rho,u_{min}\bigl(\tilde{z}_{k}^\rho\bigr),\tilde{z}_k^\rho\big).
	\end{align*}
	So by uniqueness of the minimizer $u_{min}\bigl(\tilde{z}_{k}^\rho\bigr)$ we have $u_k^\rho=u_{min}\bigl(\tilde{z}_{k}^\rho\bigr)$.
	Finally we need to show that $z_k^\rho=\tilde{z}_k^\rho$. For this purpose we employ the following index trick from \cite[Lemma~2.10; Prop.~3.2]{Almi2020}. By Lemma~\ref{lemma:monotonicity} the sequence $\left(\calE\big(t_k^\rho,u_{k,i}^\rho,z_{k,i}^\rho\big)+\calR\big(z_{k,i}^\rho-z_{k-1}^\rho\big)\right)_{i\in\N_0}$ is monotonically decreasing \textcircled{\small 4}. Together with $i_{(m-1)}\leq i_m-1$ for all $m\in\N$ \textcircled{\small 5} and the minimality of $u_{k,i_m}^\rho=u_{min}\bigl(z_{k,i_m-1}^\rho\bigr)$ \textcircled{\small 6} and $z_{k,i_{(m-1)}}^\rho=z_{min}\Bigl(u_{k,i_{(m-1)}}^\rho\Bigr)$ \textcircled{\small 7} we thus obtain
	\begin{align}\notag
		\calE\Big(t_k^\rho,u_{k,i_m}^\rho,z_{k,i_m-1}^\rho\Big)+\calR\Big(z_{k,i_m-1}^\rho-z_{k-1}^\rho\Big)
		&\overset{\textcircled{\small 6}}{\leq} \calE\Big(t_k^\rho,u_{k,i_m-1}^\rho,z_{k,i_m-1}^\rho\Big)+\calR\Big(z_{k,i_m-1}^\rho-z_{k-1}^\rho\Big)\\\notag
		&\hspace{-8.5pt}\overset{\textcircled{\small 4}\&\textcircled{\small 5}}{\leq} \calE\Big(t_k^\rho,u_{k,i_{(m-1)}}^\rho,z_{k,i_{(m-1)}}^\rho\Big)+\calR\Big(z_{k,i_{(m-1)}}^\rho-z_{k-1}^\rho\Big)\\\label{proof:css3}
		&\overset{\textcircled{\small 7}}{\leq} \calE\Big(t_k^\rho,u_{k,i_{(m-1)}}^\rho,z_k^\rho\Big)+\calR\big(z_k^\rho-z_{k-1}^\rho\big).
	\end{align}
	Now again by $u_{k,i_m}^\rho\rtom u_k^\rho$ in $\calU^{\tilde{p}}$, $z_{k,i_m-1}^\rho\rtom \tilde{z}_k^\rho$ in $\calZ$, the lower semicontinuity of $\calE\big(t_k^\rho,\cdot,\cdot\big)+\calR\big(\cdot-z_{k-1}^\rho\big)$ on $\calU^{\tilde{p}}\times\calZ$ (Lemma~\ref{lemma:continuity}) \textcircled{\small 8} and the continuity of $\calE(t_k^\rho,\cdot,z_k^\rho)$ on $\calU^{\tilde{p}}$ \textcircled{\small 9} this implies
	\begin{align*}
		\calE\big(t_k^\rho,u_k^\rho,\tilde{z}_k^\rho\big)+\calR\big(\tilde{z}_k^\rho-z_{k-1}^\rho\big)
		&\overset{\textcircled{\small 8}}{\leq}\liminf_{m\to\infty} \calE\Big(t_k^\rho,u_{k,i_m}^\rho,z_{k,i_m-1}^\rho\Big)+\calR\Big(z_{k,i_m-1}^\rho-z_{k-1}^\rho\Big)\\
		&\hspace{-5pt}\overset{\eqref{proof:css3}}{\leq} \liminf_{m\to\infty}\calE\Big(t_k^\rho,u_{k,i_{(m-1)}}^\rho,z_k^\rho\Big)+\calR\big(z_k^\rho-z_{k-1}^\rho\big)\\
		&\overset{\textcircled{\small 9}}{=}\calE\big(t_k^\rho,u_k^\rho,z_k^\rho\big)+\calR\big(z_k^\rho-z_{k-1}^\rho\big).
	\end{align*}
	So by uniqueness of the minimizer $z_k^\rho=z_{min}\bigl(u_k^\rho\bigr)$ we obtain $z_k^\rho=\tilde{z}_k^\rho$.
	In summary we proved
	\begin{align*}
		u_{k,i_m}^\rho\rtom u_k^\rho\quad\text{in }\calU^{\tilde{p}}\quad\text{and}\quad
		z_{k,i_m}^\rho\rtom z_k^\rho\quad\text{in } \calZ
	\end{align*}
	with $z_k^\rho=z_{min}\bigl(u_k^\rho\bigr)$ and $u_k^\rho=u_{min}\bigl(z_k^\rho\bigr)$, i.e. the assertion for $j=0$. 
	
	By induction we can now prove that the assertion is actually true for all $j\in\N_0$. Assume the assertion holds true for an arbitrarily fixed $j\in\N_0$. Then by repeating the above arguments with $(i_m+j+1)$ instead of $i_m$ we find $\hat{u}_k^\rho\in\calU^{\tilde{p}}$ and $\hat{z}_k^\rho\in\calZ$ such that
	\begin{align*}
		u_{k,i_m+j+1}^\rho\rtom \hat{u}_k^\rho\quad\text{in }\calU^{\tilde{p}}\quad\text{and}\quad
		z_{k,i_m+j+1}^\rho\rtom \hat{z}_k^\rho\quad\text{in } \calZ
	\end{align*}
	and also $\hat{z}_k^\rho=z_{min}\bigl(\hat{u}_k^\rho\bigr)$ and $\hat{u}_k^\rho=u_{min}\bigl(z_k^\rho\bigr)$ (notice that $z_k^\rho$ now plays the role of $\tilde{z}_k^\rho$). So $u_k^\rho=u_{min}\bigl(z_k^\rho\bigr)$ implies $\hat{u}_k^\rho=u_k^\rho$ and therefore $\hat{z}_k^\rho=z_{min}\bigl(\hat{u}_k^\rho\bigr)=z_{min}\bigl(u_k^\rho\bigr)=z_k^\rho$, what finishes the proof.
\end{proof}

Notice that the convergences are not necessarily uniform in $j\in\N_0$, i.e. for $\varepsilon>0$ there does not necessarily exists an index $M_\varepsilon\in\N$ independent of $j$ such that all subsequent elements of the sequences stay in the $\varepsilon$-neighborhood of their limits.
	
	\begin{lemma}[Accumulation points of the AM iterates]\label{lemma:accumulation points}
		Assume \eqref{assumption:AMp}, \eqref{assumption:AMkrhoz0} and \eqref{assumption:AMiterates} are satisfied.
		Either the whole sequence $\bigl(u_{k,i}^\rho,z_{k,i}^\rho\bigr)_{i\in\bbN}$ converges in $\calU^{\tilde{p}}\times\calZ$ or the set of its accumulation points has no isolated point with respect to $\calU^{\tilde{p}}\times\calZ$. 
	\end{lemma}
	\begin{proof}
		Suppose there are two accumulation points $\big(u_k^\rho,z_k^\rho\big)\neq \big(\tilde{u}_k^\rho,\tilde{z}_k^\rho\big)$ of the sequence $\bigl(u_{k,i}^\rho,z_{k,i}^\rho\bigr)_{i\in\bbN}$, where $\big(u_k^\rho,z_k^\rho\big)$ is an isolated point in the set of all accumulation points $A\subset \calU^{\tilde{p}}\times\calZ$, i.e. 
		\begin{align}\label{fixpoint_contradiction}
			\exists \delta>0,\,\forall (u,z)\in A\setminus\left\{ \big(u_k^\rho,z_k^\rho\big) \right\}:\quad \norm{\begin{pmatrix}
					u_k^\rho\\
					z_k^\rho
				\end{pmatrix}-\begin{pmatrix}
					u\\
					z
			\end{pmatrix}}_{\calU^{\tilde{p}}\times\calZ}=\bnorm{u_k^\rho-u}_{\calU^{\tilde{p}}}+\bnorm{z_k^\rho-z}_\calZ>\delta.
		\end{align}
		Now let $\varepsilon>0$ be sufficiently small such that $\big(C_1+C_1 C_2\big)\varepsilon<\delta$, where according to Lemma~\ref{lemmaA1}
		\begin{align*}
			\exists C_1>0,\,\forall z_1,z_2\in\calZ_{\left[0,1\right]}: \quad \bnorm{u_{min}(z_2)-u_{min}(z_1)}_{\calU^{\tilde{p}}}\leq C_1\bnorm{z_2-z_1}_\calZ
		\end{align*}
		and according to Lemma~\ref{lemmaA2EM}
		\begin{align*}
			\exists C_2>0,\,\forall u_1,u_2\in\calU^{\tilde{p}}:\quad \bnorm{z_{min}(u_2)-z_{min}(u_1)}_\calZ\leq C_2\bnorm{u_2-u_1}_{\calU^{\tilde{p}}}.
		\end{align*}
		Especially let $\varepsilon<\delta$. For $r>0$ let $B_r\bigl(u_k^\rho,z_k^\rho\bigr)$ denote the open ball with radius $r$ with respect to $\calU^{\tilde{p}}\times\calZ$ centered at $\bigl(u_k^\rho,z_k^\rho\bigr)$. Now since $\bigl(u_k^\rho,z_k^\rho\bigr)$ and $\bigl(\tilde{u}_k^\rho,\tilde{z}_k^\rho\bigr)$ are both accumulation points and since $\bigl(\tilde{u}_k^\rho,\tilde{z}_k^\rho\bigr)$ lies outside of $B_\delta\bigl(u_k^\rho,z_k^\rho\bigr)$, there has to be an infinite set of indices $I$, such that for all $i\in I$ we have
		\begin{align*}
			\bigl(u_{k,i}^\rho,z_{k,i}^\rho\bigr)\in B_\varepsilon\bigl(u_k^\rho,z_k^\rho\bigr)\quad\text{and}\quad \bigl(u_{k,i+1}^\rho,z_{k,i+1}^\rho\bigr)\notin B_\varepsilon\bigl(u_k^\rho,z_k^\rho\bigr).
		\end{align*}
		With Lemma~\ref{lemma:be_AM}, Lemma~\ref{lemma:closedness}, $u_{k,i+1}^\rho=u_{min}(z_{k,i}^\rho)$ and $u_k^\rho=u_{min}(z_k^\rho)$ (Prop.~\ref{lemma_convergentsubsequence}) it follows by Lemma~\ref{lemmaA1}
		\begin{align*}
			\bnorm{u_{k,i+1}^\rho-u_k^\rho}_{\calU^{\tilde{p}}}\leq C_1\bnorm{z_{k,i}^\rho-z_k^\rho}_\calZ\leq C_1\varepsilon.
		\end{align*}
		This implies together with Lemma~\ref{lemmaA2EM} noticing $z_{k,i+1}^\rho=z_{min}(u_{k,i+1}^\rho)$ and $z_k^\rho=z_{min}(u_k^\rho)$ (Prop.~\ref{lemma_convergentsubsequence})
		\begin{align*}
			\bnorm{z_{k,i+1}^\rho-z_k^\rho}_\calZ\leq C_2\bnorm{u_{k,i+1}^\rho-u_k^\rho}_{\calU^{\tilde{p}}}\leq C_1C_2\varepsilon.
		\end{align*}
		In summary we obtain
		\begin{align*}
			\norm{\begin{pmatrix}
					u_{k,i+1}^\rho\\
					z_{k,i+1}^\rho
				\end{pmatrix}-\begin{pmatrix}
					u_k^\rho\\
					z_k^\rho
			\end{pmatrix}}_{\calU^{\tilde{p}}\times\calZ}=\bnorm{u_{k,i+1}^\rho-u_k^\rho}_{\calU^{\tilde{p}}}+\bnorm{z_{k,i+1}^\rho-z_k^\rho}_\calZ\leq \left(C_1+C_1C_2\right)\varepsilon<\delta.
		\end{align*}
		Hence it holds
		\begin{align*}
			\bigl(u_{k,i+1}^\rho,z_{k,i+1}^\rho\bigr)\in B_\delta\big(u_k^\rho,z_k^\rho\big)\big\backslash B_\varepsilon\big(u_k^\rho,z_k^\rho\big).
		\end{align*}
		But as in Prop.~\ref{lemma_convergentsubsequence} we can find a convergent subsequence of $\bigl(u_{k,i+1}^\rho,z_{k,i+1}^\rho\bigr)_{i\in I}$ with respect to $\calU^{\tilde{p}}\times\calZ$. And since $\overline{B_\delta\big(u_k^\rho,z_k^\rho\big)}^{\calU^{\tilde{p}}\times\calZ}\Big\backslash B_\varepsilon\big(u_k^\rho,z_k^\rho\big)$ is closed in $\calU^{\tilde{p}}\times\calZ$, we have for its limit $\Big(\tilde{\tilde{u}}_k^\rho,\tilde{\tilde{z}}_k^\rho\Big)$
		\begin{align*}
			\varepsilon\leq 
			\norm{\begin{pmatrix}
					u_{k}^\rho\\
					z_{k}^\rho
				\end{pmatrix}-\begin{pmatrix}
					\tilde{\tilde{u}}_k^\rho\\
					\tilde{\tilde{z}}_k^\rho
			\end{pmatrix}}_{\calU^{\tilde{p}}\times\calZ}\leq\delta
		\end{align*}
		what contradicts \eqref{fixpoint_contradiction}.
	\end{proof}

	However, within an implementation we are not able to select convergent subsequences respectively their limits.
	Fortunately, the fixpoint property of the E\&M iterates $u_k^\rho$ and $z_k^\rho$ we proved in Proposition~\ref{lemma_convergentsubsequence} is sufficient to carry out our convergence analysis. 
	The fixpont set is an attractive set in the sense formulated in the next lemma.
	\begin{lemma}[Attractive fixpoint set]
		Assume \eqref{assumption:AMp}, \eqref{assumption:AMkrhoz0} and \eqref{assumption:AMiterates} are satisfied.
		The fixpoint set $A_k^\rho\subset \calU^{\tilde{p}}\times\calZ_{\left[0,1\right]}$ of the AM scheme is \textit{attractive} in the following sense
		\begin{enumerate}
			\item it is obviously \textit{forward invariant under the AM scheme}, i.e. if we start an AM step with an initial value $z_{k,0}^\rho\in\calZ_{\left[0,1\right]}$ such that $\bigl(u_{k,1}^\rho,z_{k,0}^\rho\bigr)\in A_k^\rho$, then $\bigl(u_{k,i}^\rho,z_{k,i}^\rho\bigr)\in A_k^\rho$ for all $i\in\N$,
			\item the \textit{basin of attraction} for $A_k^\rho$ given by\label{basinofattraction}
			\begin{align*}
				B\bigl(A_k^\rho\bigr)\coloneq \Set{\substack{\text{initial value }z_{k,0}^\rho\in\calZ_{\left[0,1\right]}\\\text{ for an AM step}}}{\substack{\text{for any open neighborhood } N \text{ of } A_k^\rho,\\\exists i_0\in\N,\,\forall i\geq i_0:\quad \left(u_{k,i}^\rho,z_{k,i}^\rho\right)\in N}}
			\end{align*}
			is equal to the whole set $\calZ_{\left[0,1\right]}$, i.e. no matter with what initial values we start, the iterates converge to the fixpoint set.
		\end{enumerate}
	\end{lemma}
	Property \eqref{basinofattraction} can be proven analogously to the usual subsubsequence argument in topological spaces, that a sequence converges if any subsequence has a convergent subsubsequence and all of them share the same limit.
	
	\begin{proof}[Proof of \eqref{basinofattraction}]
		Let $z_{k,0}^\rho\in\calZ_{\left[0,1\right]}$ be an arbitrary initial value for an AM step and let $N\subset \calU^{\tilde{p}}\times\calZ$ be an arbitrary open neighborhood of $A_k^\rho$. Assume there exists a subsequence $\left(i_m\right)_{m\in\N}$ such that $\bigl(u_{k,i_m}^\rho,z_{k,i_m}^\rho\bigr)\notin N$ for all $m\in\N$. Then we find as in Proposition~\ref{lemma_convergentsubsequence} a subsubsequence $\left(i_{m_l}\right)_{l\in\N}$, such that
		\begin{align*}
			\Bigl(u_{k,i_{m_l}}^\rho,z_{k,i_{m_l}}^\rho\Bigr)\rtol \bigl(u_k^\rho,z_k^\rho\bigr)\in A_k^\rho\text{ in } \calU^{\tilde{p}}\times\calZ.
		\end{align*}
		Since $\bigl(u_k^\rho,z_k^\rho\bigr)\in A_k^\rho$ and $N$ is an open neighborhood of $A_k^\rho$, we therefore obtain
		\begin{align*}
			\exists L\in\N,\,\forall l\geq L:\quad \left(u_{k,i_{m_l}}^\rho,z_{k,i_{m_l}}^\rho\right)\in N
		\end{align*}
		in contradiction to $\left(u_{k,i_m}^\rho,z_{k,i_m}^\rho\right)\notin N$ for all $m\in\N$. Therefore our assumption was wrong and we have
		\begin{align*}
			\exists i_0\in\N,\,\forall i\geq i_0:\quad \bigl(u_{k,i}^\rho,z_{k,i}^\rho\bigr)\in N.
		\end{align*}
		Since $N$ and $z_{k,0}^\rho\in\calZ_{\left[0,1\right]}$ were arbitrary this shows $B\bigl(A_k^\rho\bigr)=\calZ_{\left[0,1\right]}$.
	\end{proof}
		\begin{remark}
		Observe that by applying Proposition~\ref{lemma_convergentsubsequence} to weakly convergent subsequences of the AM sequence $\bigl(u_{k,i}^\rho,z_{k,i}^\rho\bigr)_{i\in\bbN}$, we find that every accumulation point of the AM sequence in the weak topology of $\calU^{\tilde{p}}\times\calZ$ is even an accumulation point in the strong topology of $\calU^{\tilde{p}}\times\calZ$ and particularly it is a fixpoint of the AM scheme.
	\end{remark}

	\subsection{Properties of and uniform bounds for the Efendiev\&Mielke iterates
	}
	Within the rest of this section we will have the following\\
	\noindent
	\textbf{Assumptions. } 
	\begin{align}\label{assumption:initialvalue}
		\text{Let an initial value } z_0\in\calZ_{\left[0,1\right]} \text{ be given.}
	\end{align}
	\begin{align}
		\begin{split}
			&\text{Further let } \bigl(u_{k}^\rho\bigr)_{k\in\bbN}\subset \calU \text{ and } \bigl(z_{k}^\rho\bigr)_{k\in\bbN}\subset\calZ \text{ be the sequences generated by the outer for-}\\\label{assumption:EMiterates}
			&\text{loop of Algorithm~\ref{alg1} for a given } \rho>0,\text{ i.e. let $u_k^\rho$ and $z_k^\rho$ be limits as in Proposition~\ref{lemma_convergentsubsequence}}.
		\end{split}
	\end{align}
	First we observe via induction that assumption \eqref{assumption:AMkrhoz0} from the previous section is satisfied, if the initial value $z_0$ lies in the physically relevant range $\calZ_{\left[0,1\right]}$. 
	\begin{corollary}[Basic estimates for the E\&M iterates $z_k^\rho$]\label{cor:beEM}
		Assume \eqref{assumption:initialvalue} and \eqref{assumption:EMiterates}.
		Then for all $\rho>0$ and all $k\in\bbN_0$ we have
		\begin{align}\label{prop_be1}
			z_k^\rho\in\calZ_{\left[0,1\right]}.
		\end{align}
	\end{corollary}
	\begin{proof}
		This follows via induction over $k\in\N_0$ from Lemma~\ref{lemma:be_AM}, Proposition~\ref{lemma_convergentsubsequence} and Lemma~\ref{lemma:closedness}.
	\end{proof}
	To derive the following uniform bounds we follow \cite[Prop.~2.1]{Knees2018}.
	\begin{lemma}[Uniform bounds for the E\&M iterates $u_k^\rho$ and $z_k^\rho$]\label{prop:basicestimates}
		Assume \eqref{assumption:initialvalue} and \eqref{assumption:EMiterates}.
		Then we have
		\begin{align}\label{prop_be2}
			\sup\limits_{\rho>0,\,k\in\bbN_0} \bnorm{u_k^\rho}_{\calU^{\tilde{p}}} <\infty
			\quad \text{and} \quad 
			\sup\limits_{\rho>0,\,k\in\bbN_0} \bnorm{z_k^\rho}_\calZ <\infty.
		\end{align}
		Moreover, for $u_{min}(t_0,z_0)\coloneq\argmin\Set{\calE\big(t_0,u,z_0\big)}{u\in\calU}$ it holds
			\begin{align}\label{upperboundenergyEM}
				\forall \rho>0:\quad \calE\big(t_0^\rho,u_{0}^\rho,z_{0}^\rho\big)\leq \calE\big(t_0,u_{min}(t_0,z_0),z_{0}\big).
			\end{align}
	\end{lemma}	
	\begin{proof}
		By Corollary~\ref{cor:beEM} assumption \eqref{assumption:AMkrhoz0} is satisfied. Therefore, Lemma~\ref{lemma:boundsAM} and Proposition~\ref{lemma_convergentsubsequence} imply the uniform bound for the displacement iterates $u_k^\rho$ in \eqref{prop_be2}.
		
		To prove the uniform bound for the damage iterates $z_k^\rho$ we will show that there exists $C>0$ such that $\calE\big(t_k^\rho,u_k^\rho,z_k^\rho\big)\leq C$ for all $k\in\bbN$ and then use the coercivity estimate \eqref{coercivity}.		
		With \eqref{coercivity} we find $\mu, c>0$ such that for all $t\in\left[0,T\right]$, $u\in\calU$ and $z\in\calZ$
		\begin{align}\label{proof:prop_be5}
			\abs{\partial_t\calE(t,u,z)}=\big\vert\big\langle \dot{\ell}(t),u\big\rangle\big\vert\leq \bigl\Vert\dot{\ell}(t)\bigr\Vert_{\calU^*}\norm{u}_\calU\leq c_\ell\norm{u}_\calU\leq \mu\Bigl(\calE(t,u,z)+c\Bigr).
		\end{align}
		Hence we have for all $t\in\left[0,T\right]$, $u\in\calU$ and $z\in\calZ$ and for $s\in\left[0,t\right]$
		\begin{align*}
			\calE(t,u,z)+c=\int_{s}^{t}\partial_t\calE(\tau,u,z)\dd{\tau}+\calE(s,u,z)+c\overset{\text{\eqref{proof:prop_be5}}}{\leq} \int_{s}^{t}\mu\Big(\calE(\tau,u,z)+c\Big)\dd{\tau}+\calE(s,u,z)+c.
		\end{align*}
		Now analogously to \cite[Section~2.1.1]{MielkeRoubicek2015} we obtain for all $t\in\left[0,T\right]$, $u\in\calU$ and $z\in\calZ$ and for $s\in\left[0,t\right]$ with Grönwalls inequality
		\begin{align*}
			\calE(t,u,z)+c\leq \calE(s,u,z)+c+ \int_{s}^{t}\Bigl(\calE(s,u,z)+c\Bigr)\mu e^{\mu\left(t-\tau\right)}\dd{\tau}=\Big(\calE(s,u,z)+c\Big) e^{\mu\left(t-s\right)}
		\end{align*}
		and again with \eqref{proof:prop_be5} this implies
		\begin{align}\label{proof:prop_be6}
			\abs{\partial_t\calE(t,u,z)}\leq \mu\Big(\calE(t,u,z)+c\Big)\leq \mu\Bigl(\calE(s,u,z)+c\Bigr) e^{\mu\left(t-s\right)}.
		\end{align}
		By \eqref{upperboundenergyAM} we have for all $\rho>0$, $k\in\N_0$ and $i\in\N$
		\begin{align*}
			\calE\big(t_k^\rho,u_{k,i}^\rho,z_{k,i}^\rho\big)\leq \calE\big(t_k^\rho,u_{k,1}^\rho,z_{k-1}^\rho\big),
		\end{align*}
		what together with $u_{k,1}^\rho=\argmin\Set{\calE\bigl(t_k^\rho,u,z_{k,0}^\rho\bigr)}{u\in\calU}$, $z_{k,0}^\rho=z_{k-1}^\rho$, $t_0^\rho=t_0$ and $z_{-1}^\rho=z_0$ leads to
		\begin{align*}
			\calE\bigl(t_k^\rho,u_{k,i}^\rho,z_{k,i}^\rho\bigr)\leq \calE\bigl(t_k^\rho,u_{k-1}^\rho,z_{k-1}^\rho\bigr)\quad\text{for } k\in\N\quad
			\text{and}\quad\calE\bigl(t_0^\rho,u_{0,i}^\rho,z_{0,i}^\rho\bigr)\leq \calE\bigl(t_0,u_{0,1}^\rho,z_{0}\bigr)\quad\text{for } k=0.
		\end{align*}
		Together with the convergences $u_{k,i_m}^\rho\rtom u_k^\rho$ in $\calU^{\tilde{p}}$ and $z_{k,i_m}^\rho\rtom z_k^\rho$ in $\calZ$ from Lemma~\ref{lemma_convergentsubsequence} and the continuity of $\calE\bigl(t_k^\rho,\cdot,\cdot\bigr)$ on $\calU^{\tilde{p}}\times\calZ$ from Lemma~\ref{lemma:continuity} this yields
		\begin{align}\notag
			\calE\bigl(t_k^\rho,u_k^\rho,z_k^\rho\bigr)&= \lim_{m\to\infty} \calE\bigl(t_k^\rho,u_{k,i_m}^\rho,z_{k,i_m}^\rho\bigr)\leq \calE\bigl(t_k^\rho,u_{k-1}^\rho,z_{k-1}^\rho\bigr)\quad \text{for }k\in\N\\\label{initialbounds}
			\text{and}\quad \calE\bigl(t_0^\rho,u_0^\rho,z_0^\rho\bigr)&= \lim_{m\to\infty} \calE\bigl(t_0^\rho,u_{0,i_m}^\rho,z_{0,i_m}^\rho\bigr)\leq \calE\bigl(t_0,u_{0,1}^\rho,z_{0}\bigr)\quad \text{for }k=0.
		\end{align}
		Therefore we have for $k\in\N$
		\begin{align*}
			\calE\bigl(t_k^\rho,u_k^\rho,z_k^\rho\bigr)- \calE\bigl(t_{k-1}^\rho,u_{k-1}^\rho,z_{k-1}^\rho\bigr)&\leq \calE\bigl(t_k^\rho,u_{k-1}^\rho,z_{k-1}^\rho\bigr)-\calE\bigl(t_{k-1}^\rho,u_{k-1}^\rho,z_{k-1}^\rho\bigr)\\
			&=\int_{t_{k-1}^\rho}^{t_k^\rho}\partial_t \calE\bigl(\tau,u_{k-1}^\rho,z_{k-1}^\rho\bigr)\dd{\tau},
		\end{align*}
		so \eqref{proof:prop_be6} with $s=t_{k-1}^\rho$ implies
		\begin{align*}
			\calE\bigl(t_k^\rho,u_k^\rho,z_k^\rho\bigr)&\leq \calE\bigl(t_{k-1}^\rho,u_{k-1}^\rho,z_{k-1}^\rho\bigr)+\int_{t_{k-1}^\rho}^{t_k^\rho} \mu \Bigl(\calE\bigl(t_{k-1}^\rho,u_{k-1}^\rho,z_{k-1}^\rho\bigr)+c\Bigr) e^{\mu\left(\tau-t_{k-1}^\rho\right)}\dd{\tau}\\
			&=\calE\bigl(t_{k-1}^\rho,u_{k-1}^\rho,z_{k-1}^\rho\bigr)+\Bigl( \calE\bigl(t_{k-1}^\rho,u_{k-1}^\rho,z_{k-1}^\rho\bigr)+c \Bigr)\left(e^{\mu\left(t_k^\rho-t_{k-1}^\rho\right)}-1\right)\\
			&=\Bigl(\calE\bigl(t_{k-1}^\rho,u_{k-1}^\rho,z_{k-1}^\rho\bigr)+c\Bigr)e^{\mu\left(t_k^\rho-t_{k-1}^\rho\right)}-c.
		\end{align*}
		Induction with respect to $k\in\bbN$ leads to
		\begin{align*}
			\calE\bigl(t_k^\rho,u_k^\rho,z_k^\rho\bigr)+c&\leq \Bigl(\calE\bigl(t_0^\rho,u_0^\rho,z_0^\rho\bigr)+c\Bigr)\prod_{j=1}^{k} e^{\mu\left(t_j^\rho-t_{j-1}^\rho\right)}=\Bigl(\calE\bigl(t_0^\rho,u_0^\rho,z_0^\rho\bigr)+c\Bigr)e^{\mu t_k^\rho}.
		\end{align*}
		So by \eqref{initialbounds} and if we take into account that $u_{0,1}^\rho=\argmin\Set{\calE\bigl(t_0,u,z_0\bigr)}{u\in\calU}$ is independent of $\rho$, we end up with the uniform estimate
		\begin{align*}
			\exists C>0,\,\forall k\in\bbN_0,\,\forall \rho>0:\quad \calE\bigl(t_k^\rho,u_k^\rho,z_k^\rho\bigr)\leq \Bigl(\calE\big(t_0,u_{0,1}^\rho,z_{0}\big)+c\Bigr)e^{\mu T}\leq C.
		\end{align*}
		Thus the uniform bound for the damage iterates $z_k^\rho$ in \eqref{prop_be2} follows from coercivity estimate \eqref{coercivity}.
	\end{proof}
	For the proof of the following Lemma we follow \cite[Proposition~2.2]{Knees2018} and we will utilize some tools from convex analysis collected in Appendix~\ref{appendix:convexanalysis}.
	\begin{lemma}[Optimality properties of the E\&M iterates $t_k^\rho$, $u_k^\rho$ and $z_k^\rho$]\label{prop:optimality}
		Assume \eqref{assumption:initialvalue} and \eqref{assumption:EMiterates} and let $k\in\bbN_0$ and $\rho>0$ be arbitrary. The E\&M iterates $\big(t_k^\rho,u_k^\rho,z_k^\rho\big)$ fulfill the following properties. There exists $\xi_k^\rho\in\calV^*$ such that
		\begin{gather}
			\big\Vert\xi_k^\rho\big\Vert_{\calV^*}\left(\norm{z_k^\rho-z_{k-1}^\rho}_\calV-\rho\right)=0\label{op40},
		\end{gather}
		\begin{align}
			\begin{split}
				\rho\dist_{\calV^\ast}\Bigl(-\D_z\calE\bigl(t_k^\rho,u_k^\rho,z_k^\rho\bigr),\partial^\calZ\calR(0)\Bigr)&=\bigl\langle \xi_k^\rho,z_k^\rho-z_{k-1}^\rho\bigr\rangle\\&=\rho\,\bigl\Vert\xi_k^\rho\bigr\Vert_{\calV^*}=
				\bigl\Vert\xi_k^\rho\bigr\Vert_{\calV^*}\norm{z_k^\rho-z_{k-1}^\rho}_\calV\label{op41},
			\end{split}
		\end{align}
		\begin{gather}
			\calR\bigl(z_k^\rho-z_{k-1}^\rho\bigr)+\rho\dist_{\calV^\ast}\Bigl(-\D_z\calE\big(t_k^\rho,u_k^\rho,z_k^\rho\big),\partial^\calZ\calR(0)\Bigr)=\left\langle -\D_z\calE\bigl(t_k^\rho,u_k^\rho,z_k^\rho\bigr),z_k^\rho-z_{k-1}^\rho\right\rangle \label{op42},
		\end{gather}
		\begin{multline}
			\calR\bigl(z_k^\rho-z_{k-1}^\rho\bigr)+\norm{z_k^\rho-z_{k-1}^\rho}_\calV\dist_{\calV^\ast}\Bigl(-\D_z\calE\bigl(t_k^\rho,u_k^\rho,z_k^\rho\bigr),\partial^\calZ\calR(0)\Bigr)\\
			=\left\langle -\D_z\calE\bigl(t_k^\rho,u_k^\rho,z_k^\rho\bigr),z_k^\rho-z_{k-1}^\rho\right\rangle \label{op43},
		\end{multline}
		\begin{gather}
			\text{and}\quad\forall v\in\calZ:\quad \calR(v)\geq -\left\langle \xi_k^\rho+\D_z\calE\bigl(t_k^\rho,u_k^\rho,z_k^\rho\bigr),v\right\rangle.\label{op44}
		\end{gather}
	\end{lemma}
	We emphasize that, for $k=0$, $u_0^\rho$ and $z_0^\rho$ are the iterates calculated within the first outer loop of Algorithm~\ref{alg1} and only $z_{-1}^\rho=z_0$ and $t_0^\rho=t_0$ are set to the initial values.
		\begin{proof}
			By \eqref{fixpoint:z} $z_k^\rho$ minimizes $\calE\bigl(t_k^\rho,u_k^\rho,\cdot\bigr)+\psi_\rho\bigl(\cdot-z_{k-1}^\rho\bigr)$ over $\calZ$, where
			\begin{align*}
				\psi_\rho\coloneq \calR+\chi_\rho\quad\text{with characteristic function}\quad \chi_\rho(v)=\begin{cases}
					0, & \text{ if } \norm{v}_\calV\leq \rho\\
					\infty, &\text{ otherwise}.
				\end{cases}
			\end{align*}
			Therefore we have
			\begin{align}\label{stationarity}
				-\D_z\calE\bigl(t_k^\rho,u_k^\rho,z_k^\rho\bigr)\in\partial^\calZ\psi_\rho\bigl(z_k^\rho-z_{k-1}^\rho\bigr)
			\end{align}
			respectively by Fenchel's identity (see Lemma~\ref{fenchel})
			\begin{align*}
				\psi_\rho\bigl(z_k^\rho-z_{k-1}^\rho\bigr)+\psi_\rho^{*^\calZ}\bigl(-\D_z\calE\big(t_k^\rho,u_k^\rho,z_k^\rho\big)\bigr)=\left\langle -\D_z\calE\bigl(t_k^\rho,u_k^\rho,z_k^\rho\bigr),z_k^\rho-z_{k-1}^\rho\right\rangle.
			\end{align*}
			By Lemma~\ref{lemmaA1knees2018} the convex conjugate of $\psi_\rho$ is calculated to be
			\begin{align*}
				\psi_\rho^{*^\calZ}(\eta)=\rho\,\dist_{\calV^\ast}\bigl(\eta,\partial^\calZ\calR(0)\bigr)\quad\forall \eta\in\calZ^\ast.
			\end{align*}
			In particular, $\psi_\rho^{*^\calZ}\bigl(-\D_z\calE\big(t_k^\rho,u_k^\rho,z_k^\rho\big)\bigr)<\infty$ implies that $\left(\partial^\calZ\calR(0)+\D_z\calE\big(t_k^\rho,u_k^\rho,z_k^\rho\big)\right)\cap \calV^\ast\neq\emptyset$ and $\dist_{\calV^\ast}\Bigl(-\D_z\calE\bigl(t_k^\rho,u_k^\rho,z_k^\rho\bigr),\partial^\calZ\calR(0)\Bigr)<\infty$.	
			By $\norm{z_k^\rho-z_{k-1}^\rho}_\calV\leq \rho$, it follows that $\psi_\rho\bigl(z_k^\rho-z_{k-1}^\rho\bigr)=\calR\bigl(z_k^\rho-z_{k-1}^\rho\bigr)$, so that
			\begin{align*}
				\calR\bigl(z_k^\rho-z_{k-1}^\rho\bigr)+\rho\,\dist_{\calV^\ast}\Bigl(-\D_z\calE\bigl(t_k^\rho,u_k^\rho,z_k^\rho\bigr),\partial^\calZ\calR(0)\Bigr)=\left\langle-\D_z\calE\bigl(t_k^\rho,u_k^\rho,z_k^\rho\bigr),z_k^\rho-z_{k-1}^\rho\right\rangle,
			\end{align*}
			which is \eqref{op42}. 
			Moreover, the subdifferential sum rule from Lemma~\ref{lemma:subdiffsumrule} is applicable to $\psi_\rho=\calR+\chi_\rho$, since $\calR:\calZ\to\R_\infty$ and $\chi_\rho:\calZ\to\R_\infty$ are convex, $\calR$ is lower semicontinuous on $\calZ$ by Lemma~\ref{lemma:continuity} and it holds $0\in\dom(\calR)\cap\dom(\chi_\rho)$ and $\chi_\rho$ is continuous in $0$. Therefore, in view of \eqref{stationarity} there exists a $\xi_k^\rho\in\partial^\calZ\chi_\rho\bigl(z_k^\rho-z_{k-1}^\rho\bigr)$ with
			\begin{align}\label{stationarity2}
				-\xi_k^\rho-\D_z\calE\bigl(t_k^\rho,u_k^\rho,z_k^\rho\bigr)\in\partial^\calZ\calR\bigl(z_k^\rho-z_{k-1}^\rho\bigr)
			\end{align}
			respectively by Fenchel's identity (Lemma~\ref{fenchel})
			\begin{align*}
				\calR\bigl(z_k^\rho-z_{k-1}^\rho\bigr)+\calR^{\ast^\calZ}\bigl(-\xi_k^\rho-\D_z\calE\bigl(t_k^\rho,u_k^\rho,z_k^\rho\bigr)\bigr)=\left\langle -\xi_k^\rho-\D_z\calE\bigl(t_k^\rho,u_k^\rho,z_k^\rho\bigr),z_k^\rho-z_{k-1}^\rho\right\rangle.
			\end{align*}
			If we substract (\ref{op42}) from this, we obtain
			\begin{align*}
				\calR^{\ast^\calZ}\bigl(-\xi_k^\rho-\D_z\calE\bigl(t_k^\rho,u_k^\rho,z_k^\rho\bigr)\bigr)-\rho\,\dist_{\calV^\ast}\Bigl(-\D_z\calE\bigl(t_k^\rho,u_k^\rho,z_k^\rho\bigr),\partial^\calZ\calR(0)\Bigr)=\left\langle -\xi_k^\rho,z_k^\rho-z_{k-1}^\rho\right\rangle.
			\end{align*}
			By \eqref{prop:convexconjugatepos1hom} we have
			\begin{align*}
				\forall \xi\in\calZ^\ast:\quad\calR^{*^\calZ}(\xi)=\chi_{_{\partial^\calZ\calR(0)}}(\xi)=
				\begin{cases}
					0, &\text{if } \xi\in\partial^\calZ\calR(0)\\
					\infty, &\text{ otherwise}.
				\end{cases}
			\end{align*}
			Also we know from Lemma~\ref{prop:subdifferentialpos1homogen} that $\partial^\calZ\calR\bigl(z_k^\rho-z_{k-1}^\rho\bigr)\subset\partial^\calZ\calR(0)$. Hence it follows from \eqref{stationarity2} that
			\begin{align*}
				\calR^{\ast^\calZ}\bigl(-\xi_k^\rho-\D_z\calE\bigl(t_k^\rho,u_k^\rho,z_k^\rho\bigr)\bigr)=0,
			\end{align*}
			so the first equation in \eqref{op41} is proved. Moreover, by the characterization of $\partial^\calZ\calR(0)$ from Lemma~\ref{prop:subdifferentialpos1homogen} relation \eqref{stationarity2} implies \eqref{op44}.
			
			By Lemma~\ref{lemmaA310} we have $\partial^\calZ\chi_\rho(v)=\partial^\calV\chi_\rho(v)$ for all $v\in\calZ$, thus $\xi_k^\rho\in\partial^\calV\chi_\rho\bigl(z_k^\rho-z_{k-1}^\rho\bigr)\subset\calV^\ast$. Therefore, Fenchel's identity (Lemma~\ref{fenchel}), which is applicable since $\chi_\rho$ is lower semicontinuous on $\calV$, implies
			\begin{align*}
				\chi_\rho(z_k^\rho-z_{k-1}^\rho)+\chi_\rho^{*^\calV}(\xi_k^\rho)=\left\langle \xi_k^\rho,z_k^\rho-z_{k-1}^\rho\right\rangle
			\end{align*}
			and taking into account that $\chi_\rho^{*^\calV}(\xi)=\sup\Set{\bigl\langle\xi,\frac{v}{\rho}\bigr\rangle}{v\in\calV\text{ with }\bnorm{\frac{v}{\rho}}_\calV\leq 1}=\rho\,\norm{\xi}_{\calV^*}$ for all $\xi\in\calV^\ast$, this is equivalent to
			$\bnorm{z_k^\rho-z_{k-1}^\rho}_\calV\leq \rho$ and $\rho\,\bnorm{\xi_k^\rho}_{\calV^*}=\left\langle\xi_k^\rho,z_k^\rho-z_{k-1}^\rho\right\rangle$. Thus the second equality in \eqref{op41} is proved. Moreover, this implies
			\begin{align*}
				\bnorm{\xi_k^\rho}_{\calV^*}\bnorm{z_k^\rho-z_{k-1}^\rho}_\calV\leq \rho\,\bnorm{\xi_k^\rho}_{\calV^*}=\left\langle\xi_k^\rho,z_k^\rho-z_{k-1}^\rho\right\rangle
				\leq \bnorm{\xi_k^\rho}_{\calV^*}\bnorm{z_k^\rho-z_{k-1}^\rho}_\calV,
			\end{align*}
			what shows the third equality in \eqref{op41} as well as \eqref{op40}. Finally, by \eqref{op40} $\bnorm{z_k^\rho-z_{k-1}^\rho}_\calV<\rho$ implies $\bnorm{\xi_k^\rho}_{\calV^*}=0$ and thus, by \eqref{op41}, also $\dist_{\calV^\ast}\Bigl(-\D_z\calE\bigl(t_k^\rho,u_k^\rho,z_k^\rho\bigr),\partial^\calZ\calR(0)\Bigr)=0$, such that we obtain \eqref{op43} from \eqref{op42}.
		\end{proof}
		The next proposition guarantees that the final time $T$ is reached after a finite number of Efendiev\&Mielke iterations $N(\rho)$. Moreover, it provides some crucial estimates for the further convergence analysis. Essentially, \eqref{prop:fi1} will yield a bound on the artificial time $S_\rho$, with \eqref{prop:fi2} we can apply the lower semicontinuity result from Lemma~\ref{lemma:lowersemicon} to handle the unbounded dissipation potential and \eqref{prop:fi3} is needed to derive estimate \eqref{remainder} for the remainder $r_\rho$ of the discrete energy dissipation balance \eqref{energydissipationbalance}. We basically follow \cite[Proposition~2.3]{Knees2018} respectively \cite[Proposition~3.2.12]{Sievers2020}.
		
		\begin{proposition}[Crucial estimates and finite number $N(\rho)$ of Efendiev\&Mielke iterations]\label{prop:crucial}
			Let $\tilde{p}$ and $\alpha$ be chosen as in \eqref{tildep} and \eqref{alpha}.
			Assume \eqref{assumption:initialvalue} and \eqref{assumption:EMiterates}. Additionally, assume that there exists $\tilde{u}\in\calU^{\tilde{p}}$ such that \mbox{$\D_z\calE(t_0,\tilde{u},z_0)\in\calV^\ast$}.
			Then for every $\rho>0$ there exists a minimal number $N(\rho)\in\bbN$ of E\&M iterations such that the final time $T$ is reached, i.e. $t_{N(\rho)}^\rho=T$. Moreover, there exists a constant $C>0$ such that for all $\rho>0$ and all $k\in\N_0$ with $\xi_k^\rho$ from Lemma~\ref{prop:optimality} 
			we have
			\begin{align}\label{prop:fi1}
				\big\Vert\xi_k^\rho\big\Vert_{\calV^*}+\sum_{j=0}^{k}\norm{z_{j+1}^\rho-z_j^\rho}_\calZ&\leq C,
				\\
				\dist_{\calV^\ast}\Bigl(-\D_z\calE\big(t_k^\rho,u_k^\rho,z_k^\rho\big),\partial^\calZ\calR(0)\Bigr)&=\big\Vert\xi_k^\rho\big\Vert_{\calV^*}\leq C,\label{prop:fi2}\\
				\text{and}\qquad 
				\sum_{j=0}^{k}\norm{z_{j+1}^\rho-z_j^\rho}_\calZ^2&\leq C\, \rho.\label{prop:fi3}
			\end{align}
		\end{proposition}
		\begin{remark}
			Observe that for $k=0$ the iterates are not the initial values but the iterates of the first E\&M loop.
		\end{remark}
		\begin{proof}
			Let $k\in\bbN_0$ and $\rho>0$ be arbitrary. Inserting \eqref{op41} in \eqref{op42} with $k$ replaced by $k+1$ and subtracting the resulting equation from \eqref{op44} with $v\coloneq z_{k+1}^\rho-z_k^\rho$ leads to
			\begin{multline*}
				\bnorm{\xi_{k+1}^\rho}_{\calV^*}\bnorm{z_{k+1}^\rho-z_k^\rho}_\calV-\left\langle \xi_k^\rho,z_{k+1}^\rho-z_k^\rho\right\rangle\leq \left\langle \D_z\calE\bigl(t_{k+1}^\rho,u_{k+1}^\rho,z_{k+1}^\rho\bigr)-\D_z\calE\bigl(t_{k}^\rho,u_{k}^\rho,z_{k}^\rho\bigr),z_{k}^\rho-z_{k+1}^\rho\right\rangle.
			\end{multline*}
			Thus by \eqref{frechetz} we obtain
			\begin{align*}
				\bigr\Vert\xi_{k+1}^\rho\bigl\Vert_{\calV^*}\norm{z_{k+1}^\rho-z_k^\rho}_\calV&-\left\langle \xi_k^\rho,z_{k+1}^\rho-z_k^\rho\right\rangle+\kappa\norm{z_{k+1}^\rho-z_k^\rho}_\calZ^2\\
				&\leq \int_\Omega \left(z_{k}^\rho-z_{k+1}^\rho\right)\biggl(z_{k+1}^\rho\Bigl(\boldC\strain\bigl(u_{k+1}^\rho\bigr)\Bigr)\colon\strain\bigl(u_{k+1}^\rho\bigr)-z_k^\rho \left(\boldC\strain\bigl(u_{k}^\rho\bigr)\right)\colon\strain\bigl(u_{k}^\rho\bigr)\biggr)\dd{x}\\
				&=
				\int_\Omega \left(z_{k}^\rho-z_{k+1}^\rho\right) z_{k+1}^\rho \biggl(\boldC\Bigl(\strain\bigl(u_{k+1}^\rho\bigr)+\strain\big(u_{k}^\rho\big)\Bigr)\biggr)\colon\Bigl(\strain\big(u_{k+1}^\rho\big)-\strain\big(u_{k}^\rho\big)\Bigr)\dd{x}\\
				&\hspace{5.2cm}-\int_\Omega \underbrace{\left(z_{k}^\rho-z_{k+1}^\rho\right)^2\Bigl(\boldC\strain\bigl(u_{k}^\rho\bigr)\Bigr)\colon\strain\bigl(u_{k}^\rho\bigr)}_{\geq 0}\dd{x}.
			\end{align*}
		By Hölder's inequality with $\frac{1}{r_2}+\frac{2}{\tilde{p}}=1$, \eqref{prop_be1}, the uniform bound for the displacement iterates from \eqref{prop_be2} and Lemma~\ref{lemmaA1} with $r_1=\frac{\tilde{p}p}{p-\tilde{p}}$ the first term on the right hand side can be estimated by
		\begin{align*}
				r.h.s.&\leq C\norm{z_{k}^\rho-z_{k+1}^\rho}_{L^{r_2}(\Omega)}\norm{\strain\big(u_{k+1}^\rho\big)-\strain\big(u_{k}^\rho\big)}_{L^{\tilde{p}}(\Omega)}\left(\norm{\strain\big(u_{k+1}^\rho\big)}_{L^{\tilde{p}}(\Omega)}+\norm{\strain\big(u_{k}^\rho\big)}_{L^{\tilde{p}}(\Omega)}\right)\\
				&\leq C\norm{z_{k}^\rho-z_{k+1}^\rho}_\calV\norm{u_{k+1}^\rho-u_{k}^\rho}_{\calU^{\tilde{p}}}\\
				&\leq C\norm{z_{k}^\rho-z_{k+1}^\rho}_\calV \left(t_{k+1}^\rho-t_{k}^\rho+\norm{z_{k}^\rho-z_{k+1}^\rho}_{L^{r_1}(\Omega)}\right).
			\end{align*}
			Hence we obtain together with $\left\langle \xi_k^\rho,z_{k+1}^\rho-z_k^\rho\right\rangle\leq \norm{\xi_{k}^\rho}_{\calV^*}\norm{z_{k}^\rho-z_{k+1}^\rho}_\calV$
			\begin{align}\label{proof:fi}\notag
				\norm{\xi_{k+1}^\rho}_{\calV^*}\norm{z_{k}^\rho-z_{k+1}^\rho}_\calV-\norm{\xi_{k}^\rho}_{\calV^*}&\norm{z_{k}^\rho-z_{k+1}^\rho}_\calV+\kappa\norm{z_{k}^\rho-z_{k+1}^\rho}_\calZ^2\\
				&\hspace{10pt}\leq C\norm{z_{k}^\rho-z_{k+1}^\rho}_\calV \left(t_{k+1}^\rho-t_{k}^\rho+\norm{z_{k}^\rho-z_{k+1}^\rho}_{L^{r_1}(\Omega)}\right).
			\end{align}
			Division by $\norm{z_{k}^\rho-z_{k+1}^\rho}_\calV$, where we take into account the continuous embedding $\calZ\hookrightarrow\calV$, yields
			\begin{align*}
				\norm{\xi_{k+1}^\rho}_{\calV^*}-\big\Vert\xi_{k}^\rho\big\Vert_{\calV^*}+c\norm{z_{k}^\rho-z_{k+1}^\rho}_\calZ\leq C \left(t_{k+1}^\rho-t_{k}^\rho+\norm{z_{k}^\rho-z_{k+1}^\rho}_{L^{r_1}(\Omega)}\right).
			\end{align*}
			By Ehrling's Lemma~\ref{lemma:ehrling} with $\varepsilon=\frac{c}{2C}>0$ we find a constant $C(\varepsilon)>0$ such that
			\begin{align*}
				\norm{z_{k}^\rho-z_{k+1}^\rho}_{L^{r_1}(\Omega)}\leq \varepsilon \norm{z_{k}^\rho-z_{k+1}^\rho}_\calZ+C(\varepsilon)\norm{z_{k}^\rho-z_{k+1}^\rho}_\calX,
			\end{align*}
			such that by absorption and finally by \eqref{R2} it follows
			\begin{align*}
				\norm{\xi_{k+1}^\rho}_{\calV^*}-\big\Vert\xi_{k}^\rho\big\Vert_{\calV^*}+\frac{c}{2}\norm{z_{k}^\rho-z_{k+1}^\rho}_\calZ&\leq C \left(t_{k+1}^\rho-t_{k}^\rho+C(\varepsilon)\norm{z_{k}^\rho-z_{k+1}^\rho}_{\calX}\right)\\
				&\hspace{-4pt}\overset{\text{\eqref{R2}}}{\leq}C \Big(t_{k+1}^\rho-t_{k}^\rho+\calR\big(z_{k+1}^\rho-z_{k}^\rho\big)\Big).
			\end{align*}
			Summation with respect to $k\in\N_0$ implies
			\begin{align}\label{proof:summation}
				\norm{\xi_{k+1}^\rho}_{\calV^*}+\frac{c}{2}\sum_{j=0}^{k}\norm{z_{j}^\rho-z_{j+1}^\rho}_\calZ&\leq \big\Vert\xi_{0}^\rho\big\Vert_{\calV^*}+C\left( T+\sum_{j=0}^{k}\calR\big(z_{j+1}^\rho-z_{j}^\rho\big)\right).
			\end{align}
			Moreover, by definition of the dissipation potential $\calR$ and by the minimality in \eqref{fixpoint:z}, $z_{j+1}^\rho-z_{j}^\rho\leq 0$ has to be satisfied a.e. in $\Omega$ for all $j\in\N_0$. Therefore, it follows 
			for all $k\in\N_0$
			\begin{align*}
				0\leq\sum_{j=0}^{k}\calR\big(z_{j+1}^\rho-z_{j}^\rho\big)&=\sum_{j=0}^{k} \int_\Omega \kappa\abs{z_{j+1}^\rho(x)-z_{j}^\rho(x)}\dd{x}=\kappa\sum_{j=0}^{k} \int_\Omega z_{j}^\rho(x)-z_{j+1}^\rho(x)\dd{x}\\
				&=\kappa \int_{\Omega} \underbrace{z_0^\rho(x)-z_{k+1}^\rho(x)}_{\leq 1\text{ by \eqref{prop_be1}}}\dd{x}\leq C.
			\end{align*}
			Hence, to finish the proof of \eqref{prop:fi1} based on \eqref{proof:summation}, 		
			it only remains to show that there exists $C_0>0$ with $\big\Vert\xi_{0}^\rho\big\Vert_{\calV^*}\leq C_0$. This can be obtained analogously to the above:
			Inserting \eqref{op41} in \eqref{op42} for $k=0$, where we take into account that $z_{-1}^\rho=z_0$, and adding a zero yields
			\begin{multline*}
				\calR\big(z_0^\rho-z_0\big)+\big\Vert\xi_0^\rho\big\Vert_{\calV^\ast}\big\Vert z_0^\rho-z_0\big\Vert_\calV+\big\langle \D_z\calE\bigl(t_0^\rho,u_0^\rho,z_0^\rho\bigr)-\D_z\calE\big(t_0,\tilde{u},z_0\big),z_0^\rho-z_0\big\rangle\\ =\big\langle -\D_z\calE\big(t_0,\tilde{u},z_0\big),z_0^\rho-z_0\big\rangle.
			\end{multline*}
			Thus by \eqref{frechetz} we obtain
			\begin{align*}
				\underbrace{\calR\bigl(z_0^\rho-z_0\bigr)}_{\geq 0}&+\bigl\Vert\xi_0^\rho\bigr\Vert_{\calV^\ast}\bigl\Vert z_0^\rho-z_0\bigr\Vert_\calV+\kappa\bigl\Vert z_0^\rho-z_0\bigr\Vert_\calZ^2\\
				&\hspace{-1cm}=\int_\Omega \bigl(z_0-z_0^\rho\bigr)\biggl(z_0^\rho \Bigl(\boldC\strain\bigl(u_0^\rho\bigr)\Bigr)\colon\strain\bigl(u_0^\rho\bigr)-z_0\Bigl(\boldC\strain\bigl(\tilde{u}\bigr)\Bigr)\colon\strain\bigl(\tilde{u}\bigr)\biggr)\dd{x}-\bigl\langle \D_z\calE\bigl(t_0,\tilde{u},z_0\bigr),z_0^\rho-z_0\bigr\rangle\\
				&\hspace{-1cm}=
				\int_\Omega \left(z_0-z_0^\rho\right) z_{0}^\rho \biggl(\boldC\Bigl(\strain\bigl(u_{0}^\rho\bigr)+\strain\big(\tilde{u}\big)\Bigr)\biggr)\colon\Bigl(\strain\big(u_{0}^\rho\big)-\strain\big(\tilde{u}\big)\Bigr)\dd{x}\\
				&\hspace{3.1cm}-\int_\Omega \underbrace{\left(z_0-z_0^\rho\right)^2\Bigl(\boldC\strain\bigl(\tilde{u}\bigr)\Bigr)\colon\strain\bigl(\tilde{u}\bigr)}_{\geq 0}\dd{x}-\big\langle \D_z\calE\big(t_0,\tilde{u},z_0\big),z_0^\rho-z_0\big\rangle.
			\end{align*}
		By Hölder's inequality with $\frac{1}{r_2}+\frac{2}{\tilde{p}}=1$, \eqref{prop_be1}
		and the uniform bound for the displacement iterates from \eqref{prop_be2} 
		the right hand side can be estimated by
		\begin{align*}
			r.h.s.&\leq C\bnorm{z_{0}-z_{0}^\rho}_{L^{r_2}(\Omega)}\bnorm{\strain\big(u_{0}^\rho\big)-\strain\big(\tilde{u}\big)}_{L^{\tilde{p}}(\Omega)}\left(\bnorm{\strain\big(u_{0}^\rho\big)}_{L^{\tilde{p}}(\Omega)}+\bnorm{\strain\big(\tilde{u}\big)}_{L^{\tilde{p}}(\Omega)}\right)\\
			&\hspace{8cm}+\bnorm{\D_z\calE\big(t_0,\tilde{u},z_0\big)}_{\calV^\ast}\big\Vert z_0^\rho-z_0\big\Vert_{\calV}\\
			&\leq C\bnorm{z_{0}-z_{0}^\rho}_\calV\bnorm{u_{0}^\rho-\tilde{u}}_{\calU^{\tilde{p}}}+\bnorm{\D_z\calE\big(t_0,\tilde{u},z_0\big)}_{\calV^\ast}\bigl\Vert z_0^\rho-z_0\bigr\Vert_{\calV}.
		\end{align*}
		Hence we have due to the additional assumption \mbox{$\D_z\calE(t_0,\tilde{u},z_0)\in\calV^\ast$} for a suitable $\tilde{u}\in\calU^{\tilde{p}}$
		\begin{align*}
			\bnorm{\xi_{0}^\rho}_{\calV^*}\bnorm{z_{0}-z_{0}^\rho}_\calV+\kappa\bnorm{z_{0}-z_{0}^\rho}_\calZ^2
			\leq C\bnorm{z_{0}-z_{0}^\rho}_\calV.
		\end{align*}
		So dividing by $\bnorm{z_{0}-z_{0}^\rho}_\calV$, where we take into account the continuous embedding $\calZ\hookrightarrow\calV$, yields
		\begin{align}\label{boundxi0}
			\bnorm{\xi_{0}^\rho}_{\calV^*}+\kappa\bnorm{z_{0}-z_{0}^\rho}_\calZ
			\leq C.
		\end{align}
		Altogether this proves \eqref{prop:fi1}.
		With this inequality at hand we can now show, that the final time $T$ is reached after a finite number of iterations. Suppose this is not the case, i.e. there exists $t_\ast\leq T$ such that $\lim\limits_{k\to\infty} t_k^\rho=t_\ast$. Then in particular the sequence $\left(t_{k+1}^\rho-t_k^\rho\right)_{k\in\bbN}$ tends to zero. But by construction this implies that $\left(\norm{z_{k+1}^\rho-z_k^\rho}_\calV\right)_{k\in\bbN}$ tends to $\rho$, what contradicts the convergence of the series $\sum_{j=0}^{\infty}\norm{z_{j+1}^\rho-z_j^\rho}_\calV$, guaranteed by (\ref{prop:fi1}). Therefore, there exists $N(\rho)\in\bbN$ such that $t_{N(\rho)}^\rho=T$.
		
		Moreover, \eqref{prop:fi2} is a direct consequence of \eqref{prop:fi1} and \eqref{op41}.
		Finally, estimate \eqref{prop:fi3} can be shown in an analogous way to \cite[Prop.3.2.12 (3.2.31)]{Sievers2020}. From \eqref{op40} and $\bnorm{z_k^\rho-z_{k+1}^\rho}_\calV\leq \rho$ it follows, that
			\begin{align*}
				-\bigl\Vert\xi_k^\rho\bigr\Vert_{\calV^\ast}\bnorm{z_{k-1}^\rho-z_{k}^\rho}_\calV\overset{\text{\eqref{op40}}}{=}-\bigl\Vert\xi_k^\rho\bigr\Vert_{\calV^\ast} \rho\leq 
				-\bigl\Vert\xi_k^\rho\bigr\Vert_{\calV^\ast}\bnorm{z_k^\rho-z_{k+1}^\rho}_\calV.
			\end{align*}
			Therefore, \eqref{proof:fi} implies
			\begin{align*}
				\bnorm{\xi_{k+1}^\rho}_{\calV^*}&\bnorm{z_{k}^\rho-z_{k+1}^\rho}_\calV-\bigl\Vert\xi_{k}^\rho\bigr\Vert_{\calV^*}\bnorm{z_{k-1}^\rho-z_{k}^\rho}_\calV+\kappa\bnorm{z_{k}^\rho-z_{k+1}^\rho}_\calZ^2\\
				&\leq C\bnorm{z_{k}^\rho-z_{k+1}^\rho}_\calV \left(t_{k+1}^\rho-t_{k}^\rho+\bnorm{z_{k}^\rho-z_{k+1}^\rho}_{L^{r_1}(\Omega)}\right)\leq C \rho \left(t_{k+1}^\rho-t_{k}^\rho+\bnorm{z_{k}^\rho-z_{k+1}^\rho}_{\calZ}\right).
			\end{align*}
			Summation with respect to $k$ together with \eqref{boundxi0} and \eqref{prop:fi1} yields
			\begin{align*}
				\norm{\xi_{k+1}^\rho}_{\calV^*}\norm{z_{k}^\rho-z_{k+1}^\rho}_\calV+\kappa\sum_{j=0}^{k}\norm{z_{j}^\rho-z_{j+1}^\rho}_\calZ^2&\leq \underbrace{\big\Vert\xi_{0}^\rho\big\Vert_{\calV^*}}_{\overset{\text{\eqref{boundxi0}}}{\leq} C}\underbrace{\big\Vert z_{-1}^\rho-z_{0}^\rho\big\Vert_\calV}_{\leq \rho}+C \rho\Bigl( T+\underbrace{\sum_{i=1}^{k}\norm{z_{k}^\rho-z_{k+1}^\rho}_\calZ}_{\leq C\text{ by \eqref{prop:fi1}}}\Bigr)\\
				&\leq C \rho,
			\end{align*}
			i.e. \eqref{prop:fi3}.
		\end{proof}
		\begin{corollary}[Upper bound for the number of E\&M iterations]\label{prop:Nestimate}
			The minimal number $N(\rho)\in\N$ of Efendiev\&Mielke iterations such that the final time $T$ is reached from the previous Proposition~\ref{prop:crucial} satisfies the following estimate
			\begin{align}\label{boundartificial}
				\exists C>0,\, \forall \rho>0:\quad N(\rho)\leq 1+\frac{C}{\rho}.
			\end{align}
		\end{corollary}
			\begin{proof}
				Let $N(\rho)\in\N$ be the minimal number of E\&M iterations such that $t_{N(\rho)}^\rho=T$. 
				
				Due to the minimality of $N(\rho)$ we have $t_{k}^\rho<T$ for all $k\in\{0,\ldots,N(\rho)-1\}$. Thus it holds by construction
				$t_{k}^\rho=t_{k-1}^\rho+\rho-\norm{z_{k-1}^\rho-z_{k-2}^\rho}_\calV$ for all $k\in\{1,\ldots,N(\rho)-1\}$ and summation yields
				\begin{align*}
					\bigl(N(\rho)-1\bigr)\rho=t_{N(\rho)-1}+\sum_{k=1}^{N(\rho)-1}\norm{z_{k-1}^\rho-z_{k-2}^\rho}_\calV\leq T+C,
				\end{align*}
				i.e. $N(\rho)\rho\leq \rho+T+C$ and $N(\rho)\leq 1+\frac{C}{\rho}$.

			\end{proof}

%
		
		
		\subsection{Construction and properties of interpolants}\label{sec:interpolants}
		Within this section let $\rho>0$ be given and let $N(\rho)\in\N$ be the minimal number of E\&M iterations such that $t_{N(\rho)}^\rho=T$. 
		We set $S_\rho\coloneq N(\rho)\rho$ and $s_k^\rho\coloneq k\rho$ for $k\in\N_0\cup\{-1\}$. To include the case of a jump at $t_0$, we set
		\begin{align*}
			t_{-1}^\rho\coloneq t_0\quad\text{and}\quad u_{-1}^\rho\coloneq u_{0,1}^\rho.
		\end{align*}
		Further, let us define the affine linear and constant right respectively left continuous interpolants $\hat{\square}_\rho, \underline{\square}_\rho$ and $\overline{\square}_\rho$ of the E\&M iterates, where $\square\in\{t,u,z\}$, by setting 
		for all $k\in\N_0$ and $s\in\left[ s_{k-1}^\rho,s_k^\rho \right)$
		\begin{align*}
			\hat{t}_\rho(s)&\coloneq t_{k-1}^\rho+\left(s-s_{k-1}^\rho\right)\rho^{-1}\left(t_k^\rho-t_{k-1}^\rho\right)
			,\quad \underline{t}_\rho(s)\coloneq t_{k-1}^\rho,
			\\
			\hat{u}_\rho(s)&\coloneq u_{k-1}^\rho+\left(s-s_{k-1}^\rho\right)\rho^{-1}\left(u_k^\rho-u_{k-1}^\rho\right),\quad \underline{u}_\rho(s)\coloneq u_{k-1}^\rho,
			\\
			\hat{z}_\rho(s)&\coloneq z_{k-1}^\rho+\left(s-s_{k-1}^\rho\right)\rho^{-1}\left(z_k^\rho-z_{k-1}^\rho\right),\quad \underline{z}_\rho(s)\coloneq z_{k-1}^\rho,
		\end{align*}
		and for $s\in \left(s_{k-1}^\rho,s_k^\rho\right]$
		\begin{align*}
			\overline{t}_\rho(s)\coloneq t_{k}^\rho,\quad \overline{u}_\rho(s)\coloneq u_{k}^\rho,\quad \overline{z}_\rho(s)\coloneq z_{k}^\rho.
		\end{align*}
		In particular, we have $ \hat{t}_\rho(S_\rho)=\overline{t}_\rho(S_\rho)=t_{N(\rho)}^\rho=T$ and due to \eqref{boundartificial} it holds $\sup\limits_{\rho\in\left(0,1\right)}S_\rho<\infty$.
		Note that $\hat{t}_\rho(s)$ stays constantly at $T$ for all $s\geq S_\rho$.
		Furthermore, observe that by construction 
		\begin{align*}
			\forall k\in\left\{ 1,\ldots,N(\rho)-1 \right\}:&\quad t_k^\rho-t_{k-1}^\rho+\bnorm{z_{k-1}^\rho-z_{k-2}^\rho}_\calV=\rho\\
			\text{and}&\quad t_{N(\rho)}^\rho-t_{N(\rho)-1}^\rho+\bnorm{z_{N(\rho)-1}^\rho-z_{N(\rho)-2}^\rho}_\calV\leq\rho,
		\end{align*}
		i.e. for almost all $s\in\left[0,S_\rho-\rho\right]$ it is $\hat{t}_\rho'(s)+\bnorm{\hat{z}_\rho'(s-\rho)}_\calV=1$ and for almost all $s\in\left[S_\rho-\rho,S_\rho\right]$ we have $\hat{t}_\rho'(s)+\bnorm{\hat{z}_\rho'(s-\rho)}_\calV\leq 1$.
		Moreover, for all $k\in\N_0$ it is $\abs{t_{k}^\rho-t_{k-1}^\rho}\leq\rho$ by construction, and by \eqref{fixpoint:z} we have $\bnorm{z_k^\rho-z_{k-1}^\rho}_\calV\leq \rho$. This implies by Lemma~\ref{lemmaA1}, \eqref{fixpoint:u} and with $\alpha\geq r_1$ that also
		\begin{align*}
			\bnorm{u_k^\rho-u_{k-1}^\rho}_{\calU^{\tilde{p}}}
			\leq C\Big(\underbrace{t_{k}^\rho-t_{k-1}^\rho}_{\leq \rho}+\underbrace{\bnorm{z_{k}^\rho-z_{k-1}^\rho}_{\calV}}_{\leq  \rho}\Big)\leq C\, \rho.
		\end{align*}
		This shows the Lipschitz continuity of $\hat{t}_\rho$, of $\hat{z}_\rho$ with respect to $\calV$ (both with Lipschitz constant equal to $1$) and of $\hat{u}_\rho$ with respect to $\calU^{\tilde{p}}$ with Lipschitz constants independent of $\rho$. In fact, \eqref{prop_be2} also gives us Lipschitz continuity of $\hat{z}_\rho$ with respect to $\calZ$ but with a Lipschitz constant dependent on $\rho$.
		In addition we obtain from the choice of the interpolants and \eqref{prop:fi3} for all $\tilde{S}>0$
		\begin{align*}
			\exists C>0,\,\forall \tilde{S}>0,\,\forall \rho>0:\quad \norm{\hat{z}_\rho'}_{L^2((0,\tilde{S});\calZ)}^2&=\int_{0}^{\tilde{S}}\norm{\hat{z}_\rho'(s)}_{\calZ}^2\dd{s}
			\leq \sum_{k=1}^{\lceil\tilde{S}/\rho\rceil}\int_{s_{k-1}^\rho}^{s_k^\rho}\rho^{-2}\norm{z_k^\rho-z_{k-1}^\rho}_\calZ^2\dd{s}\\
			&=\rho^{-1}\sum_{k=1}^{\lceil\tilde{S}/\rho\rceil}\norm{z_k^\rho-z_{k-1}^\rho}_\calZ^2\overset{\text{\eqref{prop:fi3}}}{\leq} C.
		\end{align*}
		
		In summary we obtain together with Lemma~\ref{prop:basicestimates} and Proposition~\ref{prop:crucial}
		\begin{align}\label{reachedT}
			\hat{t}_\rho(S_\rho)=\underline{t}_\rho(S_\rho)=T\quad\text{and}\quad\forall \tilde{S}>0:\quad\left(\hat{t}_\rho,\hat{u}_\rho,\hat{z}_\rho\right)\in W^{1,\infty}\bigl((-\rho,\tilde{S});\bbR\times\calU^{\tilde{p}}\times\calZ\bigr),
		\end{align}
		\begin{multline}
			\label{bound}
			\forall \tilde{S}>0:\quad \sup\limits_{\rho\in\left(0,1\right)} \Bigl( \norm{\hat{t}_\rho}_{W^{1,\infty}((-\rho,\tilde{S});\bbR)}+\norm{\hat{u}_\rho}_{W^{1,\infty}((-\rho,\tilde{S});\calU^{\tilde{p}})}+\norm{\hat{z}_\rho}_{W^{1,\infty}((-\rho,\tilde{S});\calV)}\\
			+\norm{\hat{z}_\rho'}_{L^{2}((0,\tilde{S});\calZ)}+\norm{\hat{z}_\rho}_{L^{\infty}((-\rho,\tilde{S});\calZ)}+S_\rho \Bigr)<\infty,
		\end{multline}
		\begin{gather}\label{normalization}
			\begin{split}
				\text{f.a.a. } s\in \left[-\rho,\infty\right):\quad \hat{t}_\rho'(s)\geq 0,\quad
				\text{f.a.a. } s\in\left[0,S_\rho-\rho\right]:\quad \hat{t}_\rho'(s)+\bnorm{\hat{z}_\rho'(s-\rho)}_\calV=1,\\
				\text{and}\qquad \text{f.a.a. } s\in\left[S_\rho-\rho,\infty\right):\quad \hat{t}_\rho'(s)+\bnorm{\hat{z}_\rho'(s-\rho)}_\calV\leq 1.
			\end{split}
		\end{gather}
		
		\begin{proposition}[Discrete energy dissipation balance]
			Let $\tilde{p}$ and $\alpha$ be chosen as in \eqref{tildep} and \eqref{alpha}.
			For all $\sigma_1,\sigma_2\in\left[-\rho,S_\rho\right]$ with $\sigma_1\leq\sigma_2$ the discrete energy dissipation balance
			\begin{multline}\label{energydissipationbalance}
				\calE\bigl(\hat{t}_\rho(\sigma_2),\hat{u}_\rho(\sigma_2),\hat{z}_\rho(\sigma_2)\bigr)+\int_{\sigma_1}^{\sigma_2}\calR\bigl(\hat{z}_\rho'(s)\bigr)+\bnorm{\hat{z}_\rho'(s)}_\calV\dist_{\calV^\ast}\Bigl(-\D_z\calE\bigl(\overline{t}_\rho(s),\overline{u}_\rho(s),\overline{z}_\rho(s)\bigr),\partial^{\calZ}\calR(0)\Bigr)\dd s\\
				=\calE\bigl(\hat{t}_\rho(\sigma_1),\hat{u}_\rho(\sigma_1),\hat{z}_\rho(\sigma_1)\bigr)+\int_{\sigma_1}^{\sigma_2}\partial_t\calE\bigl(\hat{t}_\rho(s),\hat{u}_\rho(s),\hat{z}_\rho(s)\bigr)\hat{t}_\rho'(s)\dd s+\int_{\sigma_1}^{\sigma_2}r_\rho(s)\dd s,
			\end{multline}
			 is satisfied, where the remainder $r_\rho$ is given by
			 \begin{multline*}
			 	r_\rho(s)=\left\langle \D_z\calE\bigl(\hat{t}_\rho(s),\hat{u}_\rho(s),\hat{z}_\rho(s)\bigr)-\D_z\calE\bigl(\overline{t}_\rho(s),\overline{u}_\rho(s),\overline{z}_\rho(s)\bigr),\hat{z}_\rho'(s)\right\rangle\\
			 	+\left\langle \D_u\calE\bigl(\hat{t}_\rho(s),\hat{u}_\rho(s),\hat{z}_\rho(s)\bigr),\hat{u}_\rho'(s)\right\rangle.
			 \end{multline*}
		 	Furthermore, we have the following estimate for the remainder $r_\rho$
		 	\begin{align}\label{remainder}
		 		\exists C>0,\,\forall \rho\in\left(0,1\right),\,\forall \sigma_1,\sigma_2\in\left[-\rho,S_\rho\right]\text{ with }\sigma_1<\sigma_2:\quad \int_{\sigma_1}^{\sigma_2}r_\rho(s)\dd s\leq C\rho.
		 	\end{align}
			Finally, the complementarity condition
			\begin{align}\label{prop:complementarity}
				\text{f.a.a. }s\in\left[0,S_\rho\right]:\quad\hat{t}_\rho'(s)\dist_{\calV^\ast}\Bigl(-\D_z\calE\bigl(\underline{t}_\rho(s),\underline{u}_\rho(s),\underline{z}_\rho(s)\bigr),\partial^{\calZ}\calR(0)\Bigr)=0
			\end{align}
			holds. 
			Moreover, we have
			\begin{align}\label{discrete:displacementEL}
				\text{f.a.a. } s\in\left[0,S_\rho\right],\,\forall w\in\calU:\quad\left\langle \D_u\calE\bigl(\overline{t}_\rho(s),\overline{u}_\rho(s),\overline{z}_\rho(s)\bigr),w\right\rangle=0.
			\end{align}
			
		\end{proposition}
		
		\begin{proof}
			Let $s\in\left[s_{k-1}^\rho,s_k^\rho\right)\subset\left[0,S_\rho\right]$ and assume that $\hat{t}_\rho'(s)>0$. Then by \eqref{normalization} it follows that $\bnorm{\hat{z}_\rho'(s-\rho)}_\calV<1$, i.e. $\bnorm{z_{k-1}^\rho-z_{k-2}^\rho}_\calV<\rho$. Hence, by \eqref{op40} we obtain $\bnorm{\xi_{k-1}^\rho}_{\calV^\ast}=0$ and thus by \eqref{op41} it is $\dist_{\calV^\ast}\Bigl(-\D_z\calE\bigl(\underline{t}_\rho(s),\underline{u}_\rho(s),\underline{z}_\rho(s)\bigr),\partial^{\calZ}\calR(0)\Bigr)=0$, which shows \eqref{prop:complementarity}.
		
			For the proof of the energy dissipation balance we apply a chain rule to \eqref{op43} and integrate over $s\in\left(\sigma_1,\sigma_2\right)$. Indeed, we have by \eqref{op43} and by the choice of the interpolants
			\begin{multline}\label{proof:debstern}
				\rho\,\calR\bigl(\hat{z}_\rho'(s)\bigr)+\rho\,\bnorm{\hat{z}_\rho'(s)}_\calV\dist_{\calV^\ast}\Bigl(-\D_z\calE\bigl(\overline{t}_\rho(s),\overline{u}_\rho(s),\overline{z}_\rho(s)\bigr),\partial^{\calZ}\calR(0)\Bigr)\\
				=\rho\left\langle -\D_z\calE\bigl(\overline{t}_\rho(s),\overline{u}_\rho(s),\overline{z}_\rho(s)\bigr), \hat{z}_\rho'(s)\right\rangle
			\end{multline}
			for all $k\in\N_0$ and $s\in\left(s_{k-1}^\rho,s_k^\rho\right)$.
			Analogously to \cite[Lemma~A.2.5]{Sievers2020} we obtain the following chain rule
			\begin{align}\label{proof:debstern2}
				\begin{split}
					\frac{\dd{}}{\dd{s}}\Bigl(\calE\bigl(\hat{t}_\rho(s),\hat{u}_\rho(s),&\hat{z}_\rho(s)\bigr)\Bigr)=\partial_t\calE\bigl(\hat{t}_\rho(s),\hat{u}_\rho(s),\hat{z}_\rho(s)\bigr)\hat{t}_\rho'(s)\\
					&+\left\langle \D_z\calE\bigl(\hat{t}_\rho(s),\hat{u}_\rho(s),\hat{z}_\rho(s)\bigr),\hat{z}_\rho'(s)\right\rangle + \left\langle \D_u\calE\bigl(\hat{t}_\rho(s),\hat{u}_\rho(s),\hat{z}_\rho(s)\bigr),\hat{u}_\rho'(s)\right\rangle
				\end{split}
			\end{align}
			for almost all $s\in\left(-\rho,S_\rho\right)$.
			Addition of $\eqref{proof:debstern}$ multiplied by $\frac{1}{\rho}$ and \eqref{proof:debstern2} and integrating over $s\in\left(\sigma_1,\sigma_2\right)$ leads to
			\begin{align*}
				\calE\bigl(\hat{t}_\rho(\sigma_2)&,\hat{u}_\rho(\sigma_2),\hat{z}_\rho(\sigma_2)\bigr)
				-\calE\bigl(\hat{t}_\rho(\sigma_1),\hat{u}_\rho(\sigma_1),\hat{z}_\rho(\sigma_1)\bigr)\\
				&\hspace{55pt}+\int_{\sigma_1}^{\sigma_2}\calR\bigl(\hat{z}_\rho'(s)\bigr)+\bnorm{\hat{z}_\rho'(s)}_\calV\dist_{\calV^\ast}\Bigl(-\D_z\calE\bigl(\overline{t}_\rho(s),\overline{u}_\rho(s),\overline{z}_\rho(s)\bigr),\partial^{\calZ}\calR(0)\Bigr)\dd{s}\\
				&=
				\int_{\sigma_1}^{\sigma_2}\partial_t\calE\bigl(\hat{t}_\rho(s),\hat{u}_\rho(s),\hat{z}_\rho(s)\bigr)\,\hat{t}_\rho'(s)+\left\langle \D_u\calE\bigl(\hat{t}_\rho(s),\hat{u}_\rho(s),\hat{z}_\rho(s)\bigr),\hat{u}_\rho'(s)\right\rangle \dd{s}\\
				&\hspace{55pt}+\int_{\sigma_1}^{\sigma_2}\left\langle \D_z\calE\bigl(\hat{t}_\rho(s),\hat{u}_\rho(s),\hat{z}_\rho(s)\bigr)-\D_z\calE\bigl(\overline{t}_\rho(s),\overline{u}_\rho(s),\overline{z}_\rho(s)\bigr),\hat{z}_\rho'(s)\right\rangle\dd{s}\\
				&=\int_{\sigma_1}^{\sigma_2}\partial_t\calE\bigl(\hat{t}_\rho(s),\hat{u}_\rho(s),\hat{z}_\rho(s)\bigr)\,\hat{t}_\rho'(s) \dd{s}
				+\int_{\sigma_1}^{\sigma_2} r_\rho(s)\dd{s}.
			\end{align*}
			Finally, we will show estimate \eqref{remainder} of the remainder $r_\rho$. First, note that by the choice of the interpolants we have for all $k\in\N_0$ and $s\in\left(s_{k-1}^\rho,s_{k}^\rho\right)$
			\begin{align*}
				\hat{t}_\rho(s)-\overline{t}_\rho(s)&=
				t_{k-1}^\rho-t_{k}^\rho+\left(s-s_{k-1}^\rho\right)\rho^{-1}\left(t_{k}^\rho-t_{k-1}^\rho\right)
				=\left(s-s_{k-1}^\rho-\rho\right)\rho^{-1}\left(t_{k}^\rho-t_{k-1}^\rho\right)\\
				&=\left(s-s_{k}^\rho\right)\hat{t}_\rho'(s),
			\end{align*}
			and analogously
			\begin{align*} 
				\hat{u}_\rho(s)-\overline{u}_\rho(s)=\left(s-s_{k}^\rho\right)\hat{u}_\rho'(s)\quad\text{and}\quad \hat{z}_\rho(s)-\overline{z}_\rho(s)=\left(s-s_{k}^\rho\right)\hat{z}_\rho'(s).
			\end{align*}
			Thus by \eqref{frechetz} 
			we obtain for the first term of the remainder $r_\rho$ for all $k\in\N_0$ and $s\in\left(s_{k-1}^\rho,s_{k}^\rho\right)$
			\begin{align*}
				\bigl\langle \D_z\calE&\bigl(\hat{t}_\rho(s),\hat{u}_\rho(s),\hat{z}_\rho(s)\bigr)
				-\D_z\calE\bigl(\overline{t}_\rho(s),\overline{u}_\rho(s),\overline{z}_\rho(s)\bigr),\hat{z}_\rho'(s)\bigr\rangle\\
				&=\kappa\int_\Omega \bigl(\hat{z}_\rho(s)-\overline{z}_\rho(s)\bigr)\hat{z}_\rho'(s)+ \nabla\big(\hat{z}_\rho(s)-\overline{z}_\rho(s)\big)\cdot\nabla \hat{z}_\rho'(s)\dd{x}\\
				&\hspace{50pt}+\int_\Omega \hat{z}_\rho'(s) \biggl(\hat{z}_\rho(s)\Bigl(\boldC\strain\bigl(\hat{u}_\rho(s)\bigr)\Bigr)\colon\strain\bigl(\hat{u}_\rho(s)\bigr)-\overline{z}_\rho(s)\Bigl(\boldC\strain\bigl(\overline{u}_\rho(s)\bigr)\Bigr)\colon\strain\bigl(\overline{u}_\rho(s)\bigr)\biggr)\dd{x}\\
				&=\kappa\left(s-s_{k}^\rho\right)\int_\Omega\abs{\hat{z}_\rho'(s)}^2+ \abs{\nabla\hat{z}_\rho'(s)}^2\dd{x}\\
				&\hspace{50pt}+\int_\Omega \hat{z}_\rho'(s) \biggl(\hat{z}_\rho(s)\Bigl(\left(\boldC\strain\bigl(\hat{u}_\rho(s)\bigr)\right)\colon\strain\bigl(\hat{u}_\rho(s)\bigr)-\left(\boldC\strain\bigl(\overline{u}_\rho(s)\bigr)\right)\colon\strain\bigl(\overline{u}_\rho(s)\bigr)\Bigr)\\
				&\hspace{180pt}+\bigl(\hat{z}_\rho(s)-\overline{z}_\rho(s)\bigr)\Bigl(\boldC\strain\bigl(\overline{u}_\rho(s)\bigr)\Bigr)\colon\strain\bigl(\overline{u}_\rho(s)\bigr)\biggr)\dd{x}\\
				&=\kappa \left(s-s_{k}^\rho\right) \bnorm{\hat{z}_\rho'(s)}_\calZ^2+
				\int_\Omega \hat{z}_\rho'(s) \hat{z}_\rho(s) \Bigl(\boldC\strain\bigl(\hat{u}_\rho(s)+\overline{u}_\rho(s)\bigr)\Bigr)\colon\strain\bigl(\hat{u}_\rho(s)-\overline{u}_\rho(s)\bigr)\dd{x}\\
				&\hspace{180pt}+
				\left(s-s_{k}^\rho\right)\int_\Omega \abs{\hat{z}_\rho'(s)}^2 \Bigl(\boldC\strain\bigl(\overline{u}_\rho(s)\bigr)\Bigr)\colon\strain\bigl(\overline{u}_\rho(s)\bigr)\dd{x}\\
				&=\left(s-s_{k}^\rho\right)\biggl(\kappa \bnorm{\hat{z}_\rho'(s)}_\calZ^2+
				\int_\Omega \hat{z}_\rho'(s) \hat{z}_\rho(s) \Bigl(\boldC\strain\bigl(\hat{u}_\rho(s)+\overline{u}_\rho(s)\bigr)\Bigr)\colon\strain\bigl(\hat{u}_\rho'(s)\bigr)\dd{x}\\
				&\hspace{180pt}+\int_\Omega \abs{\hat{z}_\rho'(s)}^2 \Bigl(\boldC\strain\bigl(\overline{u}_\rho(s)\bigr)\Bigr)\colon\strain\bigl(\overline{u}_\rho(s)\bigr)\dd{x}\biggr)\\
				&=\underbrace{\left(s-s_{k}^\rho\right)}_{\leq 0}\biggl(\underbrace{\kappa \bnorm{\hat{z}_\rho'(s)}_\calZ^2}_{\geq 0}+
				\int_\Omega \hat{z}_\rho'(s) \hat{z}_\rho(s) \Bigl(\boldC\strain\bigl(\hat{u}_\rho(s)+\overline{u}_\rho(s)\bigr)\Bigr)\colon\strain\bigl(\hat{u}_\rho'(s)\bigr)\dd{x}\\
				&\hspace{180pt}+\underbrace{\int_\Omega \abs{\hat{z}_\rho'(s)}^2 \Bigl(\boldC\strain\bigl(\overline{u}_\rho(s)\bigr)\Bigr)\colon\strain\bigl(\overline{u}_\rho(s)\bigr)\dd{x}}_{\geq 0}\biggr)\\
				&\leq \left(s-s_{k}^\rho\right)
				\int_\Omega \hat{z}_\rho'(s) \hat{z}_\rho(s) \Bigl(\boldC\strain\bigl(\hat{u}_\rho(s)+\overline{u}_\rho(s)\bigr)\Bigr)\colon\strain\bigl(\hat{u}_\rho'(s)\bigr)\dd{x}.
			\end{align*}
			The	right hand side can be estimated using Hölder's inequality with $\frac{1}{r_2}+\frac{2}{\tilde{p}}=1$, $\alpha\geq r_2$ and the uniform bounds from \eqref{bound}, \eqref{prop_be1} and \eqref{prop_be2} by
			\begin{align}\notag
				r.h.s.&\leq C\underbrace{\left(s_{k}^\rho-s\right)}_{\leq \rho} \underbrace{\bnorm{\hat{z}_\rho'(s)}_{\calV}}_{\leq C\text{ by \eqref{bound}}}\underbrace{\bnorm{\hat{z}_\rho(s)}_{L^\infty(\Omega)}}_{\leq C \text{ by \eqref{prop_be1}}}\Bigl(\underbrace{\bnorm{\hat{u}_\rho(s)}_{\calU^{\tilde{p}}}}_{\leq C\text{ by \eqref{bound}}}+\underbrace{\bnorm{\overline{u}_\rho(s)}_{\calU^{\tilde{p}}}}_{\leq C\text{ by \eqref{prop_be2}}}\Bigr)\underbrace{\bnorm{\hat{u}_\rho'(s)}_{\calU^{\tilde{p}}}}_{\leq C\text{ by \eqref{bound}}}\\
				&\leq C\rho\label{rhs1}
			\end{align}
			For the second term of the remainder $r_\rho$ notice first that \eqref{fixpoint:u} implies for all $k\in\N_0$ and $s\in\left(s_{k-1}^\rho,s_{k}^\rho\right)$
			\begin{align*}
				\forall w\in\calU:\quad\left\langle \D_u\calE\bigl(\overline{t}_\rho(s),\overline{u}_\rho(s),\overline{z}_\rho(s)\bigr),w\right\rangle=0\quad \text{and } w=\hat{u}_\rho'(s)\text{ is an admissible choice}.
			\end{align*}
			Using \eqref{frechetu} for $\D_u\calE$ we obtain for all $k\in\N_0$ and $s\in\left(s_{k-1}^\rho,s_{k}^\rho\right)$
			\begin{align*}
				\bigl\langle \D_u\calE\bigl(&\hat{t}_\rho(s),\hat{u}_\rho(s),\hat{z}_\rho(s)\bigr),\hat{u}_\rho'(s)\bigr\rangle
				=
				\left\langle \D_u\calE\bigl(\hat{t}_\rho(s),\hat{u}_\rho(s),\hat{z}_\rho(s)\bigr)-\D_u\calE\bigl(\overline{t}_\rho(s),\overline{u}_\rho(s),\overline{z}_\rho(s)\bigr),\hat{u}_\rho'(s)\right\rangle\\
				&=
				\int_\Omega \biggl(\Bigl(\hat{z}_\rho(s)^2+\eta\Bigr)\Bigl(\boldC\strain\bigl(\hat{u}_\rho(s)\bigr)\Bigr)-\Bigl(\overline{z}_\rho(s)^2+\eta\Bigr)\Bigl(\boldC\strain\bigl(\overline{u}_\rho(s)\bigr)\Bigr)\biggr)\colon\strain\bigl(\hat{u}_\rho'(s)\bigr)\dd{x}\\
				&\hspace{240pt}+\left\langle \ell\bigl(\overline{t}_\rho(s)\bigr)-\ell\bigl(\hat{t}_\rho(s)\bigr),\hat{u}_\rho'(s)\right\rangle\\
				&= 
				\int_\Omega \biggl(\hat{z}_\rho(s)^2\Bigl(\boldC\strain\bigl(\hat{u}_\rho(s)\bigr)\Bigr)-\overline{z}_\rho(s)^2\Bigl(\boldC\strain\bigl(\overline{u}_\rho(s)\bigr)\Bigr)\biggr)\colon\strain\big(\hat{u}_\rho'(s)\big)\dd{x}\\
				&\hspace{70pt}+ \int_\Omega\eta\Bigl(\boldC\strain\bigl(\hat{u}_\rho(s)-\overline{u}_\rho(s)\bigr)\Bigr)\colon\strain\bigl(\hat{u}_\rho'(s)\bigr)\dd{x}
				+\left\langle \ell\bigl(\overline{t}_\rho(s)\bigr)-\ell\bigl(\hat{t}_\rho(s)\bigr),\hat{u}_\rho'(s)\right\rangle\\
				&= 
				\int_\Omega \biggl(\hat{z}_\rho(s)^2\Bigl(\boldC\strain\bigl(\hat{u}_\rho(s)-\overline{u}_\rho(s)\bigr)\Bigr)+\Bigl(\hat{z}_\rho(s)^2-\overline{z}_\rho(s)^2\Bigr)\Bigl(\boldC\strain\bigl(\overline{u}_\rho(s)\bigr)\Bigr)\biggr)\colon\strain\big(\hat{u}_\rho'(s)\big)\dd{x}\\
				&\hspace{70pt}+\left(s-s_{k}^\rho\right) \int_\Omega\eta\Bigl(\boldC\strain\bigl(\hat{u}_\rho'(s)\bigr)\Bigr)\colon\strain\bigl(\hat{u}_\rho'(s)\bigr)\dd{x}
				+\left\langle \ell\bigl(\overline{t}_\rho(s)\bigr)-\ell\bigl(\hat{t}_\rho(s)\bigr),\hat{u}_\rho'(s)\right\rangle\\
				&= \underbrace{\left(s-s_{k}^\rho\right)}_{\leq 0}\biggr(
				\int_\Omega \hat{z}_\rho'(s)\Bigl(\hat{z}_\rho(s)+\overline{z}_\rho(s)\Bigr)\Bigl(\boldC\strain\bigl(\overline{u}_\rho(s)\bigr)\Bigr)\colon\strain\big(\hat{u}_\rho'(s)\big)\dd{x}\\
				&\hspace{55pt}+ \underbrace{\int_\Omega\Bigl(\hat{z}_\rho(s)^2+\eta\Bigr)\Bigl(\boldC\strain\bigl(\hat{u}_\rho'(s)\bigr)\Bigr)\colon\strain\bigl(\hat{u}_\rho'(s)\bigr)\dd{x}}_{\geq 0}\biggl)
				+\left\langle \ell\bigl(\overline{t}_\rho(s)\bigr)-\ell\bigl(\hat{t}_\rho(s)\bigr),\hat{u}_\rho'(s)\right\rangle.
			\end{align*}
			The right hand side can be estimated using \eqref{assume:external}, Hölder's inequality with $\frac{1}{r_2}+\frac{2}{\tilde{p}}=1$, $\alpha\geq r_2$ and the uniform bounds from \eqref{bound}, \eqref{prop_be1} and \eqref{prop_be2} by
			\begin{align}\notag
				r.h.s.&\leq C\underbrace{\left(s_{k}^\rho-s\right)}_{\leq \rho}\underbrace{\bnorm{\hat{z}_\rho'(s)}_\calV}_{\leq C\text{ by \eqref{bound}}}\underbrace{\bnorm{\hat{z}_\rho(s)+\overline{z}_\rho(s)}_{L^\infty(\Omega)}}_{\leq C\text{ by \eqref{prop_be1}}}\underbrace{\bnorm{\overline{u}_\rho(s)}_{\calU^{\tilde{p}}}}_{\leq C\text{ by \eqref{prop_be2}}}\underbrace{\bnorm{\hat{u}_\rho'(s)}_{\calU^{\tilde{p}}}}_{\leq C\text{ by \eqref{bound}}}
				\\
				&\hspace{150pt}+\norm{\ell\bigl(\overline{t}_\rho(s)\bigr)-\ell\bigl(\hat{t}_\rho(s)\bigr)}_{\left(\calU^{\tilde{p}}\right)^\ast}\underbrace{\bnorm{\hat{u}_\rho'(s)}_{\calU^{\tilde{p}}}}_{\leq C\text{ by \eqref{bound}}}\\
				&\leq C\rho +C\underbrace{\norm{\ell}_{C^{1,1}\left(\left[0,T\right];\left(\calU^{p'}\right)^\ast\right)}}_{\leq C\text{ by \eqref{assume:external}}}\underbrace{\abs{\overline{t}_\rho(s)-\hat{t}_\rho(s)}}_{=\left(s_{k}^\rho-s\right)\abs{\hat{t}'_\rho(s)}\leq C\rho\text{ by \eqref{bound}}}\\
				&\leq C\rho\label{rhs2}
			\end{align}
			Combining \eqref{rhs1} and \eqref{rhs2} we find
			\begin{align*}
				\exists C>0,\,\forall\rho>0,\,\forall s\in\left[-\rho,S_\rho\right]:\quad r_\rho(s)\leq C\rho.
			\end{align*}
			Hence,
			\begin{multline*}
				\exists C>0,\,\forall \rho>0,\,\forall \sigma_1,\sigma_2\in\left[-\rho,S_\rho\right]\text{ with }\sigma_1<\sigma_2:
				\int_{\sigma_1}^{\sigma_2}r_\rho(s)\dd{s}\leq C\rho\left(\sigma_2-\sigma_1\right)\leq C\rho \left(S_\rho+\rho\right),
			\end{multline*} 
		which together with the estimate for $S_\rho$ in \eqref{bound} yields \eqref{remainder}.
		\end{proof}

		\subsection{Convergence towards time-continuous parameterized BV solutions}

		According to \eqref{bound}, $S_\rho$ is bounded from above uniformly in $\rho\in\left(0,1\right)$, also we have
		\begin{align*}
			\forall\rho>0:\quad S_\rho=N(\rho)\rho\geq \sum_{k=1}^{N(\rho)}t_k^\rho-t_{k-1}^\rho=T.
		\end{align*}
		Thus there exists a subsequence $\left(\rho_n\right)_{n\in\N}$ with $\rho_n\rton 0$ and $S\geq T$ such that $S_{\rho_n}\rton S$. 

		For the proof of our subsequent main result we adapt the steps in the proof of \cite[Theorem~2.5]{Knees2018}.
		
		\begin{theorem}[Convergence result]\label{thm:convergence}
			Let the function spaces $\calZ$, $\calV$ and $\calX$, the energy functional $\calE$ and the dissipation potential $\calR$ be given according to section~\ref{sec:spaces}, section~\ref{sec:energy} with external loading $\ell$ as in \eqref{assume:external} respectively section~\ref{sec:dissipation}. Further, let $\tilde{p}$ and $\alpha$ be chosen as in \eqref{tildep} and \eqref{alpha}. Additionally, assume there are initial values $t_0=0$ and $z_0\in\calZ_{\left[0,1\right]}$ given such that there exists $\tilde{u}\in\calU^{\tilde{p}}$ with $\D_z\calE(t_0,\tilde{u},z_0)\in\calV^\ast$.
			Then there exists a sequence $\left(\rho_n\right)_{n\in\bbN}\subset\left(0,1\right)$ with $\rho_n\rton 0$  and limiting functions $\hat{t}\in W^{1,\infty}((0,S);\bbR)$, $\hat{u}\in W^{1,\infty}\big((0,S);\calU^{\tilde{p}}\big)$ and $\hat{z}\in W^{1,\infty}((0,S);\calV)\cap L^{\infty}((0,S);\calZ)\cap H^1((0,S);\calZ)$ such that the following convergences hold true for the interpolants constructed in section~\ref{sec:interpolants}
			\begin{align}
				\begin{split}
					\label{thm:cr1}
					\hat{t}_{\rho_n}&\rhupastn\hat{t}\quad\text{in }W^{1,\infty}((0,S);\bbR),\quad \hat{u}_{\rho_n}\rhupastn \hat{u}\quad \text{in } W^{1,\infty}\big((0,S);\calU^{\tilde{p}}\big)
					\\
					&\quad\text{and}\quad \hat{z}_{\rho_n}\rhupastn \hat{z}\quad \text{in } W^{1,\infty}((0,S);\calV)\text{, } L^{\infty}((0,S);\calZ)\text{ and } H^1((0,S);\calZ),
				\end{split}
			\\
				\begin{split}\label{th:cr3}
					\hat{t}_{\rho_n}&\rton\hat{t}\quad\text{ in } C([0,S];\R),\quad
					\hat{z}_{\rho_n}\rton \hat{z}\quad\text{ in } C([0,S];\calV)\\
					&\quad\text{and}\quad 
					\hat{u}_{\rho_n}\rton \hat{u}\quad\text{ in } C\bigl([0,S];L^2\bigl(\Omega,\R^d\bigr)\bigr),
				\end{split}
			\\
				\begin{split}
					\label{thm:cr2}
					\forall s\in\left[0,S\right]:\quad& \hat{t}_{\rho_n}(s),\overline{t}_{\rho_n}(s),\underline{t}_{\rho_n}(s)\rton\hat{t}(s)\quad\text{in }\R,\quad \hat{u}_{\rho_n}(s),\overline{u}_{\rho_n}(s),\underline{u}_{\rho_n}(s)\rhupn \hat{u}(s)\quad \text{in } \calU^{\tilde{p}}\\
					&\quad\text{and}\quad \hat{z}_{\rho_n}(s),\overline{z}_{\rho_n}(s),\underline{z}_{\rho_n}(s)\rhupn \hat{z}(s)\quad \text{in } \calZ.
				\end{split}
			\end{align}
			Moreover, the limiting functions $\big(\hat{t},\hat{u},\hat{z}\big)$ satisfy the initial conditions and the complementarity relations
			\begin{align}\label{thm:cr4}
				&\hat{t}(0)=t_0,\quad \quad \hat{t}(S)=T,\quad\hat{u}(0)=u_{min}(t_0,z_0)\quad\text{and}\quad \hat{z}(0)=z_0, \\\label{thm:cr5}
				\begin{split}
				\text{f.a.a. } s\in\left[0,S\right]:&\quad \hat{t}'(s)\geq 0,\quad \hat{t}'(s)+\bigl\Vert\hat{z}'(s)\bigr\Vert_\calV\leq 1,\\
					&\quad \hat{t}'(s)\dist_{\calV^\ast}\Bigl(-\D_z\calE\bigl(\hat{t}(s),\hat{u}(s),\hat{z}(s)\bigr),\partial^{\calZ}\calR(0)\Bigr)=0,
				\end{split}
			\end{align}
			the balance of linear momentum
			\begin{align}\label{balanceoflinearmomentum}
				\text{f.a.a. }s\in\left[0,S\right],\forall w\in\calU:\quad\left\langle \D_u\calE\bigl(\hat{t}(s),\hat{u}(s),\hat{z}(s)\bigr),w\right\rangle = 0
			\end{align}
			and the energy dissipation balance
			\begin{multline}\label{thm:cr6}
				\calE\bigl(\hat{t}(s_1),\hat{u}(s_1),\hat{z}(s_1)\bigr)+\int_{0}^{s_1}\calR\bigl(\hat{z}'(s)\bigr)+\bigl\Vert \hat{z}'(s)\bigr\Vert_\calV\dist_{\calV^\ast}\Bigl(-\D_z\calE\big(\hat{t}(s),\hat{u}(s),\hat{z}(s)\big),\partial^{\calZ}\calR(0)\Bigr)\dd s\\
				=\calE\bigl(\hat{t}(0),\hat{u}(0),\hat{z}(0)\bigr)+\int_{0}^{s_1}\partial_t\calE\bigl(\hat{t}(s),\hat{u}(s),\hat{z}(s)\bigr)\hat{t}'(s)\dd s\quad \forall s_1\in\left[0,S\right].
			\end{multline}
			We call a triple $\big(\hat{t},\hat{u},\hat{z}\big)\in W^{1,\infty}((0,S);\bbR)\times W^{1,\infty}\big((0,S);\calU^{\tilde{p}}\big)\times\left(W^{1,\infty}((0,S);\calV)\cap H^1((0,S);\calZ)\right)$ satisfying \eqref{thm:cr5} to \eqref{thm:cr6} \emph{$\calV$-parameterized balanced viscosity solution} of the rate-independent system driven by the functionals $\calE$ and $\calR$.
		\end{theorem}
		
		\begin{proof}
			First we will proof the convergences \eqref{thm:cr1}-\eqref{thm:cr2}. In view of the characterization of dual Bochner spaces from Lemma~\ref{thm:characterizationofdualspaces} and the subsequent Remark~\ref{remark:dualspaces} we can identify for all $\calQ\in\left\{  \R,\calU^{\tilde{p}},\calV,\calZ\right\}$ the space $L^\infty\left(\left[0,S\right];\calQ\right)$ with the space $\left(L^1\left(\left[0,S\right];\calQ^\ast\right)\right)^\ast$ as well as the space $L^2((0,S);\calZ)$ with the space $\left(L^2((0,S);\calZ^\ast)\right)^\ast$. Thus from the uniform bounds in \eqref{bound} it follows with Theorem~\ref{thm:weakastcompactness} (sequentially weak$\ast$ compactness) that there exists a sequence $\left(\rho_n\right)_{n\in\N}\subset\left(0,1\right)$ (we possibly have to extract further not relabeled subsequences) with $\rho_n\rton 0$ and limiting functions $\hat{t},\hat{s}\in L^\infty\left(\left[0,S\right];\R\right)$, $\hat{u},\hat{w}\in L^\infty\left(\left[0,S\right];\calU^{\tilde{p}}\right)$, $\hat{z}\in L^\infty\left(\left[0,S\right];\calZ\right)$, $\hat{v}\in L^\infty\left(\left[0,S\right];\calV\right)$ and $\tilde{v}\in L^2\left(\left[0,S\right];\calZ\right)$ such that
			\begin{align*}
				\hat{t}_{\rho_n}\rhupastn\hat{t}\quad\text{and}&\quad \hat{t}_{\rho_n}'\rhupastn\hat{s}\quad\text{in } L^\infty\left(\left[0,S\right];\R\right),\\
				\hat{u}_{\rho_n}\rhupastn\hat{u}\quad\text{and}&\quad \hat{u}_{\rho_n}'\rhupastn\hat{w}\quad\text{in } L^\infty\left(\left[0,S\right];\calU^{\tilde{p}}\right),\\
				\hat{z}_{\rho_n}\rhupastn\hat{z}\quad\text{in } L^\infty\left(\left[0,S\right];\calZ\right),&\quad  \hat{z}_{\rho_n}'\rhupastn\hat{v}\quad\text{in } L^\infty\left(\left[0,S\right];\calV\right)\\
				\text{and}&\quad\hat{z}_{\rho_n}'\rhupastn\tilde{v}\quad\text{in } L^2\left(\left[0,S\right];\calZ\right).
			\end{align*}
			Here $\hat{t}_{\rho_n}'$, $\hat{u}_{\rho_n}'$ and $\hat{z}_{\rho_n}'$ denote the time derivatives of $\hat{t}_{\rho_n}$, $\hat{u}_{\rho_n}$ and $\hat{z}_{\rho_n}$ in the sense of vector valued distributions. By using the compatibility of this notion of a time derivative with the notion of generalized time derivatives \cite{bochner}, we can show that $\hat{s}=\hat{t}'$, $\hat{w}=\hat{u}'$ and $\hat{v}=\hat{z}'$. 
			Also the last two weak$\ast$ convergences of $\left(\hat{z}_{\rho_n}'\right)_{n\in\N}$ imply
			\begin{align*}
				\hat{z}_{\rho_n}'\rhupastn\hat{v}\quad\text{in }L^2((0,S);\calV)\quad\text{as well as}\quad \hat{z}_{\rho_n}'\rhupastn\tilde{v}\quad\text{in }L^2((0,S);\calV), 
			\end{align*}
			such that by uniqueness of the minimizer $\tilde{v}=\hat{v}=\hat{z}'$.
			Thus, we subsume the above convergences to \eqref{thm:cr1}.
			
			Next we want to apply the Arzelà-Ascoli Theorem~\ref{thm:arzela-ascoli} to the sequences $\left(\hat{t}_{\rho_n}\right)_{n\in\N}\subset C\left(\left[0,S\right];\R\right)$, $\left(\hat{u}_{\rho_n}\right)_{n\in\N}\subset C\left(\left[0,S\right];L^2\left(\Omega,\R^d\right)\right)$ and $\left(\hat{z}_{\rho_n}\right)_{n\in\N}\subset C\left(\left[0,S\right];\calV\right)$. The uniform Lipschitz continuity from \eqref{bound} already implies the equicontinuity of these sequences. Moreover, we have for every $s\in\left[0,S\right]$ the uniform bounds
			\begin{align}\label{proof:cruniformbounds}
				\sup_{\rho\in\left(0,1\right)}\abs{\hat{t}_{\rho}(s)}\leq T,\quad \sup_{\rho\in\left(0,1\right)}\norm{\hat{u}_{\rho}(s)}_{\calU^{\tilde{p}}}\leq C\quad\text{and}\quad \sup_{\rho\in\left(0,1\right)}\norm{\hat{z}_{\rho}(s)}_\calZ\leq C,
			\end{align}
			what gives us together with the compact embeddings $\calU^{\tilde{p}}\hookrightarrow\hookrightarrow L^2\left(\Omega,\R^d\right)$ and $\calZ\hookrightarrow\hookrightarrow\calV$ the relative compactness of the sequences $\left(\hat{t}_{\rho_n}(s)\right)_{n\in\N}\subset \R$, $\left(\hat{u}_{\rho_n}(s)\right)_{n\in\N}\subset L^2\left(\Omega;\R^d\right)$ and $\left(\hat{z}_{\rho_n}(s)\right)_{n\in\N}\subset \calV$ for every $s\in\left[0,S\right]$. Thus the Arzelà-Ascoli Theorem~\ref{thm:arzela-ascoli} implies the existence of a (not relabeled) subsequence $\left(\rho_n\right)_{n\in\N}$ and limiting functions $\tilde{t}\in C\left(\left[0,S\right];\R\right)$, $\tilde{u}\in C\left(\left[0,S\right];L^2\left(\Omega,\R^d\right)\right)$ and $\tilde{z}\in C\left(\left[0,S\right];\calV\right)$ with
			\begin{align*}
				\hat{t}_{\rho_n}&\rton\tilde{t}\quad\text{ in } C([0,S];\R),\quad
				\hat{z}_{\rho_n}\rton \tilde{z}\quad\text{ in } C([0,S];\calV)\\
				\quad&\text{and}\quad 
				\hat{u}_{\rho_n}\rton \tilde{u}\quad\text{ in } C\bigl([0,S];L^2\bigl(\Omega,\R^d\bigr)\bigr).
			\end{align*}
			Additionally, due to the uniform bounds \eqref{proof:cruniformbounds}, for every $s\in\left[0,S\right]$ we can find a subsequence $\left(\rho_{n_k}\right)_{k\in\N}$ and limiting values  $\tilde{\tilde{u}}(s)\in L^2\left(\Omega,\R^d\right)$ and $\tilde{\tilde{z}}(s)\in \calV$ such that
			\begin{align*}
				\hat{u}_{\rho_{n_k}}(s)\rhupk \tilde{\tilde{u}}(s)\quad \text{in } \calU^{\tilde{p}}\quad\text{and}\quad \hat{z}_{\rho_{n_k}}(s)\rhupk \tilde{\tilde{z}}(s)\quad \text{in } \calZ.
			\end{align*}
			Due to the compact embeddings $\calU^{\tilde{p}}\hookrightarrow\hookrightarrow L^2\left(\Omega,\R^d\right)$ and $\calZ\hookrightarrow\hookrightarrow\calV$ and the uniqueness of the limiting value, it holds in particular $\tilde{\tilde{u}}(s)=\tilde{u}(s)$ and $\tilde{\tilde{z}}(s)=\tilde{z}(s)$. By the usual subsubsequence argument, we thus obtain convergence of the sequences $\left(\hat{u}_{\rho_{n}}(s)\right)_{n\in\N}$ and $\left(\hat{z}_{\rho_{n}}(s)\right)_{n\in\N}$ themself, i.e.
			\begin{align}\label{proof:crpointwiseconvergence}
				\hat{u}_{\rho_{n}}(s)\rhupn \tilde{u}(s)\quad \text{in } \calU^{\tilde{p}}\quad\text{and}\quad \hat{z}_{\rho_{n}}(s)\rhupn \tilde{z}(s)\quad \text{in } \calZ.
			\end{align}
			
			Next we will show that $\tilde{t}=\hat{t}$, $\tilde{u}=\hat{u}$ and $\tilde{z}=\hat{z}$. Let us discuss this in more detail for $\tilde{z}=\hat{z}$. From \eqref{thm:cr1} we obtain the convergence
			\begin{align*}
				\hat{z}_{\rho_n}\rhupastn\hat{z}\quad\text{in } L^\infty\left((0,S);\calV\right),
			\end{align*}
			which means in view of the characterization of the dual Bochner spaces from Lemma~\ref{thm:characterizationofdualspaces} that
			\begin{align*}
				\forall \xi\in L^1\left((0,S);\calV^\ast\right):\quad \int_{0}^{S}\langle\xi,\hat{z}_{\rho_n}\rangle\dd{s}\rton\int_{0}^{S}\langle \xi,\hat{z}\rangle\dd{s}.
			\end{align*}
			From \eqref{proof:crpointwiseconvergence} we also infer that
			\begin{align*}
				\forall \xi\in L^1\left((0,S);\calV^\ast\right),\,\text{f.a.a. }s\in (0,S):\quad
				\langle\xi(s),\hat{z}_{\rho_n}(s)\rangle\rton\langle\xi(s),\tilde{z}(s)\rangle.
			\end{align*}
			Moreover, we have
			\begin{align*}
				\forall \xi\in L^1\left((0,S);\calV^\ast\right):\quad\babs{\langle\xi(s),\hat{z}_{\rho_n}(s)\rangle}\leq \bnorm{\xi(s)}_{\calV^\ast}\bnorm{\hat{z}_{\rho_n}}_{C(\left[0,S\right];\calV)}
			\end{align*}
			and hence Lebesgue's dominated convergence theorem implies
			\begin{align*}
				\forall \xi\in L^1\left((0,S);\calV^\ast\right):\quad \int_{0}^{S}\left\langle\xi,\hat{z}_{\rho_n}\right\rangle\dd{s}\rton\int_{0}^{S}\left\langle \xi,\tilde{z}\right\rangle\dd{s},
			\end{align*}
			that is 
			\begin{align*}
				\hat{z}_{\rho_n}\rhupastn\tilde{z}\quad\text{in } L^\infty\left((0,S);\calV\right).
			\end{align*}
			Thus, by uniqueness of the limiting function, it is $\tilde{z}=\hat{z}$ in $L^\infty\left((0,S);\calV\right)$. Analogously we obtain $\tilde{t}=\hat{t}$ in $L^\infty\left((0,S);\R\right)$ and $\tilde{u}=\hat{u}$ in $L^\infty\left((0,S);\calU^{\tilde{p}}\right)$.
			Finally, the pointwise convergences in \eqref{thm:cr2} also follow for the piecewise constant interpolants. This can be seen as follows: Using the uniform bounds from \eqref{prop_be2} and the fact that $\sup\limits_{\rho>0,\,k\in\bbN_0} \babs{t_k^\rho} \leq T$, for every $s\in\left[0,S\right]$ we can find a subsequence $\left(\rho_{n_k}\right)_{k\in\N}$ and limiting values $\underline{t}(s),\overline{t}(s)\in \R$, $\underline{u}(s),\overline{u}(s)\in \calU^{\tilde{p}}$ and $\underline{z}(s),\overline{z}(s)\in \calZ$ such that
			\begin{align*}
				&\underline{t}_{\rho_{n_k}}(s)\rtok \underline{t}(s)\quad \text{in } \R,\quad
				\underline{u}_{\rho_{n_k}}(s)\rhupk \underline{u}(s)\quad \text{in } \calU^{\tilde{p}}\quad\text{and}\quad\underline{z}_{\rho_{n_k}}(s)\rhupk \underline{z}(s)\quad \text{in } \calZ\\
				\text{and}\quad
				&\overline{t}_{\rho_{n_k}}(s)\rtok \overline{t}(s)\quad \text{in } \R,\quad
				\overline{u}_{\rho_{n_k}}(s)\rhupk \overline{u}(s)\quad \text{in } \calU^{\tilde{p}}\quad\text{and}\quad\overline{z}_{\rho_{n_k}}(s)\rhupk \overline{z}(s)\quad \text{in } \calZ.
			\end{align*}
			Due to the compact embedding 
			  $\calZ\hookrightarrow\hookrightarrow\calV$ this in particular yields
			 \begin{align*}
			 	&
			 	\underline{z}_{\rho_{n_k}}(s)\rtok \underline{z}(s)\quad \text{in } \calV\quad\text{and}\quad\overline{z}_{\rho_{n_k}}(s)\rtok \overline{z}(s)\quad \text{in } \calV
			 \end{align*}
			Moreover, by the choice of the interpolants and by the uniform bounds from \eqref{bound} we obtain
			\begin{align*}
				\exists C>0,\,\forall s\in\left(s_{k-1}^\rho,s_k^\rho\right):&\quad \bnorm{\hat{z}_{\rho}(s)-\underline{z}_{\rho}(s)}_\calV=\babs{s-s_{k-1}^\rho}\,\bnorm{\hat{z}_{\rho}'(s)}_\calV\leq \rho\, C\\
				&\quad\text{and}\quad \bnorm{\hat{z}_{\rho}(s)-\overline{z}_{\rho}(s)}_\calV=\babs{s-s_k^\rho}\,\bnorm{\hat{z}_{\rho}'(s)}_\calV\leq \rho\, C,
			\end{align*}
			and hence $\rho_n\rton 0$ implies that
			\begin{align*}
				\forall s\in\left[0,S\right]:\quad \bnorm{\hat{z}_{\rho_n}(s)-\underline{z}_{\rho_n}(s)}_\calV\rton 0\quad\text{and}\quad \bnorm{\hat{z}_{\rho_n}(s)-\overline{z}_{\rho_n}(s)}_\calV\rton 0.
			\end{align*}
			Analogously 
			\begin{align*}
				\forall s\in\left[0,S\right]:&\quad \babs{\hat{t}_{\rho_n}(s)-\underline{t}_{\rho_n}(s)}\rton 0\quad\text{and}\quad \babs{\hat{t}_{\rho_n}(s)-\overline{t}_{\rho_n}(s)}\rton 0\\
				\quad\text{and}&\quad \bnorm{\hat{u}_{\rho_n}(s)-\overline{u}_{\rho_n}(s)}_{\calU^{\tilde{p}}}\rton 0\quad\text{and}\quad \bnorm{\hat{u}_{\rho_n}(s)-\overline{u}_{\rho_n}(s)}_{\calU^{\tilde{p}}}\rton 0\\
				\quad\text{and particularly}&\quad \hat{u}_{\rho_n}(s)-\underline{u}_{\rho_n}(s)\rhupn  0\quad \text{in }\calU^{\tilde{p}}\quad\text{and}\quad \hat{u}_{\rho_n}(s)-\overline{u}_{\rho_n}(s)\rhupn  0\quad \text{in }\calU^{\tilde{p}}.
			\end{align*}
			Now $\underline{z}_{\rho_{n_k}}(s)\rtok \underline{z}(s)$ in $\calV$, $\hat{z}_{\rho_{n}}(s)\rton \hat{z}(s)$ in $\calV$ and $\bnorm{\hat{z}_{\rho_n}(s)-\underline{z}_{\rho_n}(s)}_\calV\rton 0$ imply that $\underline{z}(s)=\hat{z}(s)$ and again by the usual subsubsequence argument we obtain convergence of the sequence $\left(\underline{z}_{\rho_{n}}(s)\right)_{n\in\N}$ itself. Analogously the convergences of the remaining constant interpolants follow.
			
			Next we will prove \eqref{thm:cr4}. For all $\rho>0$ we have $\hat{t}_\rho(0)=t_0^\rho=t_0$, so by $\hat{t}_{\rho_n}(0)\rton\hat{t}(0)$ we infer that $\hat{t}(0)=t_0$. 
			By \eqref{reachedT} we know that for all $\rho>0$ it is $\hat{t}_\rho(S_\rho)=T$. Moreover, for every $\rho>0$ the function $\hat{t}_\rho$ is Lipschitz continuous with Lipschitz constant $1$ and it holds $S_{\rho_n}\rton S$. Thus
			\begin{align*}
				\babs{\hat{t}_{\rho_n}(S)-T}=\babs{\hat{t}_{\rho_n}(S)-\hat{t}_{\rho_n}(S_{\rho_n})}\leq \babs{S-S_{\rho_n}}\rton 0
			\end{align*}
			and by $\hat{t}_{\rho_n}(S)\rton\hat{t}(S)$ it eventually follows that $\hat{t}(S)=T$.
			For all $\rho>0$ we have $\hat{z}_\rho(0)=z_0^\rho$, which is the E\&M iterate determined after the first AM loop has run through, and $\bnorm{z_0^\rho-z_0}_\calV=\bnorm{z_0^\rho-z_{-1}^\rho}_\calV\leq\rho$. Thus it follows with $\hat{z}_{\rho_n}(0)\rton\hat{z}(0)$ in $\calV$ that
			\begin{align*}
				\bnorm{\hat{z}(0)-z_0}_\calV\leq \bnorm{\hat{z}(0)-\hat{z}_{\rho_n}(0)}_\calV+\bnorm{z_0^{\rho_n}-z_0}_{\calV}\rton 0,\quad \text{i.e. } \hat{z}(0)=z_0.
			\end{align*}
			By the choice of the interpolants and by Proposition~\ref{lemma_convergentsubsequence} we have $\hat{u}_\rho(0)=u_0^\rho=u_{min}(t_0^\rho,z_0^\rho)$ for all $\rho>0$. Thus Lemma~\ref{lemmaA1} implies that
			\begin{align*}
				\bnorm{\hat{u}_{\rho_n}(0)-u_{min}(t_0,z_0)}_{\calU^{\tilde{p}}}\leq C\underbrace{\babs{t_0^{\rho_n}-t_0}}_{=0}+C\underbrace{\bnorm{z_0^{\rho_n}-z_0}_{L^{r_1}(\Omega)}}_{\leq\rho_n}\leq C\rho_n\rton 0
			\end{align*}
			and hence, by $\hat{u}_{\rho_n}(0)\rhupn \hat{u}(0)$ in $\calU^{\tilde{p}}$, it is $\hat{u}(0)=u_{min}(t_0,z_0)$.
		
			Next we will proof \eqref{thm:cr5}. By \eqref{normalization} we have $\hat{t}_\rho'(s)\geq 0$ for all $\rho>0$ and f.a.a. $s\in\left[0,S\right]$. Thus $\hat{t}_{\rho_n}'\rhupastn \hat{t}'$ in $L^\infty\left((0,S);\R\right)$ implies that
			\begin{align*}
				\forall \phi\in C_0^\infty\left((0,S);\R\right)\text{ with }\phi\geq 0:\quad 0\leq \int_{0}^{S}\phi(s)\hat{t}_{\rho_n}'(s)\rton\int_{0}^{S}\phi(s)\hat{t}'(s)\dd{s}.
			\end{align*}
			Thus we obtain by the fundamental lemma of the calculus of variations (cf. \cite[Satz~5.1]{Dobro2010}) $\hat{t}'(s)\geq 0$ f.a.a. $s\in\left[0,S\right]$, and hence also $\hat{t}'(s+\rho)+\norm{\hat{z}_\rho'(s)}_\calV\geq 0$ for all $\rho>0$ and f.a.a. $s\in\left[0,S\right]$.
			Further we have $\hat{t}_\rho'(s+\rho)+\norm{\hat{z}_\rho'(s)}_\calV\leq 1$ for all $\rho>0$ and f.a.a. $s\in\left[0,S\right]$. Consider the set
			\begin{align*}
				M\coloneq \Set{\left(\tau,v\right)\in L^2((0,S);\R)\times L^2((0,S);\calV)}{\tau(s)+\norm{v(s)}_\calV\leq 1\text{ f.a.a. } s\in (0,S)}.
			\end{align*}
			It is clearly convex and by contradiction we can proof, that it is closed in $L^2((0,S);\R)\times L^2((0,S);\calV)$: Let $\left(\tau_n,v_n\right)_{n\in\N}\subset M$ with
			\begin{align*}
				\tau_n\rton \tau\text{ in } L^2((0,S);\R)\quad \text{and}\quad v_n\rton v\text{ in } L^2((0,S);\calV).
			\end{align*}
			Assume there exists $I\subset(0,S)$ with $\abs{I}>0$ such that
			\begin{align*}
				\tau(s)+\norm{v(s)}_{\calV}>1\quad\forall s\in I.
			\end{align*}
			Then there exists $c>0$ such that
			\begin{align*}
				0<c&\leq
				\int_{I}\underbrace{\left(\tau(s)+\norm{v(s)}_{\calV}\right)}_{>1}-\underbrace{\left(\tau_n(s)+\norm{v_n(s)}_{\calV}\right)}_{\leq 1}\dd{s}\leq \int_{I}\babs{\tau(s)-\tau_n(s)}+\babs{\norm{v(s)}_{\calV}-\norm{v_n(s)}_{\calV}}\dd{s}\\
				&\leq \int_{0}^{S}\babs{\tau(s)-\tau_n(s)}+\norm{v(s)-v_n(s)}_{\calV}\dd{s}\rton 0,
			\end{align*}
			a contradiction. 
			Eventually, we infer that $M$ is weakly closed in $L^2((0,S);\R)\times L^2((0,S);\calV)$.
			The weak$\ast$ convergence $\hat{z}_{\rho_n}'\rhupastn\hat{z}'$ in $L^\infty((0,S);\calV)$ in particular implies that $\hat{z}_{\rho_n}'\rhupn\hat{z}'$ in $L^2((0,S);\calV)$. 
			Now let $\phi\in C_0^\infty\left((0,S);\R\right)$ be arbitrary and continuously extend it on $\left(-\infty,S\right)$ by setting $\phi(s)\coloneq 0$ for $s\leq 0$. In particular $\phi$ is uniformly continuous due to the Heine-Cantor theorem. Together with the uniform bound $\sup_{\rho>0,\,s\in\left[0,S\right]}\abs{\hat{t}_\rho'(s)}\leq 1$ this yields
			\begin{align*}
				\abs{\int_{0}^{S}\phi(s)\hat{t}_{\rho_n}'(s)\dd{s}-\int_{\rho_n}^{S+\rho_n}\phi(s-\rho_n)\hat{t}_{\rho_n}'(s)\dd{s}}\rton 0.
			\end{align*}
			Therefore, by the weak$\ast$ convergence $\hat{t}_{\rho_n}'\rhupastn\hat{t}'$ in $L^\infty\left(\left[0,S\right];\R\right)$, we obtain for all $\phi\in C_0^\infty\left((0,S);\R\right)$
			\begin{align}\label{rhoshift}
 				\begin{split}
 					\Biggl\lvert\int_{0}^{S}&\phi(s)\left(\hat{t}'(s)-\hat{t}_{\rho_n}'(s+\rho_n)\right)\dd{s}\Biggr\rvert\\
 					&\leq
 					\abs{\int_{0}^{S}\phi(s)\left(\hat{t}'(s)-\hat{t}_{\rho_n}'(s)\right)\dd{s}}
 					+\abs{\int_{0}^{S}\phi(s)\hat{t}_{\rho_n}'(s)\dd{s}-\int_{\rho_n}^{S+\rho_n}\phi(s-\rho_n)\hat{t}_{\rho_n}'(s)\dd{s}}\rton 0.
 				\end{split}
			\end{align}
			For all $\phi,\psi\in L^1((0,S);\R)$ and all $\tau\in L^\infty((0,S);\R)$ it holds 
			\begin{align*}
				\abs{\int_{0}^{S}\left(\phi(s)-\psi(s)\right)\tau(s)\dd{s}}\leq \norm{\phi-\psi}_{L^1((0,S);\R)}\norm{\tau}_{L^\infty((0,S);\R)}.
			\end{align*}
			Moreover, $C_0^\infty((0,S);\R)$ is densely embedded in $L^2((0,S);\R)$ and we have the uniform bound $\sup_{\rho>0}\norm{\hat{t}_\rho'}_{L^\infty((0,S);\R)}\leq 1$. Thus \eqref{rhoshift} stays valid for all $\phi\in L^2((0,S);\R)$, i.e. $\hat{t}_{\rho_n}'(\cdot+\rho_n)\rhupn \hat{t}'$ in $L^2((0,S);\R)$.
			Altogether we conclude the second relation in \eqref{thm:cr5}.
			
			Next we prove the last relation in \eqref{thm:cr5}. Therefore we use the lower semicontinuity properties from Lemma~\ref{Knees18:lemmaB2} with $\delta_n=\dist_{\calV^\ast}(-\D_z\calE\big(\underline{t}_{\rho_n}(\cdot),\underline{u}_{\rho_n}(\cdot),\underline{z}_{\rho_n}(\cdot)\big),\partial^\calZ\calR(0))$ and $\tau_n=\hat{t}_{\rho_n}'$. In advance we have to check the assumptions thereof using Lemma~\ref{lemma:lowersemicon}.
			By \eqref{prop:fi2}, we have
			\begin{align*}
				\exists C>0,\, \forall \rho>0,\,\forall s\in\left[0,S\right]:\quad \dist_{\calV^\ast}(-\D_z\calE\big(\underline{t}_{\rho_n}(s),\underline{u}_{\rho_n}(s),\underline{z}_{\rho_n}(s)\big),\partial^\calZ\calR(0))\leq C.
			\end{align*}
			If we additionally show that 
			\begin{multline}\label{liminfDzE}
				\text{f.a.a. } s\in\left[0,S\right],\,\forall v\in\calZ\text{ with } \calR(v)<\infty:\quad \\
				\liminf_{n\to\infty}\left\langle \D_z\calE\big(\underline{t}_{\rho_n}(s),\underline{u}_{\rho_n}(s),\underline{z}_{\rho_n}(s)\big),-v\right\rangle
				\geq \left\langle \D_z\calE\big(\hat{t}(s),\hat{u}(s),\hat{z}(s)\big),-v\right\rangle,
			\end{multline}
			then Lemma~\ref{lemma:lowersemicon} gives us
			\begin{multline}\label{proof:crdist}
				\text{f.a.a. } s\in\left[0,S\right]:\quad\dist_{\calV^\ast}\bigl(-\D_z\calE\bigl(\hat{t}(s),\hat{u}(s),\hat{z}(s)\bigr),\partial^\calZ\calR(0)\bigr)\\
				\leq \liminf_{n\to\infty} \dist_{\calV^\ast}\bigl(-\D_z\calE\bigl(\underline{t}_{\rho_n}(s),\underline{u}_{\rho_n}(s),\underline{z}_{\rho_n}(s)\bigr),\partial^\calZ\calR(0)\bigr).
			\end{multline}
			For this purpose let $s\in\left[0,S\right]$ and $v\in\calZ$ with $\calR(v)<\infty$ be arbitrary.  By \eqref{frechetz} we have
			\begin{multline}\label{secondtermofthis}
				\langle \D_z\calE\big(\underline{t}_{\rho_n}(s),\underline{u}_{\rho_n}(s),\underline{z}_{\rho_n}(s)\big),-v\rangle\\
				=
				-\kappa \int_{\Omega} \underline{z}_{\rho_n}(s) v+\nabla \underline{z}_{\rho_n}(s)\cdot\nabla v\dd{x}-\int_{\Omega} \underline{z}_{\rho_n}(s) v\left(\boldC\strain(\underline{u}_{\rho_n}(s))\right)\colon\strain(\underline{u}_{\rho_n}(s))\dd{x}
			\end{multline}
			Due to $\underline{z}_{\rho_n}(s)\rhupn \hat{z}(s)$ in $\calZ$ we obtain for the first term on the right hand side
			\begin{align*}
				-\kappa \int_{\Omega} \underline{z}_{\rho_n}(s) v+\nabla \underline{z}_{\rho_n}(s)\cdot\nabla v\dd{x}\rton -\kappa \int_{\Omega} \hat{z}(s) v+\nabla \hat{z}(s)\cdot\nabla v\dd{x}.
			\end{align*}
			For the second term notice first that by Hölder's inequality with $\frac{1}{r_2}+\frac{2}{\tilde{p}}=1$ and the compact embedding $\calZ\hookrightarrow\hookrightarrow L^{2r_2}(\Omega)$ we have
			\begin{align*}
				\Big\vert\int_{\Omega} \left(\underline{z}_{\rho_n}(s)-\hat{z}(s)\right) (-v)\left(\boldC\strain(\underline{u}_{\rho_n}(s))\right)&\colon\strain(\underline{u}_{\rho_n}(s))\dd{x}\Big\vert\\
				&\leq 
				\underbrace{\big\Vert\underline{z}_{\rho_n}(s)-\hat{z}(s)\big\Vert_{L^{2r_2}(\Omega)}}_{\rton 0}\underbrace{\norm{v}_{L^{2r_2}(\Omega)}}_{\leq C\norm{v}_{\calZ}}\underbrace{\big\Vert \underline{u}_{\rho_n}(s)\big\Vert_{\calU^{\tilde{p}}}^2}_{\leq C}\\
				&\rton 0.
			\end{align*}
			Then consider for fixed $z,v\in\calZ$ with $\calR(v)<\infty$ and $z\geq 0$ a.e. in $\Omega$ the mapping 
			\begin{align*}
				F_{z,v}:\calU^{\tilde{p}}\to\bbR,\quad F_{z,v}(u)\coloneq \int_\Omega z(-v) \left(\boldC\strain(u)\right)\colon\strain(u)\dd{x}.
			\end{align*}
			It is twice Frèchet differentiable with
			\begin{align*}
				\forall w,u\in\calU^{\tilde{p}}:\quad \left\langle \left(\D_u^2F_{z,v}(u)\right)(w),w\right\rangle =2\int_{\Omega} z(-v) \left(\boldC\strain(w)\right)\colon\strain(w)\dd{x}.
			\end{align*}
			By assumption \eqref{assume:elastictensor} on the elastic tensor and since we assumed $z(-v)\geq 0$ a.e. in $\Omega$, $\D_u^2F_{z,v}(u)$ is positive semidefinite in every $u\in\calU^{\tilde{p}}$. Therefore $F_{z,v}$ is convex. 
			Moreover, $F_{z,v}$ is continuous on $\calU^{\tilde{p}}$. This can be seen by using again Hölder's inequality with $\frac{1}{r_2}+\frac{2}{\tilde{p}}=1$ and the compact embedding $\calZ\hookrightarrow\hookrightarrow L^{2r_2}(\Omega)$. Thus $F_{z,v}$ is weakly lower semicontinuous on $\calU^{\tilde{p}}$. Eventually, by $\underline{u}_{\rho_n}(s)\rhupn \hat{u}(s)$ in $\calU^{\tilde{p}}$, the limes inferior of the second term in \eqref{secondtermofthis} can be estimated as
			\begin{align*}
				\liminf_{n\to\infty}\int_{\Omega} &\underline{z}_{\rho_n}(s) (-v)\left(\boldC\strain(\underline{u}_{\rho_n}(s))\right)\colon\strain(\underline{u}_{\rho_n}(s))\dd{x}\\
				&=\liminf_{n\to\infty}\left(\int_{\Omega} \left(\underline{z}_{\rho_n}(s)-\hat{z}(s)\right) (-v)\left(\boldC\strain(\underline{u}_{\rho_n}(s))\right)\colon\strain(\underline{u}_{\rho_n}(s))\dd{x}+F_{\hat{z}(s),v}\bigl(\underline{u}_{\rho_n}(s)\bigr)\right)\\
				&=\liminf_{n\to\infty}F_{\hat{z}(s),v}\bigl(\underline{u}_{\rho_n}(s)\bigr)\\
				&\geq F_{\hat{z}(s),v}\bigl(\hat{u}(s)\bigr)=\int_\Omega \hat{z}(s)(-v) \bigl(\boldC\strain(\hat{u}(s))\bigr)\colon\strain(\hat{u}(s))\dd{x}.
			\end{align*}
			Altogether this shows \eqref{liminfDzE} and thus \eqref{proof:crdist}.
			
			\noindent
			The weak$\ast$ convergence $\hat{t}_{\rho_n}'\rhupastn\hat{t}'$ in $L^\infty((0,S);\R)\cong\left(L^1((0,S);\R)\right)^\ast$ in particular implies
			$\hat{t}_{\rho_n}'\rhupn\hat{t}'$ weakly in $L^1((0,S);\R)$.
			By Lemma~\ref{Knees18:lemmaB2} we therefore obtain together with \eqref{prop:complementarity}
			\begin{align*}
				0&\leq\int_{0}^{S}\hat{t}'(s)\dist_{\calV^\ast}\bigl(-\D_z\calE\bigl(\hat{t}(s),\hat{u}(s),\hat{z}(s)\bigr),\partial^\calZ\calR(0)\bigr) \dd{s}\\
				&\leq \liminf_{n\to\infty}\int_{0}^{S}\hat{t}_{\rho_n}'(s) \dist_{\calV^\ast}\bigl(-\D_z\calE\bigl(\underline{t}_{\rho_n}(s),\underline{u}_{\rho_n}(s),\underline{z}_{\rho_n}(s)\bigr),\partial^\calZ\calR(0)\bigr) \dd{s}\overset{\text{\eqref{prop:complementarity}}}{=}0.
			\end{align*}
			Eventually, due to the non-negativity of the integrand in the first line for almost all $s\in(0,S)$, this gives us the last relation in \eqref{thm:cr5}.

			Next we show that the balance of linear momentum \eqref{balanceoflinearmomentum} is satisfied. By \eqref{discrete:displacementEL} it is
			\begin{align*}
				\text{f.a.a. } s\in\left[0,S\right],\,\forall w\in\calU :\quad\left\langle \D_u\calE\bigl(\overline{t}_{\rho_n}(s),\overline{u}_{\rho_n}(s),\overline{z}_{\rho_n}(s)\bigr),w\right\rangle= 0.
			\end{align*}
			Employing \eqref{frechetu}, the pointwise convergences from \eqref{thm:cr2}
			\begin{align*}
				\forall s\in\left[0,S\right]:\quad \overline{t}_{\rho_n}(s)\rton\hat{t}(s),\quad \overline{u}_{\rho_n}(s)\rhupn \hat{u}(s)\text{ in }\calU^{\tilde{p}},\quad \overline{z}_{\rho_n}(s)\rton\hat{z}(s)\text{ in }L^{4r_2}(\Omega)
			\end{align*}
			and the assumptions \eqref{assume:external} on the external loading $\ell$, we obtain f.a.a. $s\in\left[0,S\right]$ and all $w\in\calU$
			\begin{align*}
				0&=\left\langle \D_u\calE\bigl(\overline{t}_{\rho_n}(s),\overline{u}_{\rho_n}(s),\overline{z}_{\rho_n}(s)\bigr),w\right\rangle\\
				&=\int_{\Omega}\bigl(\bigl(\overline{z}_{\rho_n}(s)\bigr)^2+\eta\bigr)\left(\normalfont\boldC\strain(\overline{u}_{\rho_n}(s))\right)\colon\strain(w)\dd{x}-\langle \ell(\overline{t}_{\rho_n}(s)),w\rangle\\
				&\rton\int_{\Omega}\bigl(\bigl(\hat{z}(s)\bigr)^2+\eta\bigr)\left(\normalfont\boldC\strain(\hat{u}(s))\right)\colon\strain(w)\dd{x}-\langle \ell(\hat{t}(s)),w\rangle=\left\langle \D_u\calE\big(\hat{t}(s),\hat{u}(s),\hat{z}(s)\big),w\right\rangle,
			\end{align*}
			i.e. \eqref{balanceoflinearmomentum}.				
				
			For the proof of \eqref{thm:cr6} we pass to the limit in the discrete energy dissipation balance \eqref{energydissipationbalance}.
			By Lemma~\ref{lemma:continuity} the energy functional $\calE$ is weakly lower semicontinuous on $\left[0,S\right]\times\calU^{\tilde{p}}\times\calZ$. Thus the pointwise weak convergences from \eqref{thm:cr2} imply
			\begin{align}\label{claim1}
				\forall s\in\left[0,S\right]:\quad \calE\big(\hat{t}(s),\hat{u}(s),\hat{z}(s)\big)\leq \liminf_{n\to\infty}\calE\big(\hat{t}_{\rho_n}(s),\hat{u}_{\rho_n}(s),\hat{z}_{\rho_n}(s)\big).
			\end{align}
			
			As in the proof of the continuity properties of $\calE$ (Lemma~\ref{lemma:continuity}), by the assumptions \eqref{assume:external} on the external loading $\ell$ and the pointwise convergences from \eqref{thm:cr2}, the convergence principle from \cite[p.~58, Lemma~0.3(ii)]{Ruzicka2004}) yields 
			\begin{align*}
				\forall s\in\left[0,S\right]:\quad\bigl\langle\dot{\ell}\bigl(\hat{t}_{\rho_n}(s)\bigr),\hat{u}_{\rho_n}(s)\bigr\rangle\rton\bigl\langle\dot{\ell}\bigl(\hat{t}(s)\bigr),\hat{u}(s)\bigr\rangle,
			\end{align*}
			that is the pointwise convergence
			\begin{align*}
				\forall s\in\left[0,S\right]:\quad \partial_t\calE\bigl(\hat{t}_{\rho_n}(s),\hat{u}_{\rho_n}(s),\hat{z}_{\rho_n}(s)\bigr)\rton \partial_t\calE\bigl(\hat{t}(s),\hat{u}(s),\hat{z}(s)\bigr).
			\end{align*}
			Furthermore, by \eqref{assume:external} and \eqref{bound} it holds
			\begin{align*}
				\exists C>0,\,\forall s\in\left[0,S\right]:\quad\abs{\partial_t\calE\bigl(\hat{t}_{\rho_n}(s),\hat{u}_{\rho_n}(s),\hat{z}_{\rho_n}(s)\bigr)}&=\abs{\bigl\langle\dot{\ell}\bigl(\hat{t}_{\rho_n}(s)\bigr),\hat{u}_{\rho_n}(s)\bigr\rangle}\\
				&\leq \norm{\ell}_{C^{1,1}([0,T];\left(\calU^{\tilde{p}}\right)^\ast)}\bnorm{\hat{u}_{\rho_n}(s)}_{\calU^{\tilde{p}}}\leq C.
			\end{align*}
			Thus Lebesgue's dominated convergence theorem implies 
			\begin{align*}
				\partial_t\calE\bigl(\hat{t}_{\rho_n}(\cdot),\hat{u}_{\rho_n}(\cdot),\hat{z}_{\rho_n}(\cdot)\bigr)\rton \partial_t\calE\bigl(\hat{t}(\cdot),\hat{u}(\cdot),\hat{z}(\cdot)\bigr)\quad\text{in }L^1\left((0,s_1);\R\right)
			\end{align*}
			and hence we obtain by $\hat{t}_{\rho_n}'\rhupastn\hat{t}'$ in $L^\infty\left((0,s_1);\R\right)\cong \left(L^1\left((0,s_1);\R\right)\right)^\ast$ (and by using a convergence principle similar to the above one but for weak$\ast$ convergence) that
			\begin{align*}
				\int_{0}^{s_1}\partial_t\calE\bigl(\hat{t}_{\rho_n}(s),\hat{u}_{\rho_n}(s),\hat{z}_{\rho_n}(s)\bigr)\hat{t}_{\rho_n}'(s)\dd{s}\rton \int_{0}^{s_1}\partial_t\calE\bigl(\hat{t}(s),\hat{u}(s),\hat{z}(s)\bigr)\hat{t}'(s)\dd{s}.
			\end{align*}
			By the choice of the interpolants and by \eqref{upperboundenergyEM} we have
			\begin{align*}
				\liminf_{n\to\infty}\calE\bigl(\hat{t}_{\rho_n}(0),\hat{u}_{\rho_n}(0),\hat{z}_{\rho_n}(0)\bigr)&=\liminf_{n\to\infty}\calE\bigl(t_0,\hat{u}_{0}^{\rho_n},\hat{z}_{0}^{\rho_n}\bigr)\\
				&\leq \calE\bigl(t_0,u_{min}(t_0,z_0),z_0\bigr)\overset{\text{\eqref{thm:cr4}}}{=}\calE\bigl(\hat{t}(0),\hat{u}(0),\hat{z}(0)\bigr).
			\end{align*}
		
			Next, since $\calR$ is convex and continuous on $\calV$ and since $\hat{z}_{\rho_n}'\rhupastn\hat{z}'$ in $L^\infty((0,S);\calV)$, we obtain by the generalized version of Ioffe's theorem from \cite[Theorem~21]{Valadier1990}
			\begin{align*}
				\int_{0}^{s_1}\calR(\hat{z}'(s))\dd{s}\leq\liminf_{n\to\infty} \int_{0}^{s_1}\calR(\hat{z}_{\rho_n}'(s))\dd{s}.
			\end{align*}
		Lastly, by the weak$\ast$ convergence $\hat{z}_{\rho_n}'\rhupastn \hat{z}'$ in $L^\infty((0,S);\calV)$ and \eqref{proof:crdist}, Lemma~\ref{Knees18:lemmaB1} implies
		\begin{multline*}
			\int_{0}^{s_1}\bnorm{\hat{z}'(s)}_{\calV}\dist_{\calV^\ast}\Bigl(-\D_z\calE\bigl(\hat{t}(s),\hat{u}(s),\hat{z}(s)\bigr),\partial\calR(0)\Bigr)\dd{s}\\
			\leq \liminf_{n\to\infty}\int_{0}^{s_1} \bnorm{\hat{z}_{\rho_n}'(s)}_{\calV}\dist_{\calV^\ast}\bigl(-\D_z\calE\bigl(\underline{t}_{\rho_n}(s),\underline{u}_{\rho_n}(s),\underline{z}_{\rho_n}(s)\bigr),\partial^\calZ\calR(0)\bigr)\dd{s}.
		\end{multline*}
		In summary we can take the limit inferior in the discrete energy dissipation balance \eqref{energydissipationbalance}, this gives
		\begin{align*}
			\calE&\bigl(\hat{t}(s_1),\hat{u}(s_1),\hat{z}(s_1)\bigr)
			+\int_{0}^{s_1}\calR\bigl(\hat{z}'(s)\bigr)+\bnorm{\hat{z}'(s)}_{\calV}\dist_{\calV^\ast}\Bigl(-\D_z\calE\bigl(\hat{t}(s),\hat{u}(s),\hat{z}(s)\bigr),\partial\calR(0)\Bigr)\dd{s}\\
			&\leq \liminf_{n\to\infty} \biggl(
			\calE\bigl(\hat{t}_{\rho_n}(s_1),\hat{u}_{\rho_n}(s_1),\hat{z}_{\rho_n}(s_1)\bigr)\\
			&\hspace{60pt}+\int_{0}^{s_1}\calR\big(\hat{z}_{\rho_n}'(s)\big)+\bnorm{\hat{z}_{\rho_n}'(s)}_\calV\dist_{\calV^\ast}\Bigl(-\D_z\calE\bigl(\underline{t}_{\rho_n}(s),\underline{u}_{\rho_n}(s),\underline{z}_{\rho_n}(s)\bigr),\partial\calR(0)\Bigr)\dd{s}
			\biggr)\\
			&\hspace{-6pt}\overset{\text{\eqref{energydissipationbalance}}}{=}\liminf_{n\to\infty}
			\biggl(
			\calE\bigl(\hat{t}_{\rho_n}(0),\hat{u}_{\rho_n}(0),\hat{z}_{\rho_n}(0)\bigr)
			+\int_{0}^{s_1}\partial_t\calE\bigl(\hat{t}_{\rho_n}(s),\hat{u}_{\rho_n}(s),\hat{z}_{\rho_n}(s)\bigr)\hat{t}_{\rho_n}'(s)\dd{s}
			+\underbrace{\int_{0}^{s_1}r_{\rho_n}(s)\dd{s}}_{\leq C\rho_n\rton 0}
			\biggr)\\
			&\leq \calE(t_0,u_0,z_0)+\int_{0}^{s_1}\partial_t\calE\bigl(\hat{t}(s),\hat{u}(s),\hat{z}(s)\bigr)\hat{t}'(s)\dd{s}.
		\end{align*}
	
			Finally we need to show that we indeed have an equality. This is done analogously to \cite[Lemma~2.4.6]{Sievers2020} resp. \cite[Lemma~6.6]{KRZ} in the subsequent Lemma~\ref{lemma:equality}.
	\end{proof}
			\begin{lemma}[Equality]\label{lemma:equality}
				Let $\big(\hat{t},\hat{u},\hat{z}\big)$ be a triple with  $\hat{t}\in W^{1,\infty}((0,S);\bbR)$, $\hat{u}\in W^{1,\infty}\big((0,S);\calU^{\tilde{p}}\big)$ and $\hat{z}\in W^{1,\infty}((0,S);\calV)\cap H^{1}((0,S);\calZ)$ fulfilling the initial conditions \eqref{thm:cr4}, the complementarity relations \eqref{thm:cr5} and the balance of linear momentum \eqref{balanceoflinearmomentum}.
				Then $\big(\hat{t},\hat{u},\hat{z}\big)$ is a $\calV$-parameterized balanced viscosity solution if and only if the following energy dissipation inequality is fulfilled
				\begin{multline}\label{energyinequality}
					\calE\bigl(\hat{t}(s_1),\hat{u}(s_1),\hat{z}(s_1)\bigr)+\int_{0}^{s_1}\calR\bigl(\hat{z}'(s)\bigr)+\bigl\Vert\hat{z}'(s)\bigr\Vert_\calV\dist_{\calV^\ast}\Bigl(-\D_z\calE\bigl(\hat{t}(s),\hat{u}(s),\hat{z}(s)\bigr),\partial^{\calZ}\calR(0)\Bigr)\dd s\\
					\leq\calE\bigl(\hat{t}(0),\hat{u}(0),\hat{z}(0)\bigr)+\int_{0}^{s_1}\partial_t\calE\bigl(\hat{t}(s),\hat{u}(s),\hat{z}(s)\bigr)\hat{t}'(s)\dd s\quad \forall s_1\in\left[0,S\right].
				\end{multline}
			\end{lemma}
			
			\begin{proof}
				The first implication is trivial. Now let only the inequality be satisfied. 
				By a parametrized chain rule (cf. \cite[Lemma~A.2.5]{Sievers2020}) it holds
				\begin{multline*}
					\text{f.a.a. }s\in\left[0,S\right]:\quad
					\frac{\dd{}}{\dd{s}}\Bigl(\calE\bigl(\hat{t}(s),\hat{u}(s),\hat{z}(s)\bigr)\Bigr)
					=
					\partial_t\calE\bigr(\hat{t}(s),\hat{u}(s),\hat{z}(s)\bigl)\hat{t}'(s)\\
					+\left\langle \D_u\calE\bigl(\hat{t}(s),\hat{u}(s),\hat{z}(s)\bigr),\hat{u}'(s)\right\rangle
					+\left\langle \D_z\calE\bigl(\hat{t}(s),\hat{u}(s),\hat{z}(s)\bigr),\hat{z}'(s)\right\rangle.
				\end{multline*}
				Further, since $\hat{u}'(s)\in\calU^{\tilde{p}}$ is f.a.a. $s\in\left[0,S\right]$ an admissible testfunction in the balance of linear momentum \eqref{balanceoflinearmomentum}, we have
				\begin{align*}
					\text{f.a.a. }s\in\left[0,S\right]:\quad\left\langle \D_u\calE\bigl(\hat{t}(s),\hat{u}(s),\hat{z}(s)\bigr),\hat{u}'(s)\right\rangle = 0.
				\end{align*}
				
				Hence we obtain f.a.a. $s\in\left[0,S\right]$
				\begin{align}\label{CR}
					\frac{\dd{}}{\dd{s}}\Bigl(\calE\bigl(\hat{t}(s),\hat{u}(s),\hat{z}(s)\bigr)\Bigr)
					=\partial_t\calE\bigl(\hat{t}(s),\hat{u}(s),\hat{z}(s)\bigr)\hat{t}(s)
					+\left\langle \D_z\calE\bigl(\hat{t}(s),\hat{u}(s),\hat{z}(s)\bigr),\hat{z}'(s)\right\rangle.
				\end{align}
				Moreover, it holds
				\begin{align*}
					\int_{0}^{s_1}\partial_t\calE\bigl(\hat{t}(s),\hat{u}(s),\hat{z}(s)\bigr)\hat{t}'(s)\dd{s}
					&=\int_{0}^{s_1} -\bigl\langle \dot{\ell}\bigl(\hat{t}(s)\bigr),\hat{u}(s)\bigr\rangle\hat{t}'(s)\dd{s}\\
					&\leq \int_{0}^{s_1}
					\bnorm{\ell}_{C^{1,1}([0,T];\left(\calU^{\tilde{p}}\right)^\ast)}\bnorm{\hat{u}(s)}_{\calU^{\tilde{p}}}\babs{\hat{t}'(s)}\dd{s}\leq C.
				\end{align*}
				Hence \eqref{energyinequality} implies that
				\begin{align*}
					\exists C>0:\quad\int_{0}^{s_1}\bigl\Vert\hat{z}'(s)\bigr\Vert_\calV\dist_{\calV^\ast}\Bigl(-\D_z\calE\bigl(\hat{t}(s),\hat{u}(s),\hat{z}(s)\bigr),\partial^{\calZ}\calR(0)\Bigr)\dd s\leq C,
				\end{align*}
				in particular 
				\begin{align*}
					\text{f.a.a. }s\in\left[0,S\right]:\quad\bigl\Vert\hat{z}'(s)\bigr\Vert_\calV\dist_{\calV^\ast}\Bigl(-\D_z\calE\bigl(\hat{t}(s),\hat{u}(s),\hat{z}(s)\bigr),\partial^{\calZ}\calR(0)\Bigr)<\infty.
				\end{align*}
				Let $s\in\left[0,S\right]$ be arbitrary. If $\norm{\hat{z}'(s)}_\calV=0$, \eqref{CR} implies
				\begin{multline*}
					-\frac{\dd{}}{\dd{s}}\Bigl(\calE\bigl(\hat{t}(s),\hat{u}(s),\hat{z}(s)\bigr)\Bigr)+\partial_t\calE\bigl(\hat{t}(s),\hat{u}(s),\hat{z}(s)\bigr)\hat{t}(s)=0\\
					=\bigl\Vert\hat{z}'(s)\bigr\Vert_\calV\dist_{\calV^\ast}\Bigl(-\D_z\calE\bigl(\hat{t}(s),\hat{u}(s),\hat{z}(s)\bigr),\partial^{\calZ}\calR(0)\Bigr)+\calR\bigl(\hat{z}'(s)\bigr).
				\end{multline*}
			Otherwise, it is $\dist_{\calV^\ast}\Bigl(-\D_z\calE\bigl(\hat{t}(s),\hat{u}(s),\hat{z}(s)\bigr),\partial^{\calZ}\calR(0)\Bigr)<\infty$, and hence by Lemma~\ref{lemma:attaineddist} the distance is attained, i.e. there exists $\sigma(s)\in\partial^\calZ\calR(0)\subset\calZ^\ast$ with
			\begin{align*}
				\dist_{\calV^\ast}\Bigl(-\D_z\calE\bigl(\hat{t}(s),\hat{u}(s),\hat{z}(s)\bigr),\partial^{\calZ}\calR(0)\Bigr)=\bnorm{\D_z\calE\bigl(\hat{t}(s),\hat{u}(s),\hat{z}(s)\bigr)+\sigma(s)}_{\calV^\ast}.
			\end{align*}
			By \eqref{CR} and the characterization of the subdifferential $\partial^\calZ\calR(0)$ from \eqref{prop:sdphdoSubdiffin0} we thus obtain
			\begin{align*}
				-\frac{\dd{}}{\dd{s}}\Bigl(\calE\bigl(\hat{t}(s),\hat{u}(s),\hat{z}(s)\bigr)\Bigr)+\partial_t\calE&\bigl(\hat{t}(s),\hat{u}(s),\hat{z}(s)\bigr)\hat{t}(s)\\
				&=\left\langle -\D_z\calE\bigl(\hat{t}(s),\hat{u}(s),\hat{z}(s)\bigr)-\sigma(s),\hat{z}'(s)\right\rangle+\left\langle\sigma(s),\hat{z}'(s)\right\rangle_{\calZ^\ast,\calZ}\\
				&\leq \bnorm{\D_z\calE\bigl(\hat{t}(s),\hat{u}(s),\hat{z}(s)\bigr)+\sigma(s)}_{\calV^\ast}\bnorm{\hat{z}'(s)}_\calV+\calR(\hat{z}'(s))\\
				&=\bigl\Vert\hat{z}'(s)\bigr\Vert_\calV\dist_{\calV^\ast}\Bigl(-\D_z\calE\bigl(\hat{t}(s),\hat{u}(s),\hat{z}(s)\bigr),\partial^{\calZ}\calR(0)\Bigr)+\calR\bigl(\hat{z}'(s)\bigr).
			\end{align*} 
			So in any case this inequality is valid f.a.a. $s\in\left[0,S\right]$. Finally integration with respect to $s\in\left[0,s_1\right]$ yields equality in \eqref{energyinequality}.
			\end{proof}
		

	\section{Numerical experiments}

Within the previous sections, the convergence properties of an alternate minimization algorithm, combined with an adaptive scheme proposed by Efendiev and Mielke, were elaborated. As a prototype rate-independent constitutive model, a phase-field approximation of fracture mechanics was chosen. This section deals with numerical experiments highlighting the influence of arc-length parameter $\rho$ and the choice of the norm involved in the algorithm by Efendiev and Mielke. These experiments are based on finite element simulations in 2D (plane strain conditions).
 
\subsection{Algorithmic implementation}

In order to calculate a numerical solution of the underlying minimization problems in Algorithm \ref{alg1}, the energy \eqref{energy functional F} is approximated by means of finite elements. For that purpose, the fields $u^\rho_k$ and $z^\rho_k$ are interpolated by using bi-linear shape functions. While the displacement field $u^\rho_k$ is not subject to any constraint (except for those related to the prescribed boundary conditions), the phase-field $z^\rho_k$ has to satisfy the irreversibility condition $z^\rho_k\leq z^\rho_{k-1}$. Furthermore and in line with the adaptive scheme by Efendiev and Mielke, the growth of $z$ is also constrained. To be more precise, $\norm{z-z^\rho_{k-1}}_\calV\leq \rho$. Both constraints are implemented by means of an Augmented-Lagrange method, cf. \cite{geiger2002}. Since the energy to be minimized is separately convex in $u^\rho_k$ and $z^\rho_k$, the respective minima can be conveniently computed by identifying the extrema. For this reason, points characterized by a vanishing derivative are searched for ($\norm{D\mathcal F(t^\rho_k,u^\rho_k,z^\rho_k)}\leq tol$). According to alternate-minimization, this is done for $u^\rho_k$ as well as $z^\rho_k$ in a sequential manner. Within each of the respective steps, a Newton-type iteration is employed. The final algorithm was implemented in $\text{MATLAB}^\copyright$ and Hooke's model was chosen in order to capture the elastic response. It depends on Young's modulus $E$ and Poisson's ratio $\nu$. Clearly, Hooke's model satisfies both the coercivity as well as the symmetry conditions \eqref{assume:elastictensor} and \eqref{assume:symmetrictensor}.

\subsection{Numerical examples}

The modeling capabilities of the algorithm are analyzed and highlighted by means of two different representative examples. The first of those is a Compact-Tension test (Subsection~\ref{sec:ct}), while the second one is an L-shaped plate (Subsection~\ref{sec:L-shape}).

\subsubsection{Compact-Tension test (CT)}
\label{sec:ct}

\begin{figure}[htbp]
\centering
\subcaptionbox{}{
\begin{psfrags}%
\psfrag{s1}[cl][cl][1]{\color[rgb]{0,0,0}\setlength{\tabcolsep}{0pt}\begin{tabular}{c}$\bar u$\end{tabular}}%
\psfrag{l1}[cr][cr][1]{\color[rgb]{0,0,0}\setlength{\tabcolsep}{0pt}\begin{tabular}{c}$1\,\text{mm}$\end{tabular}}%
\psfrag{l2}[cm][cm][1]{\color[rgb]{0,0,0}\setlength{\tabcolsep}{0pt}\begin{tabular}{c}$1\,\text{mm}$\end{tabular}}%
\psfrag{s2}[cl][cl][1]{\color[rgb]{0,0,0}\setlength{\tabcolsep}{0pt}\begin{tabular}{c}$E = 100$ MPa\end{tabular}}%
\psfrag{s3}[cl][cl][1]{\color[rgb]{0,0,0}\setlength{\tabcolsep}{0pt}\begin{tabular}{c}$\nu = 0.3$\end{tabular}}%
\psfrag{s4}[cl][cl][1]{\color[rgb]{0,0,0}\setlength{\tabcolsep}{0pt}\begin{tabular}{c}$g_c = 1\frac{\text{N}}{\text{mm}} $\end{tabular}}%
\includegraphics[scale=1]{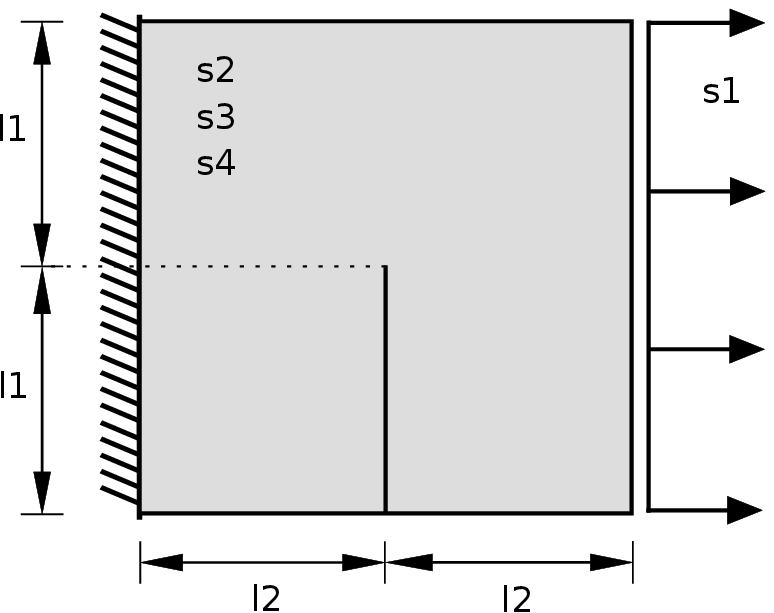}
\end{psfrags}}
\subcaptionbox{}{\includegraphics[scale=0.17]{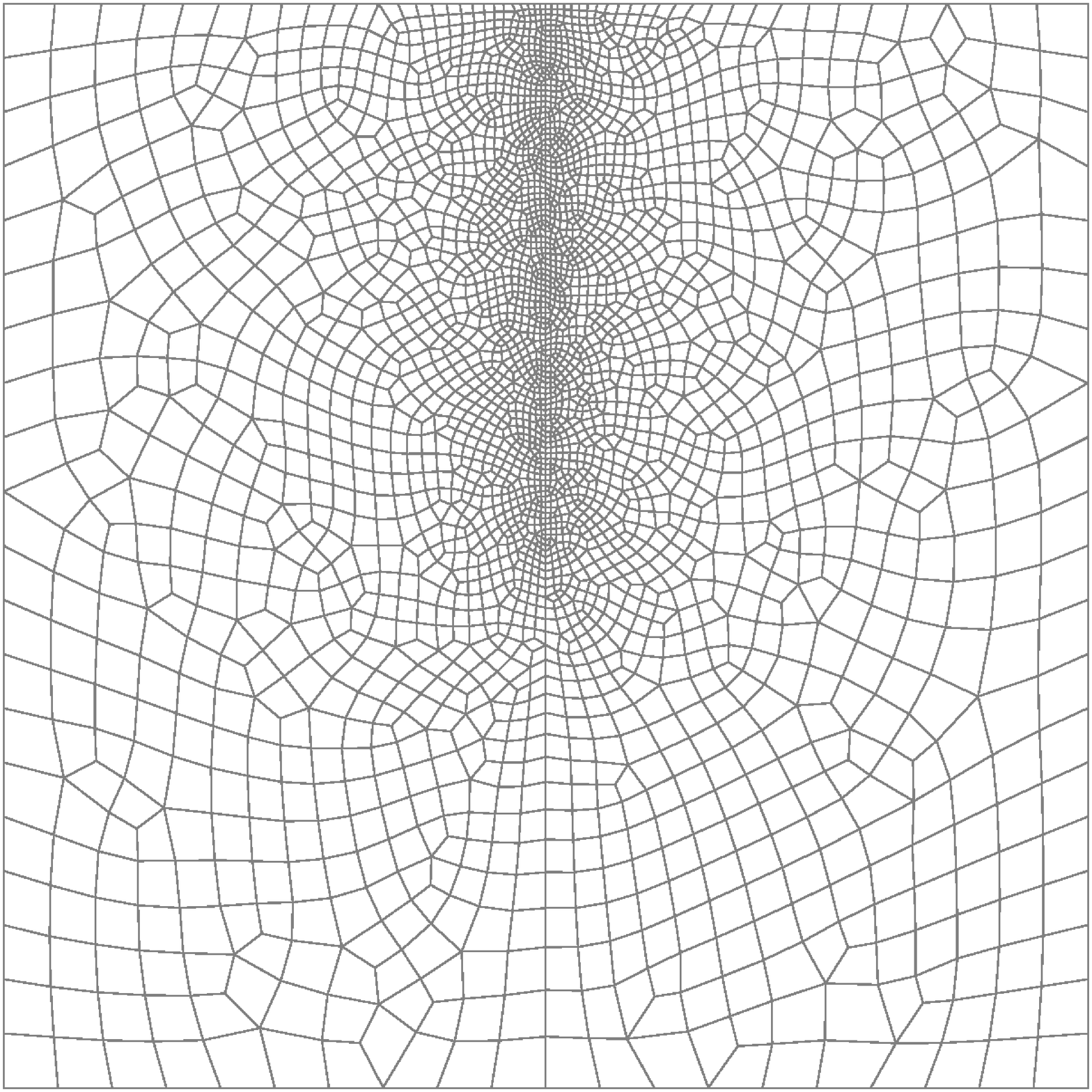}}
\caption{Numerical analysis of a Compact-Tension test (CT). The specimen shows an initial notch: (left) mechanical system; (right) finite element triangulation.}
\label{fig:CT_Specimen_Sketch}
\end{figure}

A widely used experiment in the context of brittle fracture is the Compact-Tension experiment, see Fig.~\ref{fig:CT_Specimen_Sketch}. 
The considered specimen is a squared plate with a pre-existing notch. The specimen is clamped on the left hand side and stretched by prescribing the horizontal displacement $\bar u$ at the right hand side. The finite element mesh used within the numerical simulations is also shown in Fig.~\ref{fig:CT_Specimen_Sketch}. It consists of $3021$ quadrilateral elements whose maximal side-length is $0.1\,\text{mm}$, while the elements' minimal length is $0.01\,\text{mm}$. According to Fig.~\ref{fig:CT_Specimen_Sketch}, the finite element mesh is refined in the region where the crack is expected to propagate. The length-scale parameter controlling the width of the diffuse crack is chosen as $\theta=0.025\,\text{mm}$ and plane strain conditions are assumed.

{\bf Influence of the arc-length increment $\boldsymbol{\rho}$: }
In line with the convergence analysis presented in this paper, the arc-length increment $\rho$ has to be sufficiently small. However and as far as numerical computations are concerned, a small $\rho$ leads to a large number of increments and thus, the resulting numerical costs might be prohibitive. For this reason, the influence of $\rho$ is numerically investigated here --- for a fixed norm $\calV$. Within the numerical study, $\calV=L^4(\Omega)$ is chosen and $\rho$ is varied in the range of $[0.001,\,0.5]$. The loading is prescribed by increasing $\bar u (t^\rho_k)=t^\rho_k \frac{u_{max}}{100\rho}$ with $t_k^\rho\in[0;100\rho]$ and $u_{max}=0.3\,\text{mm}$, leading to a minimum of $100$ discrete time steps for each simulation.

\begin{figure}[htbp]
	\centering
	\subcaptionbox{\label{fig:CT_rho_small}}
	{\begin{psfrags}%
			\psfrag{s1}[tc][tc][1]{\color[rgb]{0,0,0}\setlength{\tabcolsep}{0pt}\begin{tabular}{c}$\bar u$ [mm]\end{tabular}}%
			\psfrag{s2}[bc][bc][1]{\color[rgb]{0,0,0}\setlength{\tabcolsep}{0pt}\begin{tabular}{c}$F$ [N] \end{tabular}}%
			\psfrag{s0001}[cl][cl][0.5]{\color[rgb]{0,0,0}\setlength{\tabcolsep}{0pt}\begin{tabular}{c}$\rho=0.001$ \end{tabular}}%
			\psfrag{s0002}[cl][cl][0.5]{\color[rgb]{0,0,0}\setlength{\tabcolsep}{0pt}\begin{tabular}{c}$\rho=0.002$ \end{tabular}}%
			\psfrag{s0003}[cl][cl][0.5]{\color[rgb]{0,0,0}\setlength{\tabcolsep}{0pt}\begin{tabular}{c}$\rho=0.005$ \end{tabular}}%
			\psfrag{s0004}[cl][cl][0.5]{\color[rgb]{0,0,0}\setlength{\tabcolsep}{0pt}\begin{tabular}{c}$\rho=0.01$ \end{tabular}}%
			\psfrag{s0005}[cl][cl][0.5]{\color[rgb]{0,0,0}\setlength{\tabcolsep}{0pt}\begin{tabular}{c}\textcolor{red}{Quelle!} \end{tabular}}%
			\includegraphics[scale=0.4]{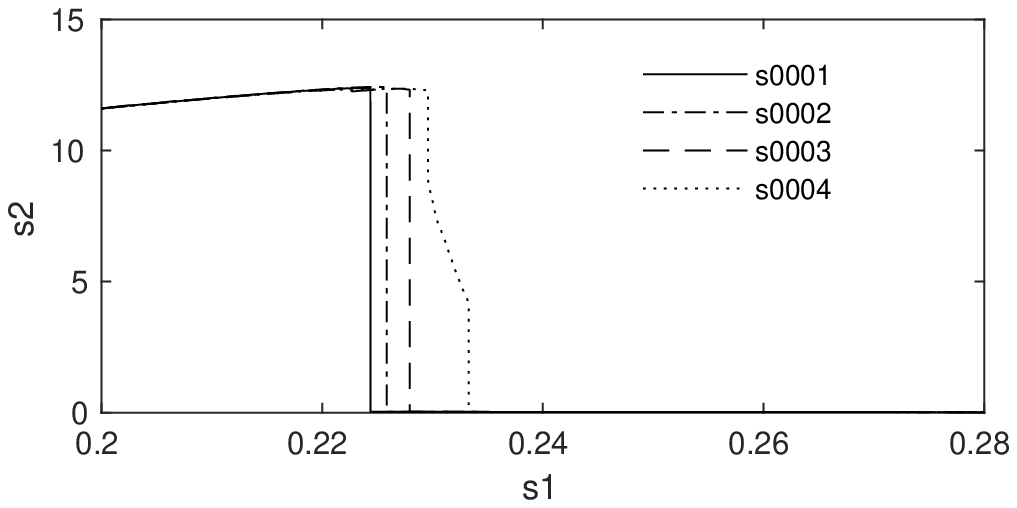}
	\end{psfrags}}
	\subcaptionbox{\label{fig:CT_FU_rho_medium}}
	{\begin{psfrags}%
			\psfrag{s1}[tc][tc][1]{\color[rgb]{0,0,0}\setlength{\tabcolsep}{0pt}\begin{tabular}{c}$\bar u$ [mm]\end{tabular}}%
			\psfrag{s2}[bc][bc][1]{\color[rgb]{0,0,0}\setlength{\tabcolsep}{0pt}\begin{tabular}{c}$F$ [N] \end{tabular}}%
			\psfrag{s0001}[cl][cl][0.5]{\color[rgb]{0,0,0}\setlength{\tabcolsep}{0pt}\begin{tabular}{c}$\rho=0.02$ \end{tabular}}%
			\psfrag{s0002}[cl][cl][0.5]{\color[rgb]{0,0,0}\setlength{\tabcolsep}{0pt}\begin{tabular}{c}$\rho=0.03$ \end{tabular}}%
			\psfrag{s0003}[cl][cl][0.5]{\color[rgb]{0,0,0}\setlength{\tabcolsep}{0pt}\begin{tabular}{c}$\rho=0.04$ \end{tabular}}%
			\psfrag{s0004}[cl][cl][0.5]{\color[rgb]{0,0,0}\setlength{\tabcolsep}{0pt}\begin{tabular}{c}$\rho=0.05$ \end{tabular}}%
			\includegraphics[scale=0.4]{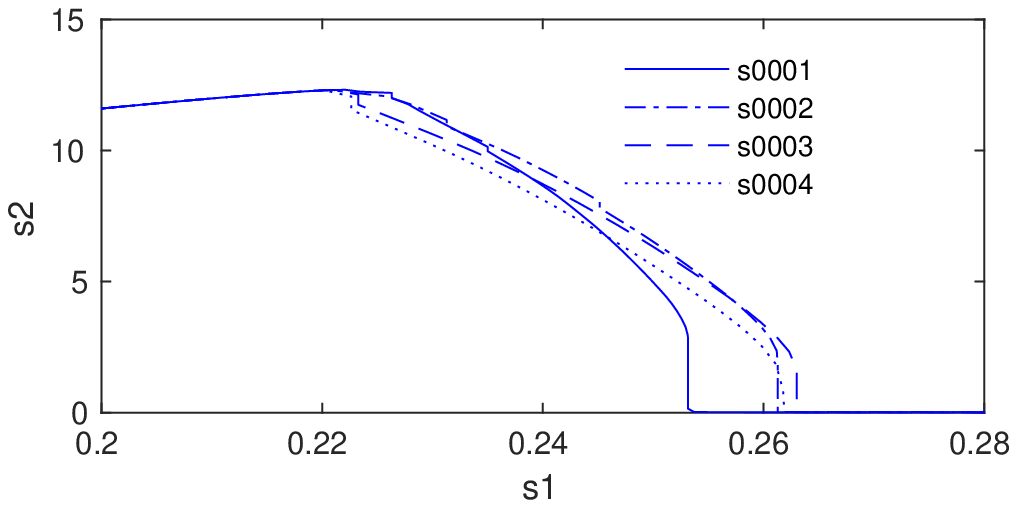}
	\end{psfrags}}
	\subcaptionbox{\label{fig:CT_FU_rho_large}}
	{\begin{psfrags}%
			\psfrag{s1}[tc][tc][1]{\color[rgb]{0,0,0}\setlength{\tabcolsep}{0pt}\begin{tabular}{c}$\bar u$ [mm]\end{tabular}}%
			\psfrag{s2}[bc][bc][1]{\color[rgb]{0,0,0}\setlength{\tabcolsep}{0pt}\begin{tabular}{c}$F$ [N] \end{tabular}}%
			\psfrag{s0001}[cl][cl][0.5]{\color[rgb]{0,0,0}\setlength{\tabcolsep}{0pt}\begin{tabular}{c}$\rho=0.1$ \end{tabular}}%
			\psfrag{s0002}[cl][cl][0.5]{\color[rgb]{0,0,0}\setlength{\tabcolsep}{0pt}\begin{tabular}{c}$\rho=0.2$ \end{tabular}}%
			\psfrag{s0003}[cl][cl][0.5]{\color[rgb]{0,0,0}\setlength{\tabcolsep}{0pt}\begin{tabular}{c}$\rho=0.5$ \end{tabular}}%
			\psfrag{s0004}[cl][cl][0.5]{\color[rgb]{0,0,0}\setlength{\tabcolsep}{0pt}\begin{tabular}{c}AM\end{tabular}}%
			\includegraphics[scale=0.4]{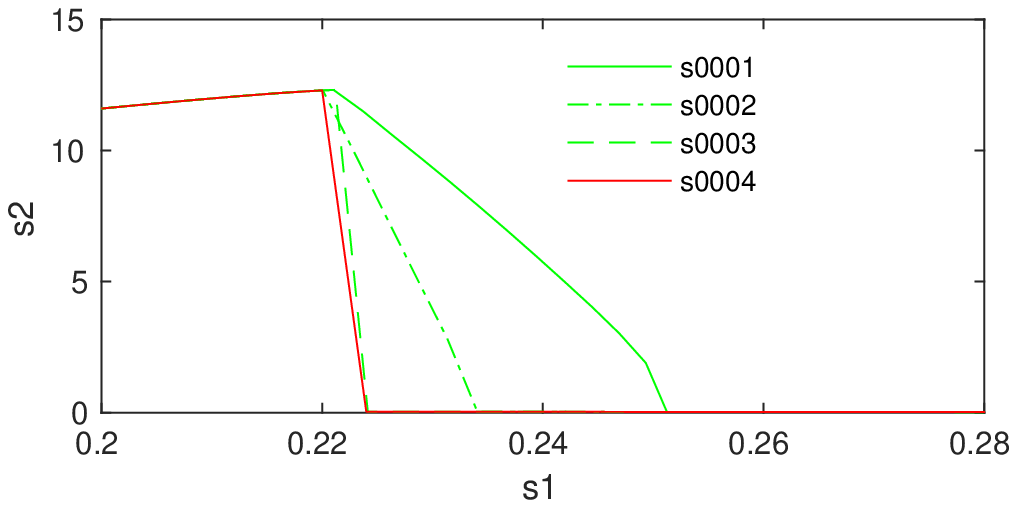}
	\end{psfrags}}
	\subcaptionbox{\label{fig:CT_FU_rho_diff}}
	{\begin{psfrags}%
			\psfrag{s1}[tc][tc][1]{\color[rgb]{0,0,0}\setlength{\tabcolsep}{0pt}\begin{tabular}{c}$\bar u$ [mm]\end{tabular}}%
			\psfrag{s2}[bc][bc][1]{\color[rgb]{0,0,0}\setlength{\tabcolsep}{0pt}\begin{tabular}{c}$F$ [N] \end{tabular}}%
			\psfrag{s0001}[cl][cl][0.5]{\color[rgb]{0,0,0}\setlength{\tabcolsep}{0pt}\begin{tabular}{c}$\rho=0.001$ \end{tabular}}%
			\psfrag{s0002}[cl][cl][0.5]{\color[rgb]{0,0,0}\setlength{\tabcolsep}{0pt}\begin{tabular}{c}$\rho=0.01$ \end{tabular}}%
			\psfrag{s0003}[cl][cl][0.5]{\color[rgb]{0,0,0}\setlength{\tabcolsep}{0pt}\begin{tabular}{c}$\rho=0.04$ \end{tabular}}%
			\psfrag{s0004}[cl][cl][0.5]{\color[rgb]{0,0,0}\setlength{\tabcolsep}{0pt}\begin{tabular}{c}$\rho=0.5$ \end{tabular}}%
			\psfrag{s0005}[cl][cl][0.5]{\color[rgb]{0,0,0}\setlength{\tabcolsep}{0pt}\begin{tabular}{c}AM \end{tabular}}%
			\hspace{-0.5em}\includegraphics[scale=0.5]{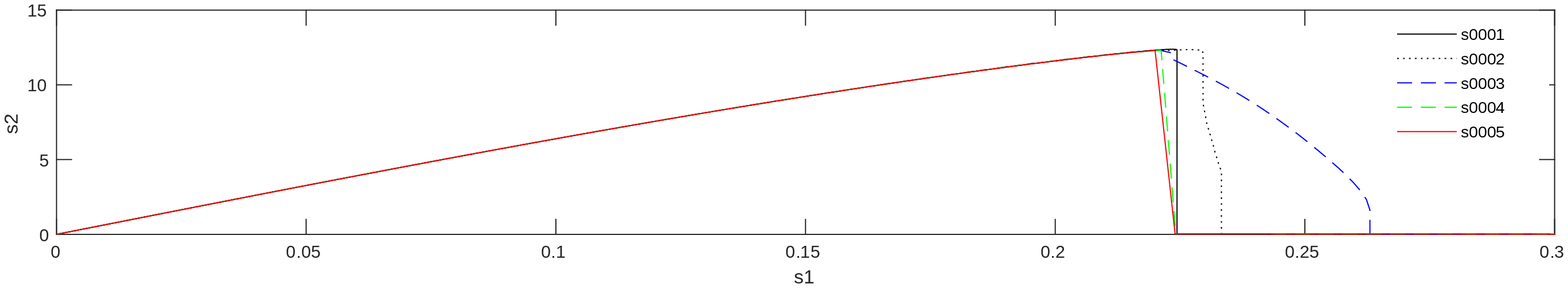}
	\end{psfrags}}
	\caption{Numerical analysis of a Compact-Tension test (CT): Force-displacement curves for different arc-length increments $\rho$. For the sake of comparison, results obtained from the standard non-adaptive alternate-minimization algorithm are also shown. Upper row: zoom in of the force-displacement curves within the time interval in which brutal crack growths occurs. Lower diagram: complete force-displacement curves for some $\rho$.}\label{fig:CT_rhoexperiment}
\end{figure}

The force-displacement curves (reaction force $F$ where $\bar u$ is applied) computed by a finite element implementation based on alternate-minimization combined with the adaptive scheme by Efendiev and Mielke are summarized in Figure \ref{fig:CT_rhoexperiment} for different arc-length increments $\rho$.
For the sake of comparison, the results obtained from the alternate-minimization algorithm without the adaptive scheme (cf. \cite{KN2017}, \cite{Miehe2010}) are also included in the diagrams (highlighted in red color). According to the right diagram in the upper row of Fig.~\ref{fig:CT_rhoexperiment} and as expected, the adaptive scheme is inactive if $\rho$ is sufficiently large. Consequently, the results converge to those of the classic alternate-minimization algorithm. Likewise, the numerical scheme seems to converge for $\rho\to 0$ --- in line with the convergence analysis. This can be seen in the left diagram in the upper row of Fig.~\ref{fig:CT_rhoexperiment}. In the case of the analyzed example, the converged state seems to be close to that corresponding to standard alternate minimization. However, this observation is not universally valid. The structural response associated with the arc-length increments between the aforementioned two limiting cases are given in the middle diagram in the upper row of Fig.~\ref{fig:CT_rhoexperiment}. Here, no physically sound brutal crack growth is predicted any more. It is interesting to note that the dependence of the results on the arc-length increment $\rho$ are not monotonic. To be more explicit, sufficiently small as well as a sufficiently large $\rho$ lead to results similar to those of standard non-adaptive alternate-minimization. By way of contrast, a completely different structural response is predicted for values of $\rho$ in between.

\begin{figure}[htbp]
	\centering
	\subcaptionbox{$\rho=0.05$\label{fig:deltaz_rho8}}
	{\begin{psfrags}%
			\psfrag{s1}[tc][tc][1]{\color[rgb]{0,0,0}\setlength{\tabcolsep}{0pt}\begin{tabular}{c}$k$\end{tabular}}%
			\psfrag{s2}[bc][bc][1]{\color[rgb]{0,0,0}\setlength{\tabcolsep}{0pt}\begin{tabular}{c}$\Delta t^\rho_{k+1}$ \end{tabular}}%
			\includegraphics[scale=0.4]{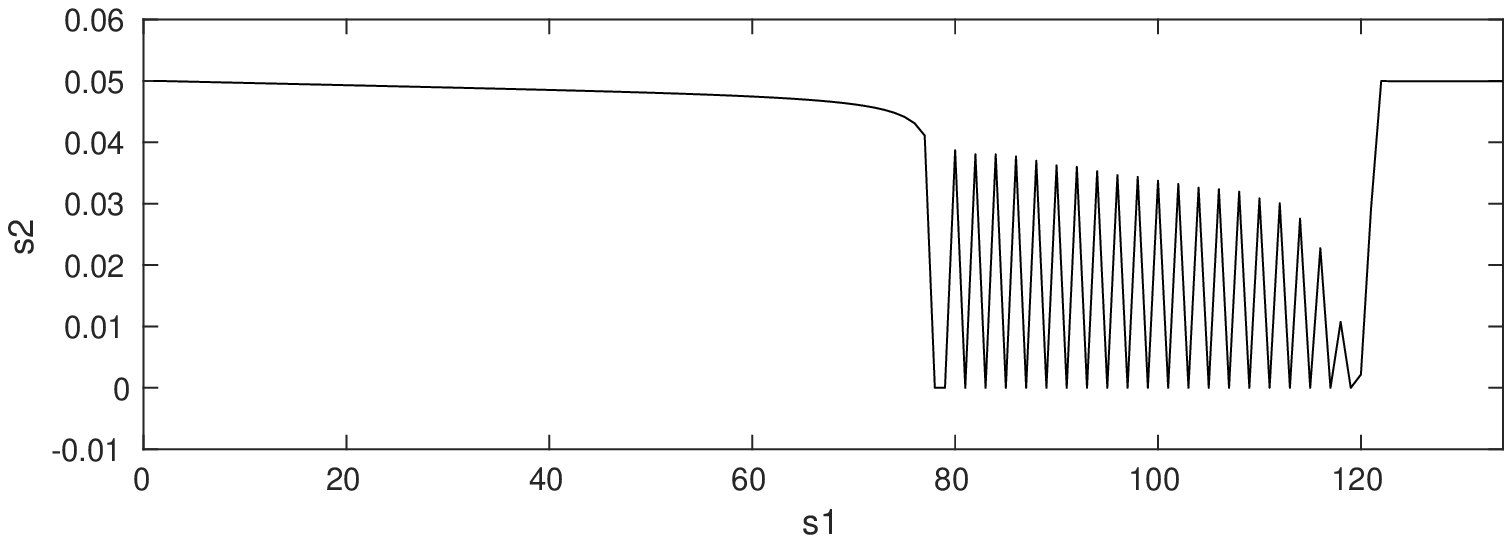}
	\end{psfrags}}
	\subcaptionbox{$\rho=0.005$\label{fig:CT_deltaz_rho3}}
	{\begin{psfrags}%
			\psfrag{s1}[tc][tc][1]{\color[rgb]{0,0,0}\setlength{\tabcolsep}{0pt}\begin{tabular}{c}$k$\end{tabular}}%
			\psfrag{s2}[bc][bc][1]{\color[rgb]{0,0,0}\setlength{\tabcolsep}{0pt}\begin{tabular}{c}$\Delta t^\rho_{k+1}$ \end{tabular}}%
			\includegraphics[scale=0.4]{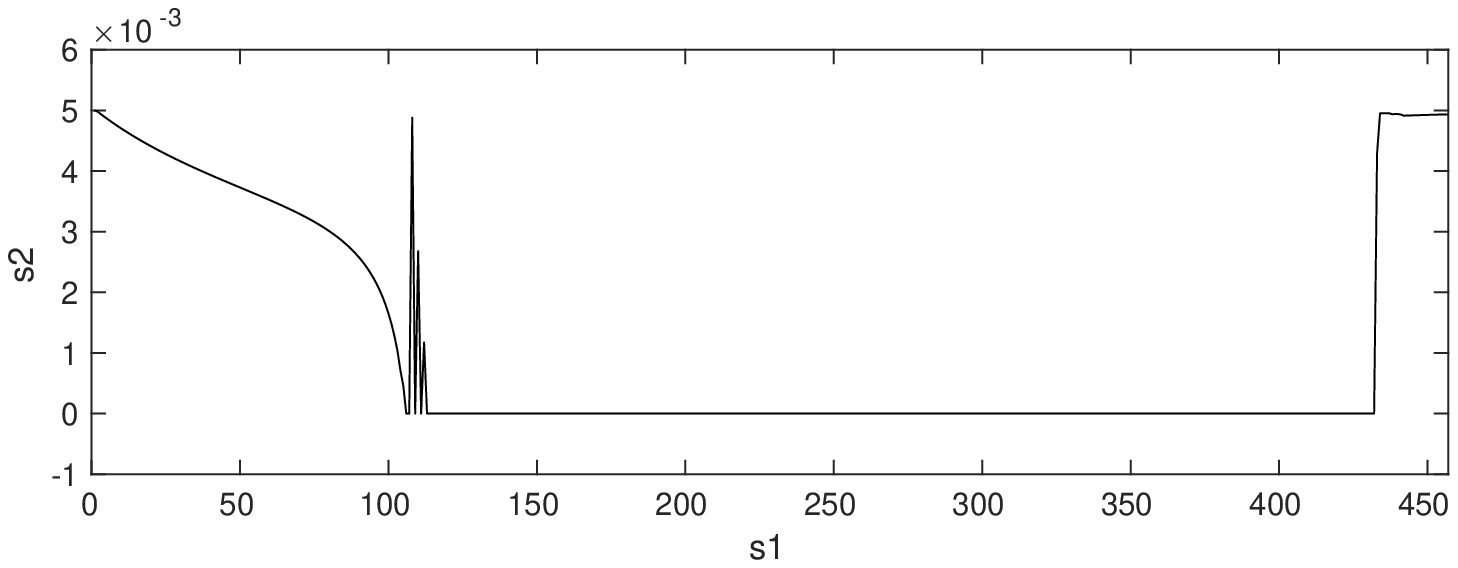}
	\end{psfrags}}
	\caption{Numerical analysis of a Compact-Tension test (CT): Time increments $\Delta t^\rho_{k+1}$ for each iteration computed from the adaptive scheme according to Efendiev Mielke.}\label{fig:CT_rhoexperiment_2}
\end{figure}

\begin{figure}[htbp]
	\centering
	\subcaptionbox{$\bar u(t^\rho_{1})=0$\label{fig:CT_Plots1}}
	{\includegraphics[scale=0.13]{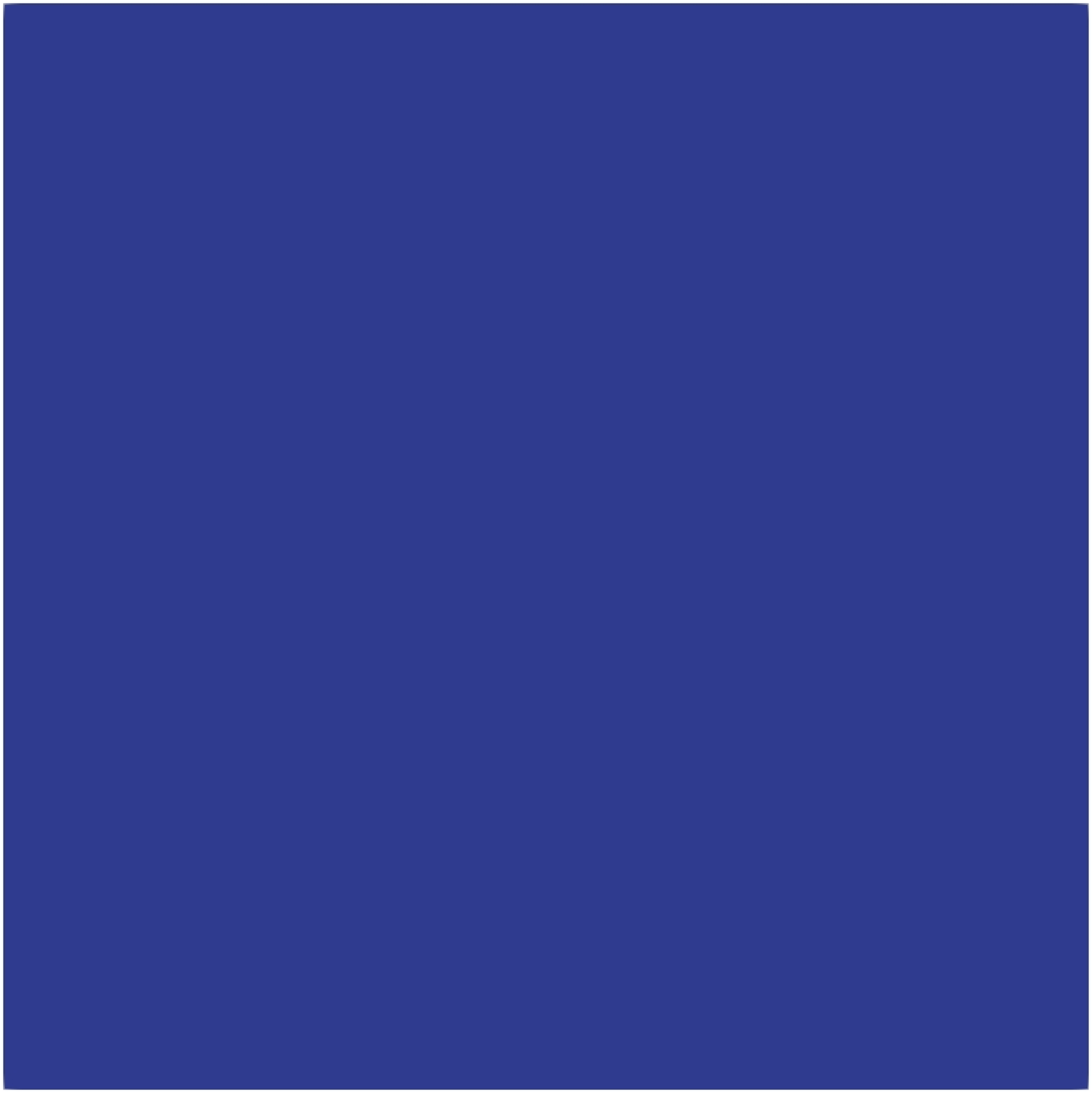}}
	\subcaptionbox{$\bar u(t^\rho_{121})=0.228\,\text{[mm]}$\label{fig:CT_Plots2}}
	{\includegraphics[scale=0.13]{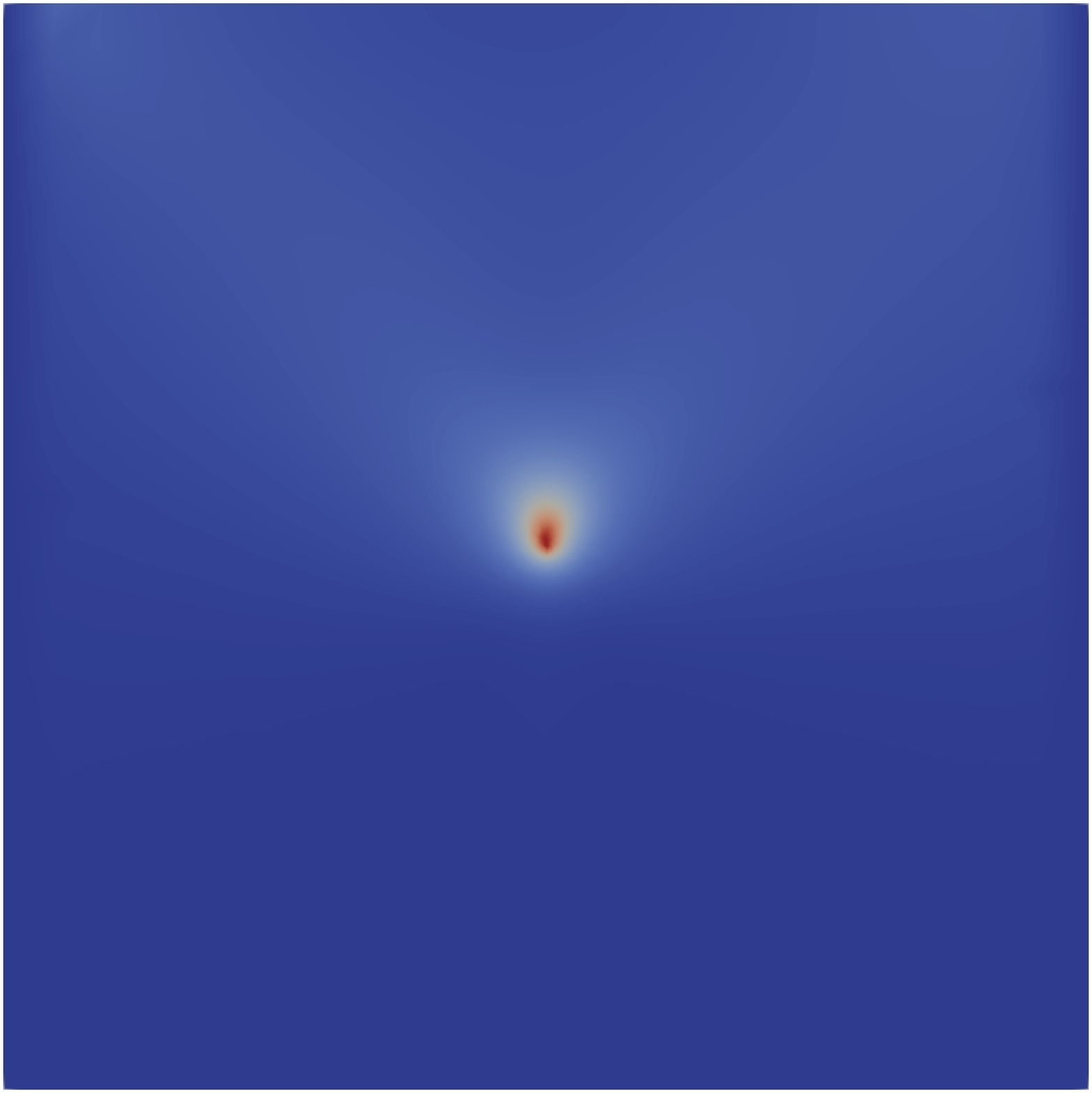}}
	\subcaptionbox{$\bar u(t^\rho_{180})=0.228\,\text{[mm]}$\label{fig:CT_Plots3}}
	{\includegraphics[scale=0.13]{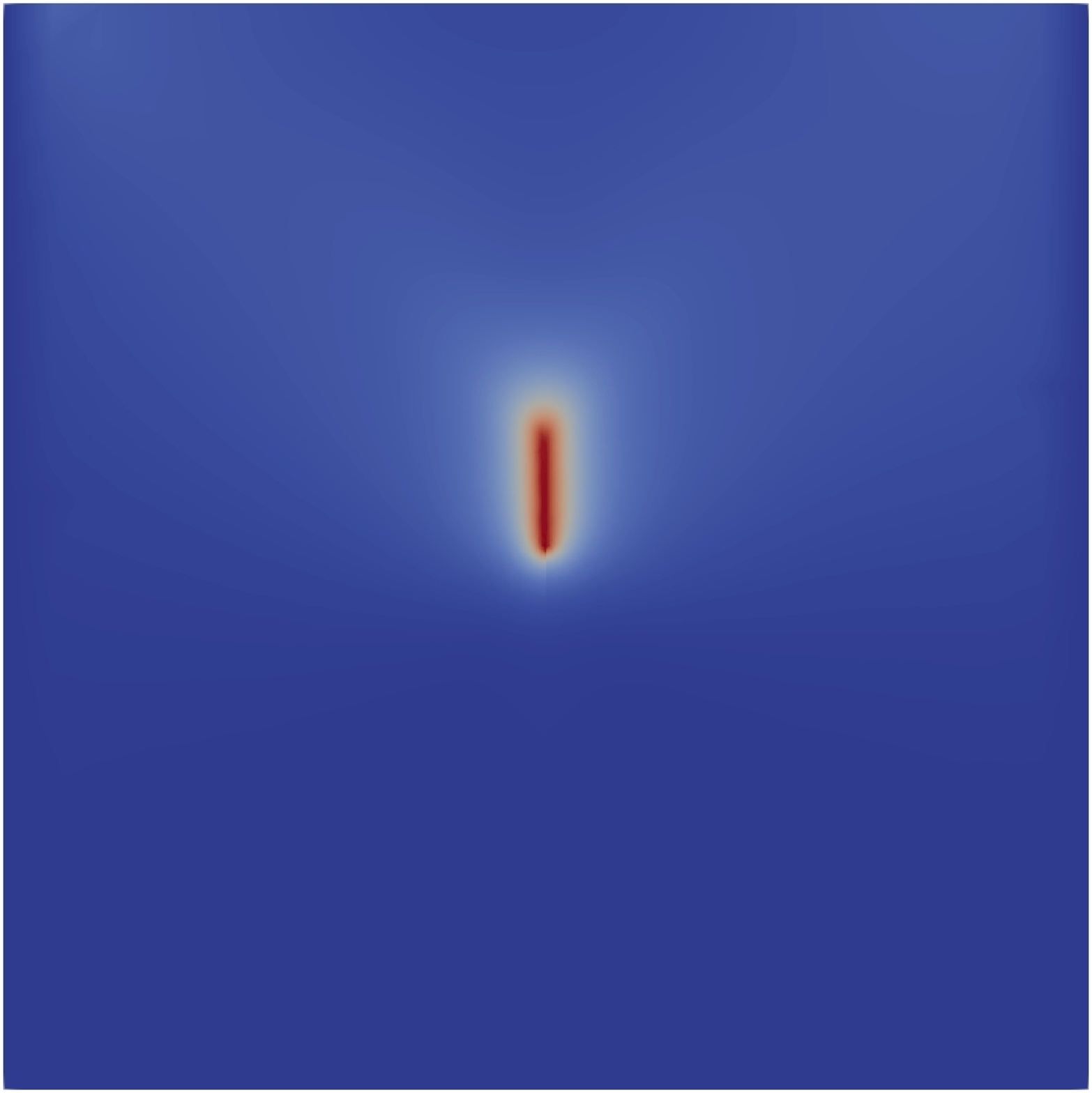}}
	\subcaptionbox{$\bar u(t^\rho_{230})=0.228\,\text{[mm]}$\label{fig:CT_Plots4}}
	{\includegraphics[scale=0.13]{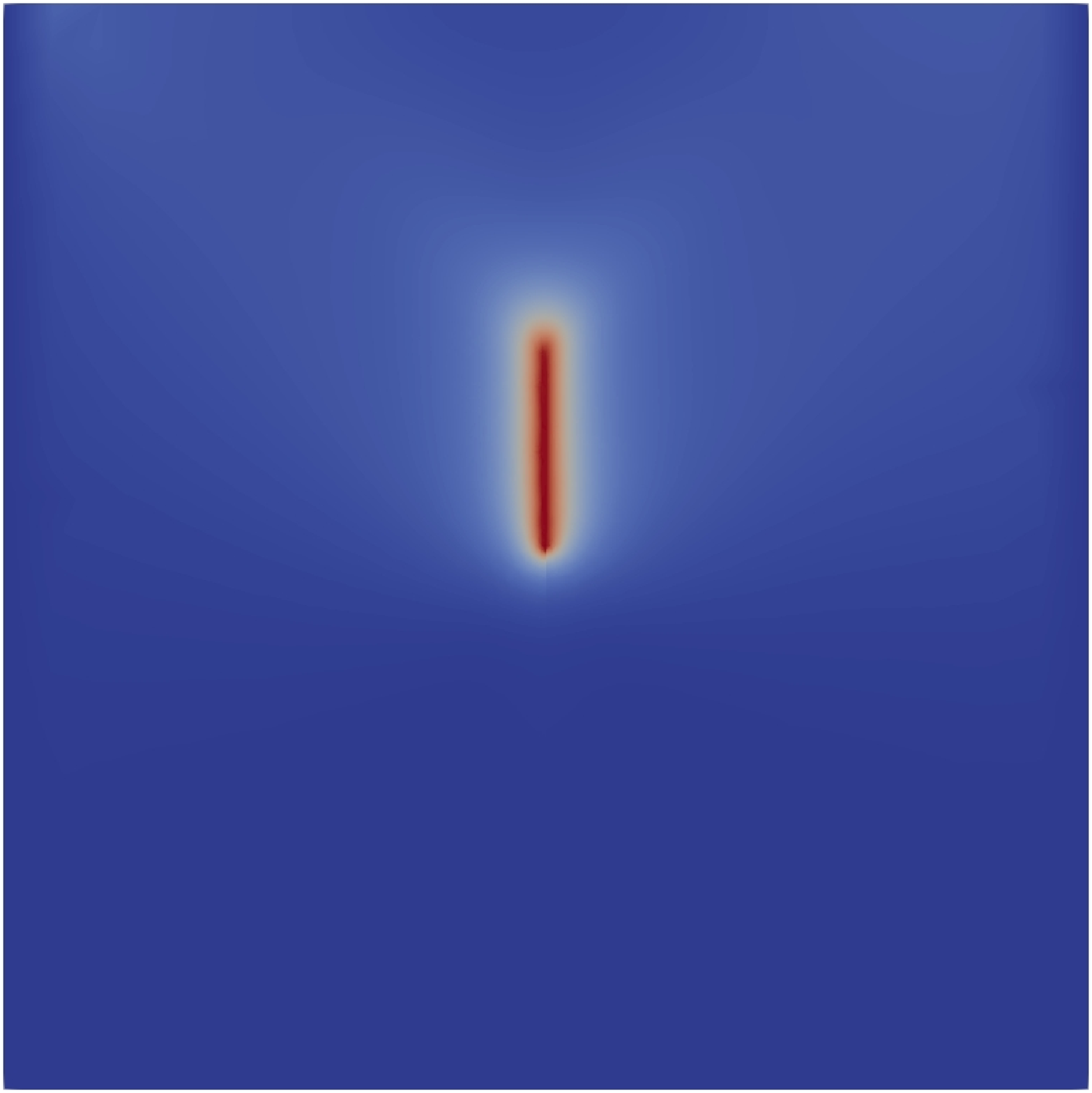}}
	\subcaptionbox{$\bar u(t^\rho_{280})=0.228\,\text{[mm]}$\label{fig:CT_Plots5}}
	{\includegraphics[scale=0.13]{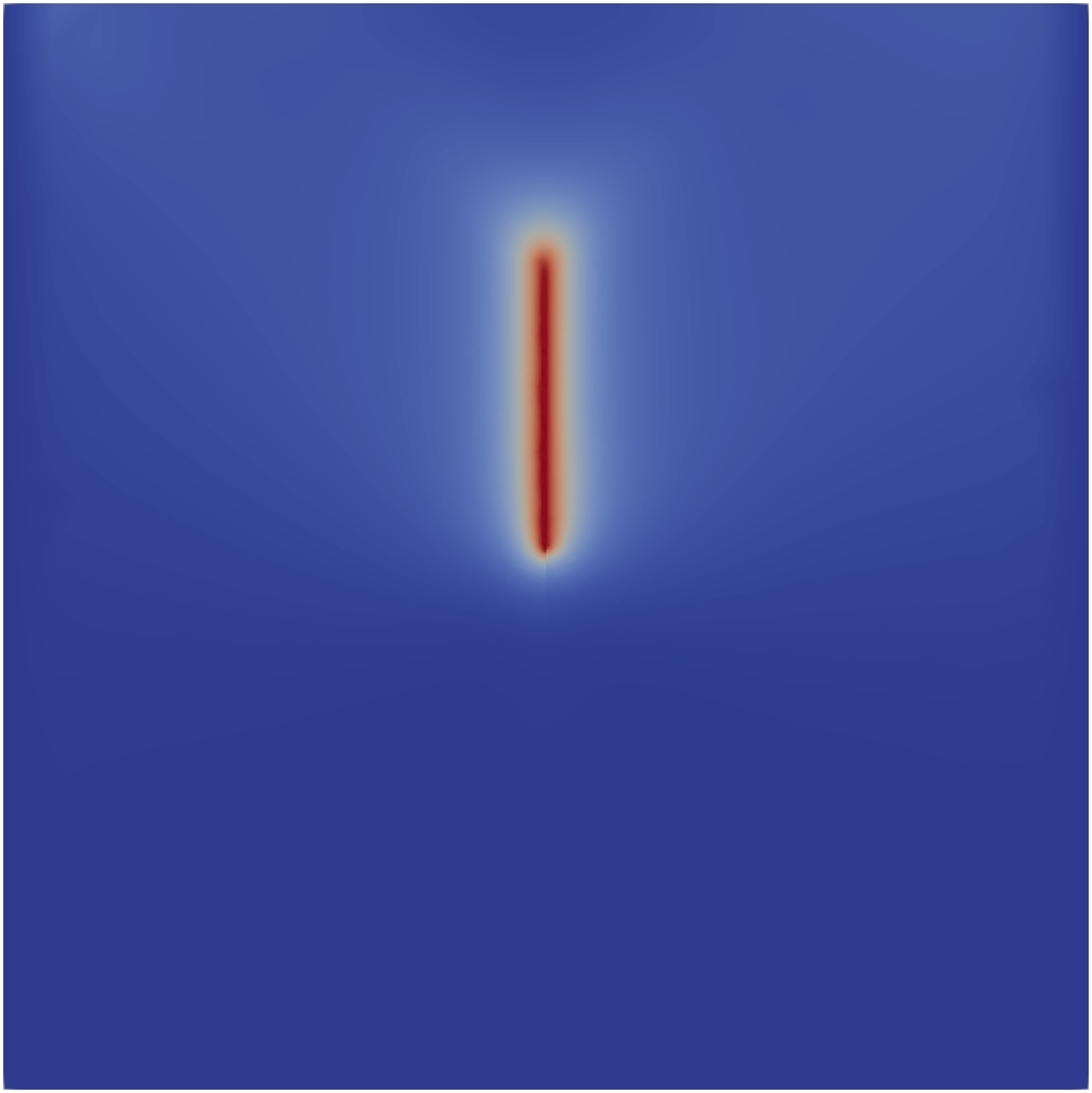}}
	\subcaptionbox{$\bar u(t^\rho_{330})=0.228\,\text{[mm]}$\label{fig:CT_Plots6}}
	{\includegraphics[scale=0.13]{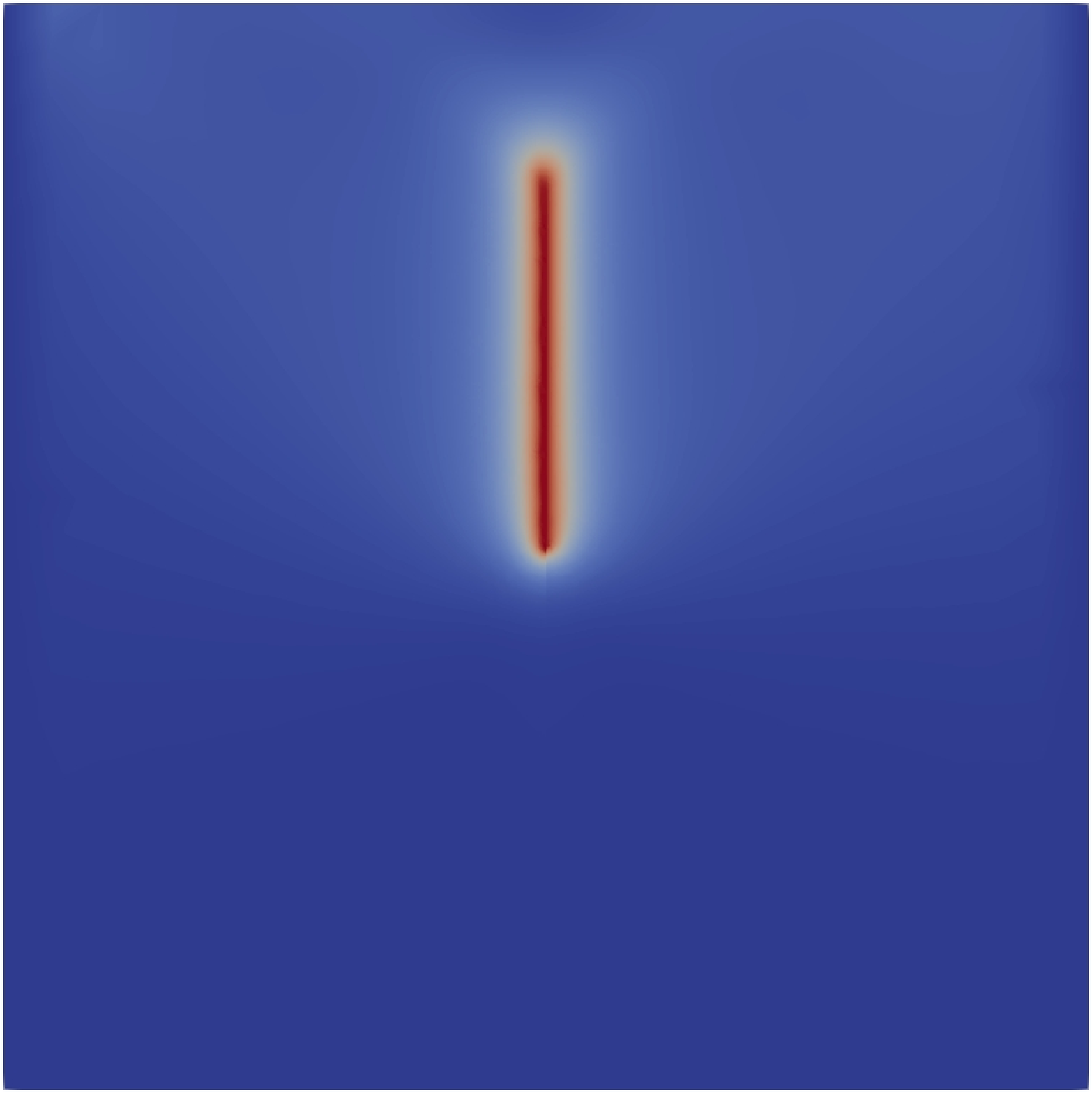}}
	\subcaptionbox{$\bar u(t^\rho_{380})=0.228\,\text{[mm]}$\label{fig:CT_Plots7}}
	{\includegraphics[scale=0.13]{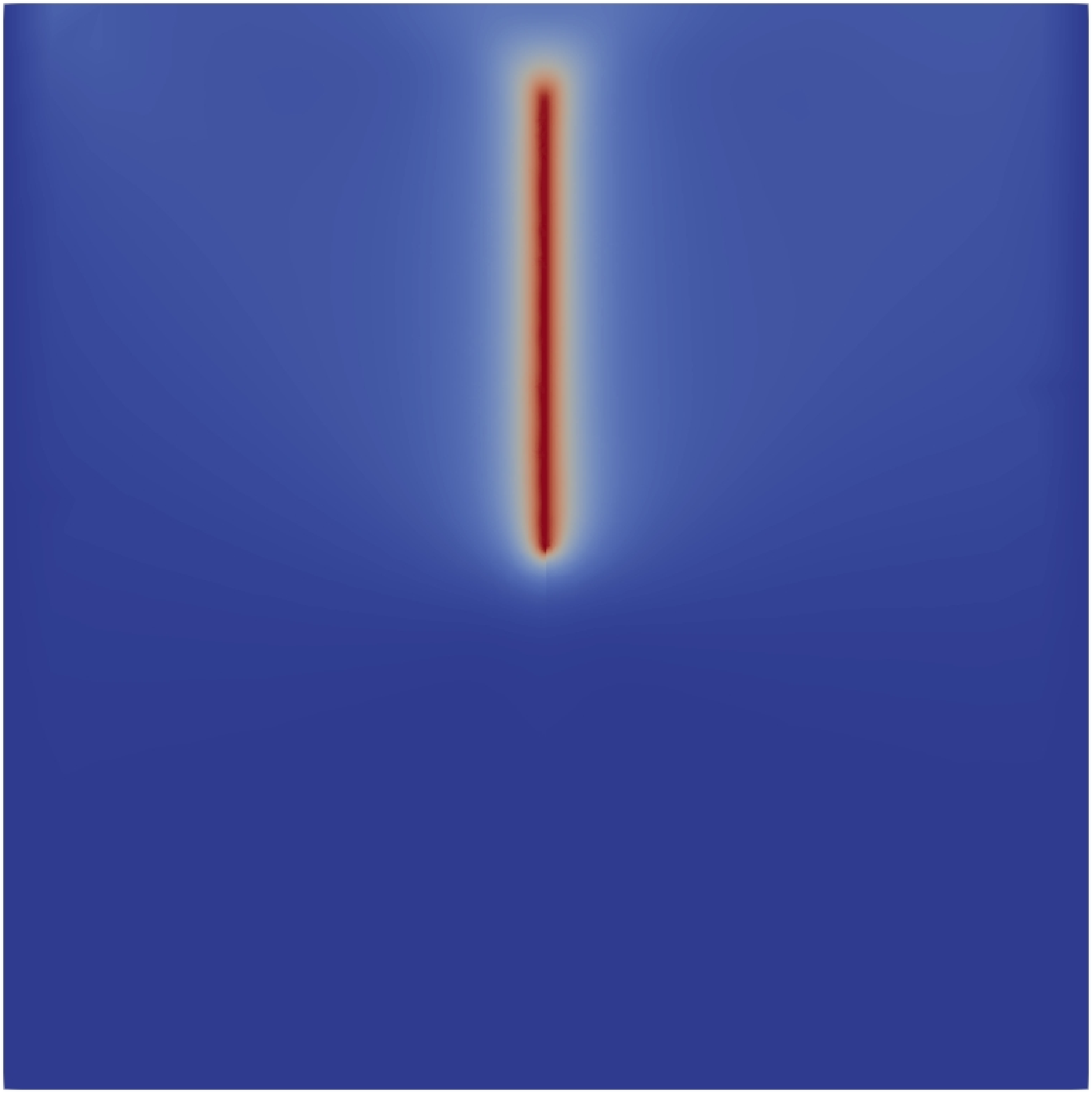}}
	\subcaptionbox{$\bar u(t^\rho_{433})=0.228\,\text{[mm]}$\label{fig:CT_Plots8}}
	{\includegraphics[scale=0.13]{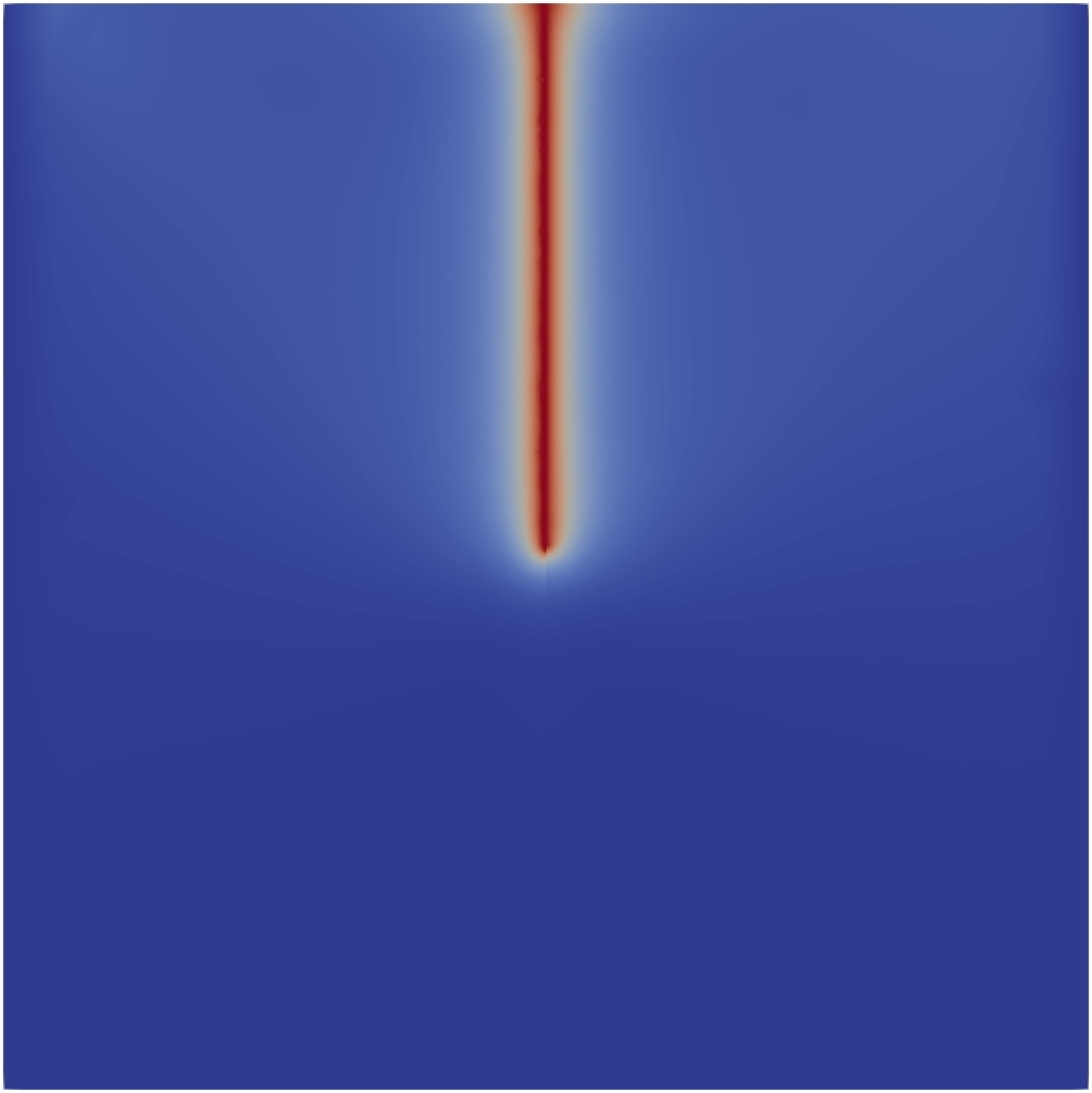}}
	\caption{Numerical analysis of a Compact-Tension test (CT): Distribution of the phase-field variable $z$ computed from the Efendiev \& Mielke adaptive scheme combined with alternate minimization for $\rho=0.005$ and $\alpha=4$. Brutal crack-growth occurs at $\bar u=0.228$ [mm].}\label{fig:CT_Plots}
\end{figure}

For further insights on the influence of $\rho$ on the numerical results, the physical time step size as predicted by the adaptive scheme by Efendiev and Mielke is provided in Fig.~\ref{fig:CT_rhoexperiment_2} for two different arc-length increments $\rho$.
While $\rho=0.05$ (left diagram in Fig.~\ref{fig:CT_rhoexperiment_2}) corresponds to the increment range in which physically sound brutal crack growth is not observed (middle diagram in the upper row of Fig.~\ref{fig:CT_rhoexperiment}), $\rho=0.005$ (right diagram in Fig.~\ref{fig:CT_rhoexperiment_2}) is associated with computations presumably relatively close to the balanced-viscosity solution characterized by the Efendiev \& Mielke algorithm. Accordingly, brutal crack growth can be identified by zero time increments (right diagram in Fig.~\ref{fig:CT_rhoexperiment_2}). By way of contrast, the unphysical range of the increment $\rho$ leads to oscillations and brutal crack growth is not observed (left diagram in Fig.~\ref{fig:CT_rhoexperiment_2}). The temporal evolution of the phase-field for the physically sound computation based on $\rho=0.005$ is provided in Fig.~\ref{fig:CT_Plots}.
As it can be seen from Fig.~\ref{fig:CT_Plots}, the adaptive scheme proposed by Efendiev \& Mielke predicts brutal crack growth at a prescribed displacement $\bar u^\rho_k=0.228$~mm.

\begin{figure}[htbp]
	\centering
	\subcaptionbox{$\rho=0.02$}
	{\begin{psfrags}%
			\psfrag{steps}[tc][tc][1]{\color[rgb]{0,0,0}\setlength{\tabcolsep}{0pt}\begin{tabular}{c}$k$\end{tabular}}%
			\psfrag{deltaT}[bc][bc][1]{\color[rgb]{0,0,0}\setlength{\tabcolsep}{0pt}\begin{tabular}{c}$\Delta t^\rho_{k+1}$ \end{tabular}}%
			\includegraphics[scale=0.5]{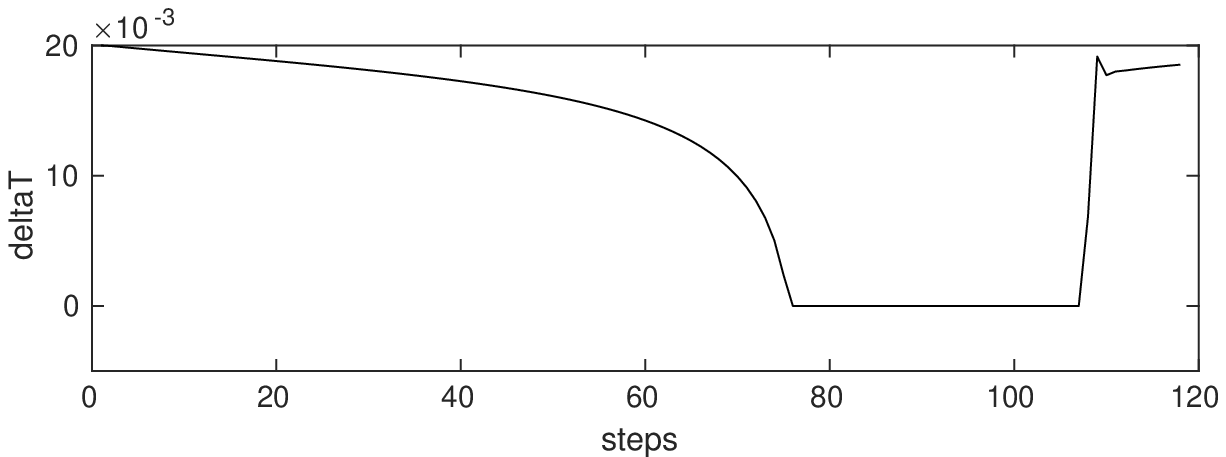}
	\end{psfrags}}
	\subcaptionbox{$\rho=0.002$}
	{\begin{psfrags}%
			\psfrag{steps}[tc][tc][1]{\color[rgb]{0,0,0}\setlength{\tabcolsep}{0pt}\begin{tabular}{c}$k$\end{tabular}}%
			\psfrag{deltaT}[bc][bc][1]{\color[rgb]{0,0,0}\setlength{\tabcolsep}{0pt}\begin{tabular}{c}$\Delta t^\rho_{k+1}$ \end{tabular}}%
			\includegraphics[scale=0.5]{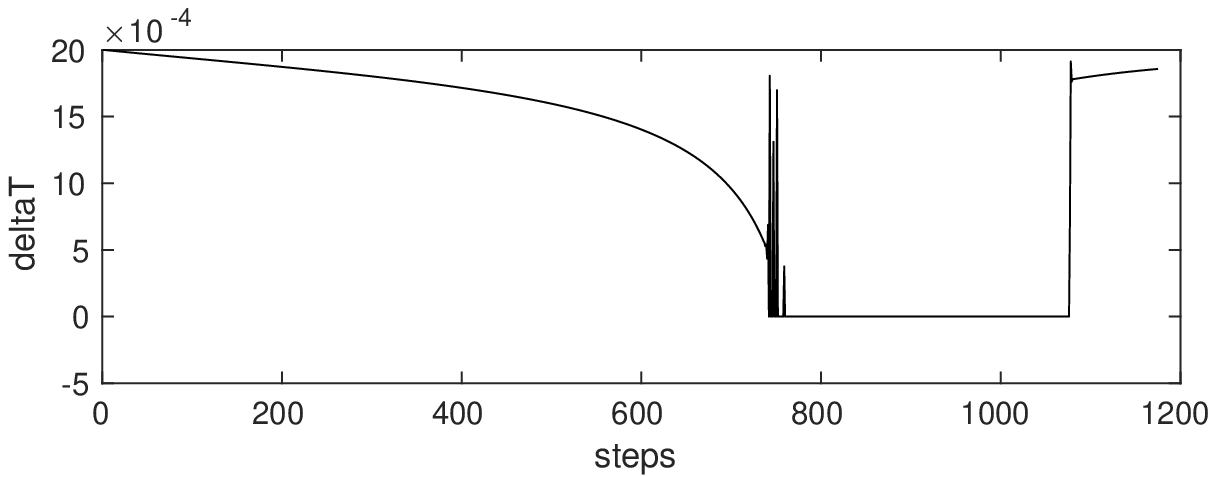}
	\end{psfrags}}
	\caption{Numerical analysis of the $0$-dimensional example investigated in \cite{KN2017}: Time increments $\Delta t^\rho_{k+1}$ for each E\&M-iterate $k$ for two different arc-length increments $\rho$.}\label{fig:0D_time_increments}
\end{figure}

By further analyzing the physically sound solution corresponding to the right diagram in Fig.~\ref{fig:CT_rhoexperiment_2} (see also Fig.~\ref{fig:CT_Plots}), a short artificial oscillation is noted at the beginning of the brutal crack growth (see right diagram in Fig.~\ref{fig:CT_rhoexperiment_2}). In order to trace back the reason for this undesired behavior, the zero-dimensional example investigated in \cite{KN2017} is re-analyzed. By doing so, the influence of the finite element discretization can be ruled out and different solvers, e.g. alternate-minimization and others, can be used for solving the underlying minimization problems. The predicted time increments are summarized in Fig.~\ref{fig:0D_time_increments}.
Again, oscillations can be seen --- although less pronounced for $\rho=0.02$. It can be thus concluded that these oscillations are neither an effect of the finite-element-discretization nor an effect of the minimization algorithm. One possible reason for these oscillations is the explicit nature of the time update within the adaptive Efendiev \& Mielke algorithm. However, further analysis of this effect are beyond the scope of the present paper, but will be analyzed in a forthcoming publication.
 
Next, the influence of the norm occurring in the Efendiev \& Mielke adaptive scheme is analyzed. For that purpose, the family of $L^\alpha$ norms is considered. The loading $\bar u$ is prescribed by $\bar u (t^\rho_k)=t^\rho_k u_{max}$ with $u_{max}=0.3\,\text{mm}$. Furthermore, the arc-length increments are set to $\rho=0.01$ and the maximal time is $T=1$. As a consequence, $\bar u$ is increased up to $u_{max}$ in at least 100 steps.

\begin{figure}[htbp]
	\centering
	\subcaptionbox{\label{fig:CT_FUcurves_lp}}
	{\begin{psfrags}%
			\psfrag{s1}[tc][tc][1]{\color[rgb]{0,0,0}\setlength{\tabcolsep}{0pt}\begin{tabular}{c}$\bar u$ [mm]\end{tabular}}%
			\psfrag{s2}[bc][bc][1]{\color[rgb]{0,0,0}\setlength{\tabcolsep}{0pt}\begin{tabular}{c}$F$ [N] \end{tabular}}%
			\psfrag{s0001}[cl][cl][0.5]{\color[rgb]{0,0,0}\setlength{\tabcolsep}{0pt}\begin{tabular}{c}$\alpha=4$\end{tabular}}%
			\psfrag{s0002}[cl][cl][0.5]{\color[rgb]{0,0,0}\setlength{\tabcolsep}{0pt}\begin{tabular}{c}$\alpha=6$\end{tabular}}%
			\psfrag{s0003}[cl][cl][0.5]{\color[rgb]{0,0,0}\setlength{\tabcolsep}{0pt}\begin{tabular}{c}$\alpha=8$\end{tabular}}%
			\psfrag{s0004}[cl][cl][0.5]{\color[rgb]{0,0,0}\setlength{\tabcolsep}{0pt}\begin{tabular}{c}$\alpha=10$\end{tabular}}%
			\psfrag{s0005}[cl][cl][0.5]{\color[rgb]{0,0,0}\setlength{\tabcolsep}{0pt}\begin{tabular}{c}$\alpha=12$\end{tabular}}%
			\psfrag{s0006}[cl][cl][0.5]{\color[rgb]{0,0,0}\setlength{\tabcolsep}{0pt}\begin{tabular}{c}$\alpha=14$\end{tabular}}%
			\includegraphics[scale=0.4]{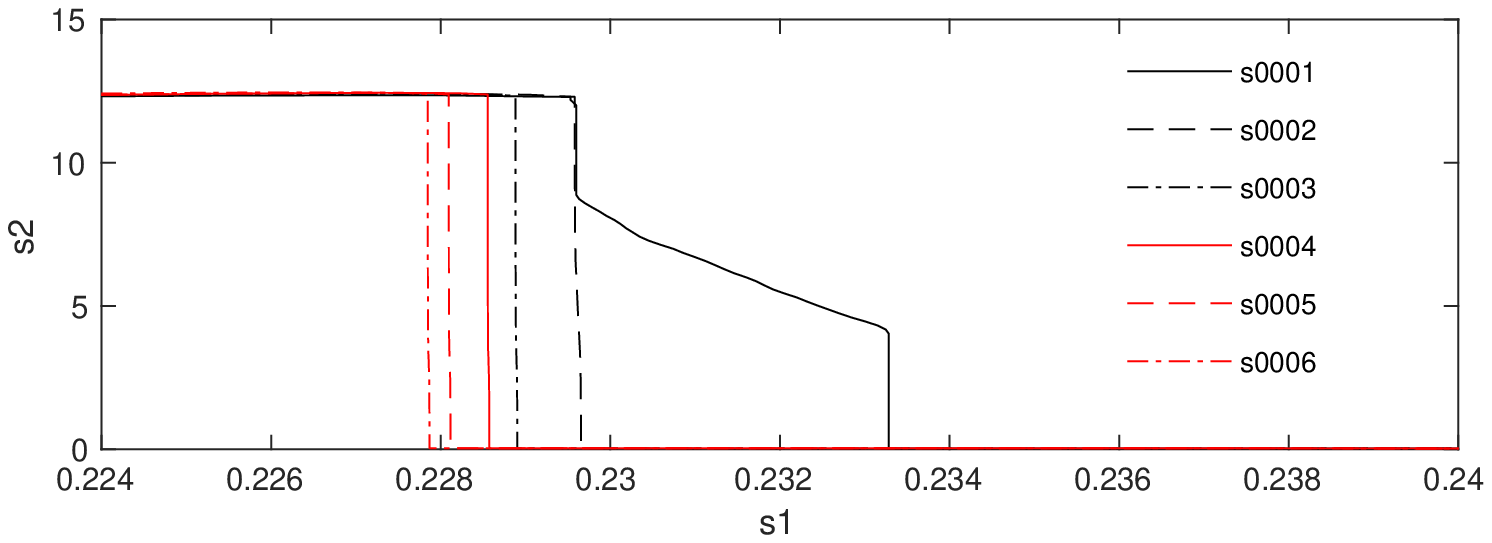}
	\end{psfrags}}
	\subcaptionbox{\label{fig:CT_timesteps_lp}}
	{\begin{psfrags}%
			\psfrag{s1}[tc][tc][1]{\color[rgb]{0,0,0}\setlength{\tabcolsep}{0pt}\begin{tabular}{c}$k$\end{tabular}}%
			\psfrag{s2}[bc][bc][1]{\color[rgb]{0,0,0}\setlength{\tabcolsep}{0pt}\begin{tabular}{c}$t^\rho_k$ \end{tabular}}%
			\psfrag{s0001}[cl][cl][0.5]{\color[rgb]{0,0,0}\setlength{\tabcolsep}{0pt}\begin{tabular}{c}$\alpha=4$\end{tabular}}%
			\psfrag{s0002}[cl][cl][0.5]{\color[rgb]{0,0,0}\setlength{\tabcolsep}{0pt}\begin{tabular}{c}$\alpha=6$\end{tabular}}%
			\psfrag{s0003}[cl][cl][0.5]{\color[rgb]{0,0,0}\setlength{\tabcolsep}{0pt}\begin{tabular}{c}$\alpha=8$\end{tabular}}%
			\psfrag{s0004}[cl][cl][0.5]{\color[rgb]{0,0,0}\setlength{\tabcolsep}{0pt}\begin{tabular}{c}$\alpha=10$\end{tabular}}%
			\psfrag{s0005}[cl][cl][0.5]{\color[rgb]{0,0,0}\setlength{\tabcolsep}{0pt}\begin{tabular}{c}$\alpha=12$\end{tabular}}%
			\psfrag{s0006}[cl][cl][0.5]{\color[rgb]{0,0,0}\setlength{\tabcolsep}{0pt}\begin{tabular}{c}$\alpha=14$\end{tabular}}%
			\includegraphics[scale=0.4]{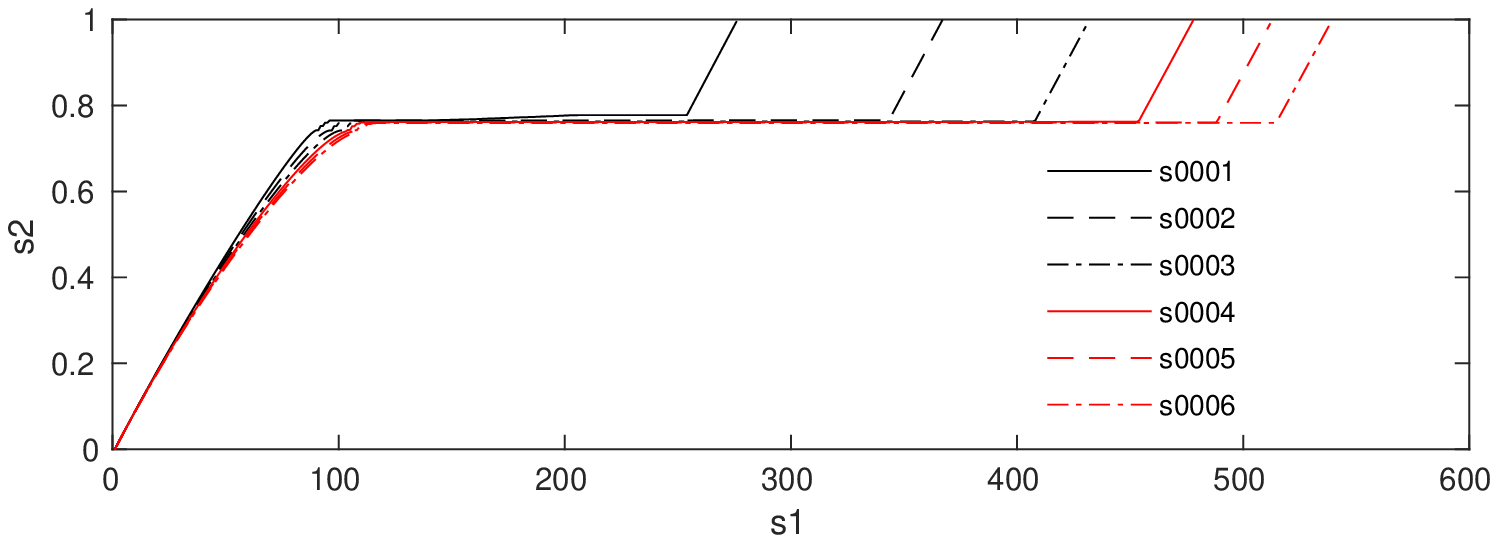}
	\end{psfrags}}
	\caption{Numerical analysis of a Compact-Tension test (CT) --- influence of the norm within the E\&M-algorithm: (left) zoom-in of the force-displacement-diagram; (right) evolution of the physical time as a function in terms of the E\&M-iteration $k$.\label{fig:CT_lpexperiment_1}}
\end{figure}

The structural response computed from the adaptive scheme is plotted in Fig.~\ref{fig:CT_lpexperiment_1} for different $L^\alpha$ norms.
To be more precise, $\alpha$ is varied between $4$ and $14$. According to the right diagram in Fig.~\ref{fig:CT_lpexperiment_1}, the arc-length increment $\rho=0.01$ has been chosen sufficiently small such that for all $\alpha$ brutal crack growth is predicted (the physical time is constant for several iterations). Furthermore the larger $\alpha$, the larger the number of E\&M-iterations. As a consequence, increasing $\alpha$ has the same effect as decreasing $\rho$ for the analyzed example. Hence, the left diagram in Fig.~\ref{fig:CT_lpexperiment_1} shows the same trend as the left diagram in Fig.~\ref{fig:CT_rhoexperiment}. Finally, it is noted that a small $\alpha<4$ as well as a too large $\alpha$ results in numerical problems.

\subsubsection{L-shaped plate}
\label{sec:L-shape}

\begin{figure}[htbp]
\centering
\subcaptionbox{}{
\begin{psfrags}%
\psfrag{s1}[cl][cl][1]{\color[rgb]{0,0,0}\setlength{\tabcolsep}{0pt}\begin{tabular}{c}$\bar u$\end{tabular}}%
\psfrag{l1}[cr][cr][1]{\color[rgb]{0,0,0}\setlength{\tabcolsep}{0pt}\begin{tabular}{c}$250\,\text{mm}$\end{tabular}}%
\psfrag{l2}[cm][cm][1]{\color[rgb]{0,0,0}\setlength{\tabcolsep}{0pt}\begin{tabular}{c}$250\,\text{mm}$\end{tabular}}%
\psfrag{s2}[cl][cl][1]{\color[rgb]{0,0,0}\setlength{\tabcolsep}{0pt}\begin{tabular}{c}$E = 25840$ MPa\end{tabular}}%
\psfrag{s3}[cl][cl][1]{\color[rgb]{0,0,0}\setlength{\tabcolsep}{0pt}\begin{tabular}{c}$\nu = 0.18$\end{tabular}}%
\psfrag{s4}[cl][cl][1]{\color[rgb]{0,0,0}\setlength{\tabcolsep}{0pt}\begin{tabular}{c}$g_c = 6.5\cdot 10^{-4}\frac{\text{N}}{\text{mm}} $\end{tabular}}%
\includegraphics[scale=1]{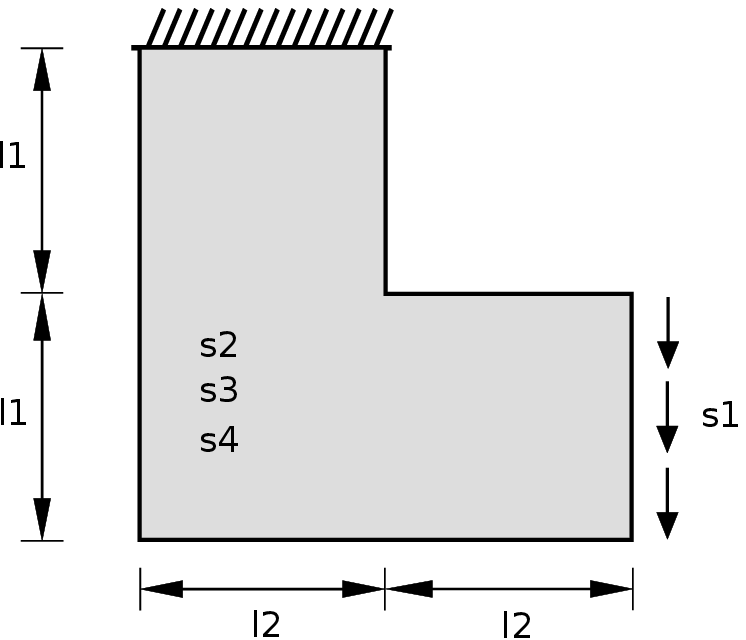}
\end{psfrags}}
\subcaptionbox{}{\includegraphics[scale=0.17]{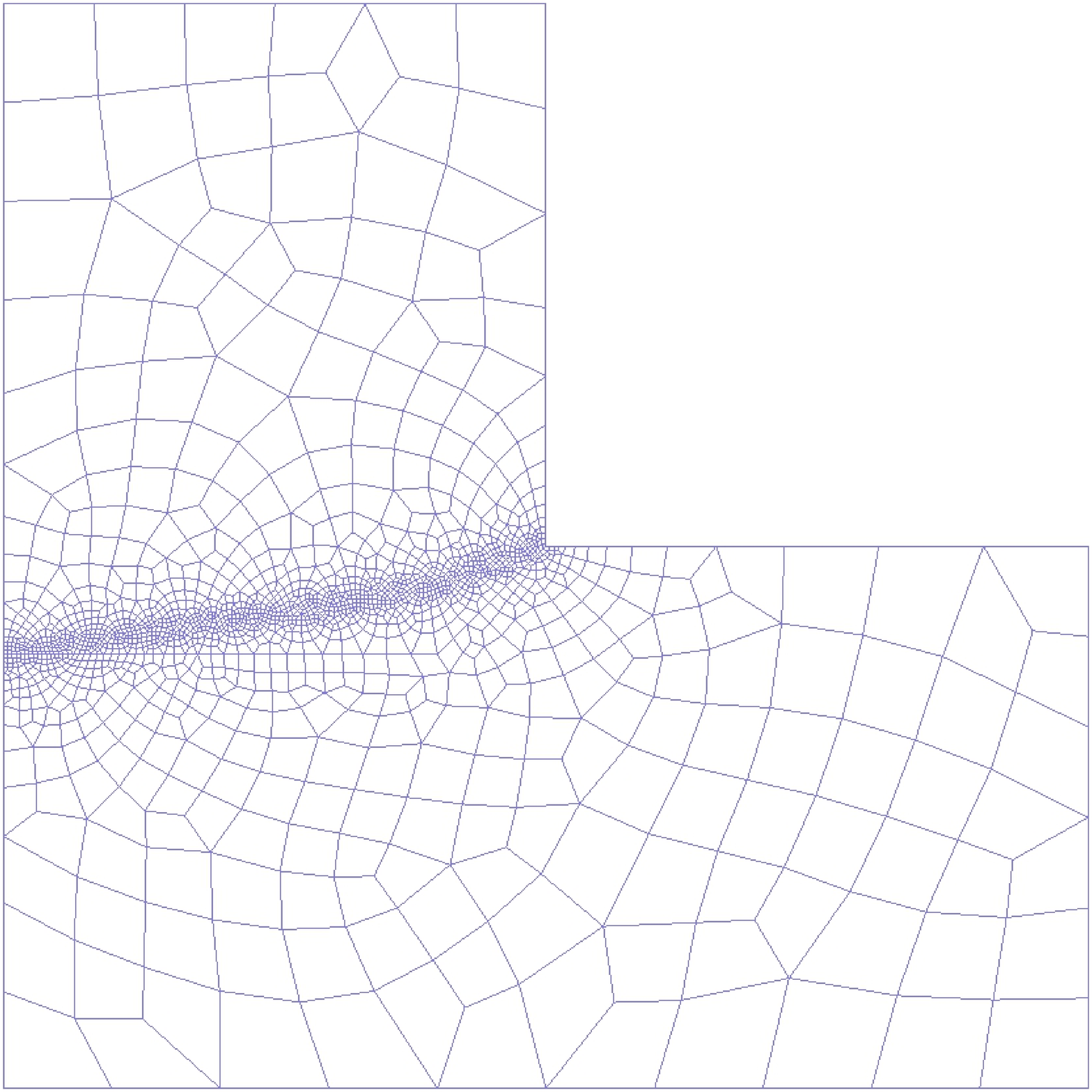}}
\caption{Numerical analysis of an L-shaped plate: (left) mechanical system; (right) finite element triangulation.}
\label{fig:LShape_Sketch}
\end{figure}

The second example is the L-shaped plate shown in Fig.~\ref{fig:LShape_Sketch}.
In contrast to the CT-specimen, the L-shaped plate is characterized by a more complex, curved crack, starting at the inner corner. The finite element mesh used within the numerical analysis is also given in Fig.~\ref{fig:LShape_Sketch}. Again, the mesh is refined in the region where the crack is expected to grow. Certainly and unlike the CT-specimen, this could lead to a small mesh bias. However, since convergence in the sense of the underlying finite element discretization is not the focus of this work and furthermore, the same mesh is chosen in all computation, this minor numerical artifact seems to be justifiable. The mesh consists of $1694$ quadrilateral elements with a maximal side-length of $50\,\text{mm}$ and a minimal length of $2\,\text{mm}$. Furthermore, plane strain conditions are assumed. The displacement-load is prescribed by the linear function $\bar u(t^\rho_k)=u_{max} \,\frac{t^\rho_k}{8.658}$ with the maximal time $T=8.658$. The arc-length increment is set to $\rho=0.08658$ such that, again, a minimum of $100$ time-steps is guaranteed and the length-scale parameter is chosen to $\theta=10$~mm.

\begin{figure}[htbp]
\centering
\subcaptionbox{\label{fig:Lshape_FUcurves_lp}}
{\begin{psfrags}%
\psfrag{s1}[tc][tc][1]{\color[rgb]{0,0,0}\setlength{\tabcolsep}{0pt}\begin{tabular}{c}$\bar u$ [mm]\end{tabular}}%
\psfrag{s2}[bc][bc][1]{\color[rgb]{0,0,0}\setlength{\tabcolsep}{0pt}\begin{tabular}{c}$F$ [N] \end{tabular}}%
\psfrag{s0001}[cl][cl][0.5]{\color[rgb]{0,0,0}\setlength{\tabcolsep}{0pt}\begin{tabular}{c}$\alpha=4$\end{tabular}}%
\psfrag{s0002}[cl][cl][0.5]{\color[rgb]{0,0,0}\setlength{\tabcolsep}{0pt}\begin{tabular}{c}$\alpha=6$\end{tabular}}%
\psfrag{s0003}[cl][cl][0.5]{\color[rgb]{0,0,0}\setlength{\tabcolsep}{0pt}\begin{tabular}{c}$\alpha=8$\end{tabular}}%
\psfrag{s0004}[cl][cl][0.5]{\color[rgb]{0,0,0}\setlength{\tabcolsep}{0pt}\begin{tabular}{c}$\alpha=10$\end{tabular}}%
\psfrag{s0005}[cl][cl][0.5]{\color[rgb]{0,0,0}\setlength{\tabcolsep}{0pt}\begin{tabular}{c}$\alpha=12$\end{tabular}}%
\psfrag{s0006}[cl][cl][0.5]{\color[rgb]{0,0,0}\setlength{\tabcolsep}{0pt}\begin{tabular}{c}$\alpha=14$\end{tabular}}%
\includegraphics[scale=0.4]{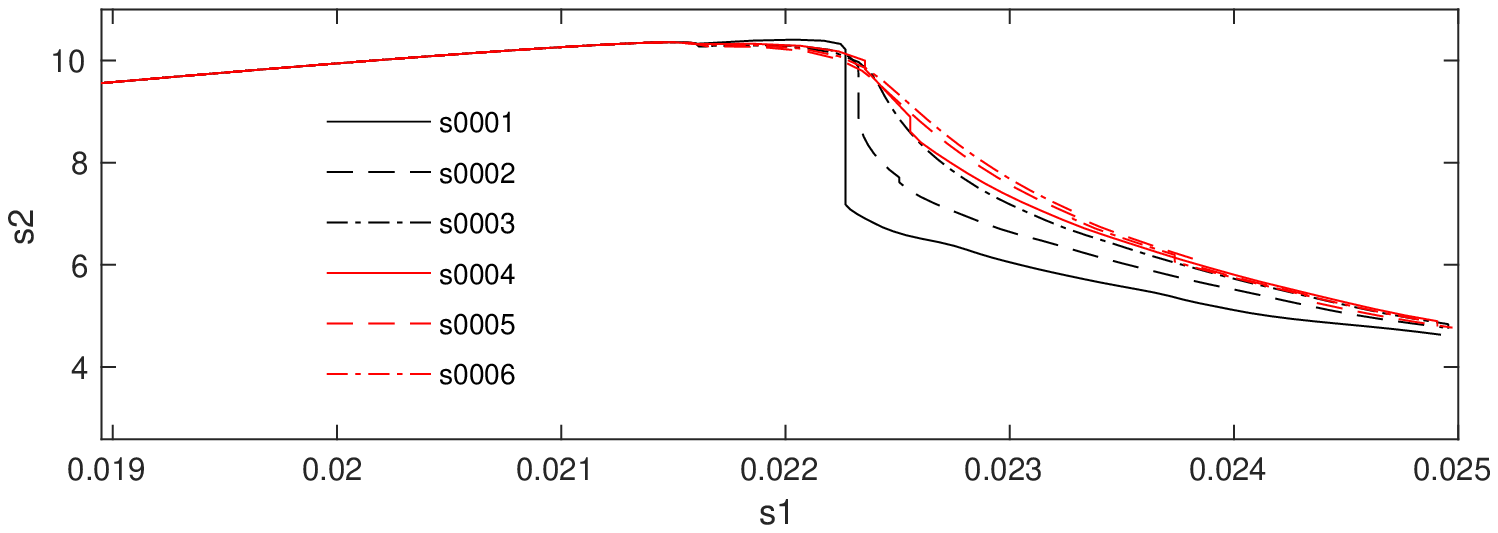}
\end{psfrags}}
\subcaptionbox{\label{fig:Lshape_timesteps_lp}}
{\begin{psfrags}%
\psfrag{s1}[tc][tc][1]{\color[rgb]{0,0,0}\setlength{\tabcolsep}{0pt}\begin{tabular}{c}$k$\end{tabular}}%
\psfrag{s2}[bc][bc][1]{\color[rgb]{0,0,0}\setlength{\tabcolsep}{0pt}\begin{tabular}{c}$t^\rho_k$ \end{tabular}}%
\psfrag{s0001}[cl][cl][0.5]{\color[rgb]{0,0,0}\setlength{\tabcolsep}{0pt}\begin{tabular}{c}$\alpha=4$\end{tabular}}%
\psfrag{s0002}[cl][cl][0.5]{\color[rgb]{0,0,0}\setlength{\tabcolsep}{0pt}\begin{tabular}{c}$\alpha=6$\end{tabular}}%
\psfrag{s0003}[cl][cl][0.5]{\color[rgb]{0,0,0}\setlength{\tabcolsep}{0pt}\begin{tabular}{c}$\alpha=8$\end{tabular}}%
\psfrag{s0004}[cl][cl][0.5]{\color[rgb]{0,0,0}\setlength{\tabcolsep}{0pt}\begin{tabular}{c}$\alpha=10$\end{tabular}}%
\psfrag{s0005}[cl][cl][0.5]{\color[rgb]{0,0,0}\setlength{\tabcolsep}{0pt}\begin{tabular}{c}$\alpha=12$\end{tabular}}%
\psfrag{s0006}[cl][cl][0.5]{\color[rgb]{0,0,0}\setlength{\tabcolsep}{0pt}\begin{tabular}{c}$\alpha=14$\end{tabular}}%
\includegraphics[scale=0.4]{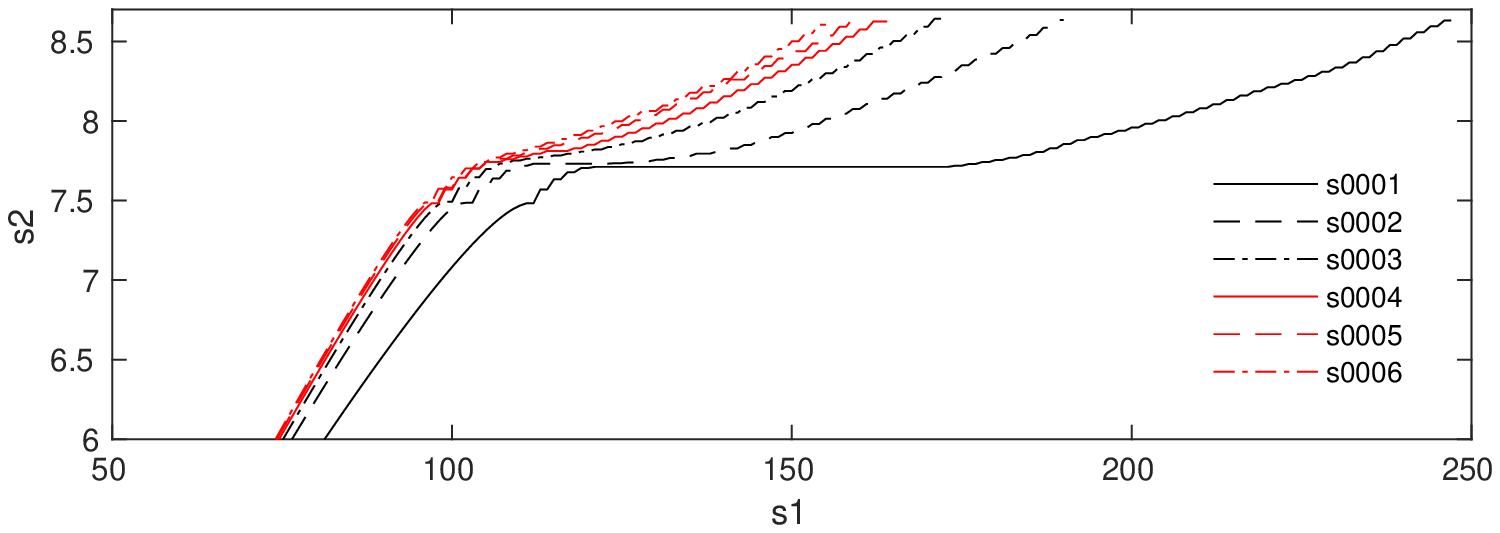}
\end{psfrags}}
\caption{Numerical analysis of an L-shaped plate: (left) zoom-in of the force-displacement curves; (right) evolution of the physical time as a function in terms of the E\&M-iteration $k$. \label{fig:Lshape_lpexperiment}}
\end{figure}

In contrast to the CT specimen analyzed in Subsection~\ref{sec:ct}, only the influence of the norm occurring within the adaptive scheme by Efendiev \& Mielke is analyzed here. However and in addition to the $L^\alpha$ norm the $H^1$ norm is also considered. The structural response computed from the adaptive algorithm is given in Fig.~\ref{fig:Lshape_lpexperiment} for different $L^\alpha$ norms.
While the left figure shows a zoom in of the force-displacement curve close to the point where brutal crack growth occurs, the evolution of the physical time as a function in terms of the E\&M-iteration $k$ is given in the right diagram --- for different $\alpha$ in each diagram. According to the right diagram in Fig.~\ref{fig:Lshape_lpexperiment}, the larger $\alpha$, the larger the time steps. It bears emphasis that this trend is in contrast to the Compact-Tension experiment. For this reason, estimating the influence of $\alpha$ on the structural response is very challenging. Due to the aforementioned trend and unlike the CT-specimen, the smaller $\alpha$, the more accurate is the prediction of brutal crack growth. Clearly, the relation just described is also influenced by the choice of the arc-length increment. To be more precise, brutal crack growth can also be predicted by using $\alpha=14$. However, the arc-length increment $\rho$ has to be chosen smaller in this case.

\begin{figure}[htbp]
\centering
\subcaptionbox{$t^\rho_{0}$=0}
{\includegraphics[scale=0.14]{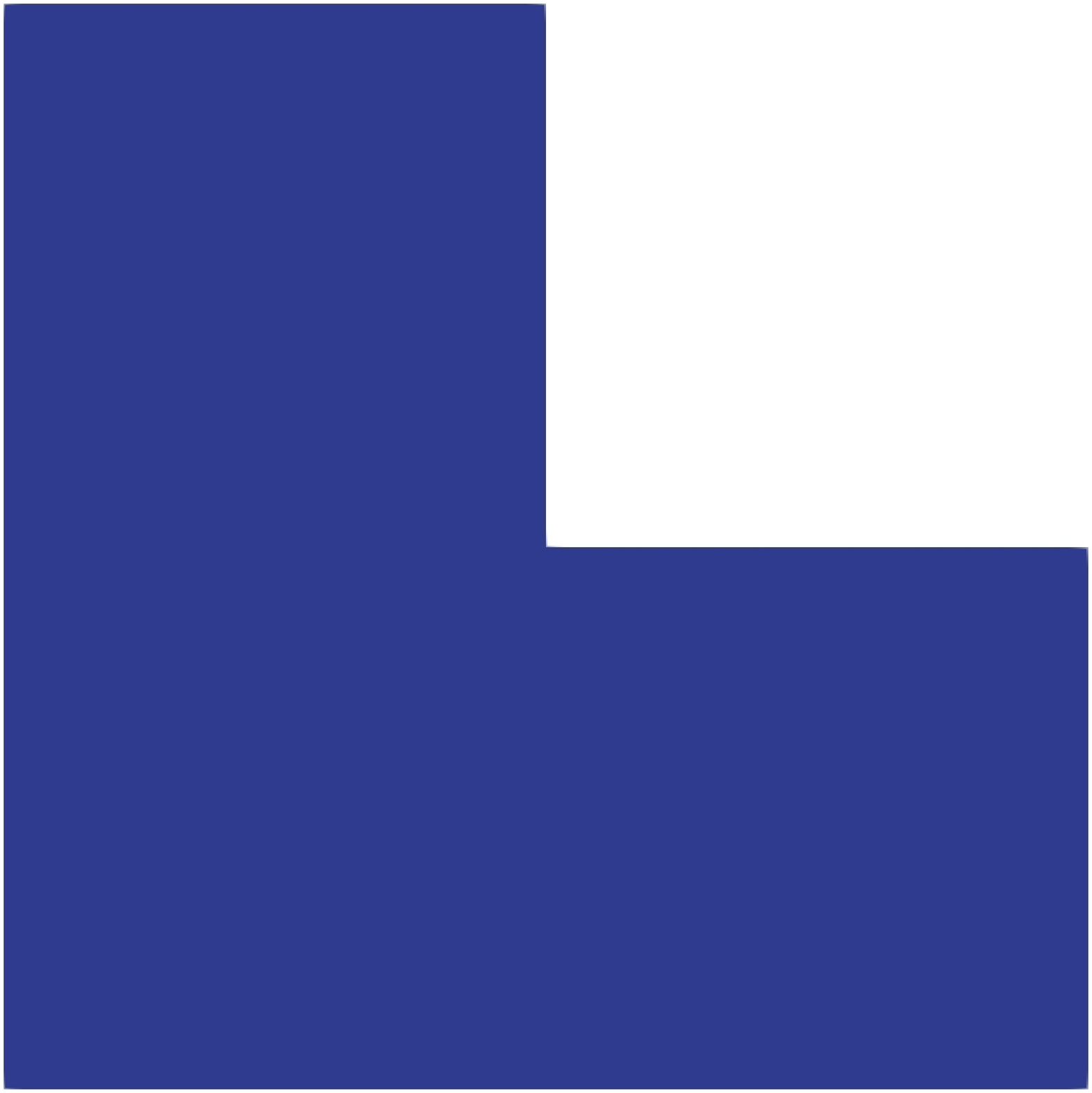}}
\subcaptionbox{$t^\rho_{121}=7.712$}
{\includegraphics[scale=0.14]{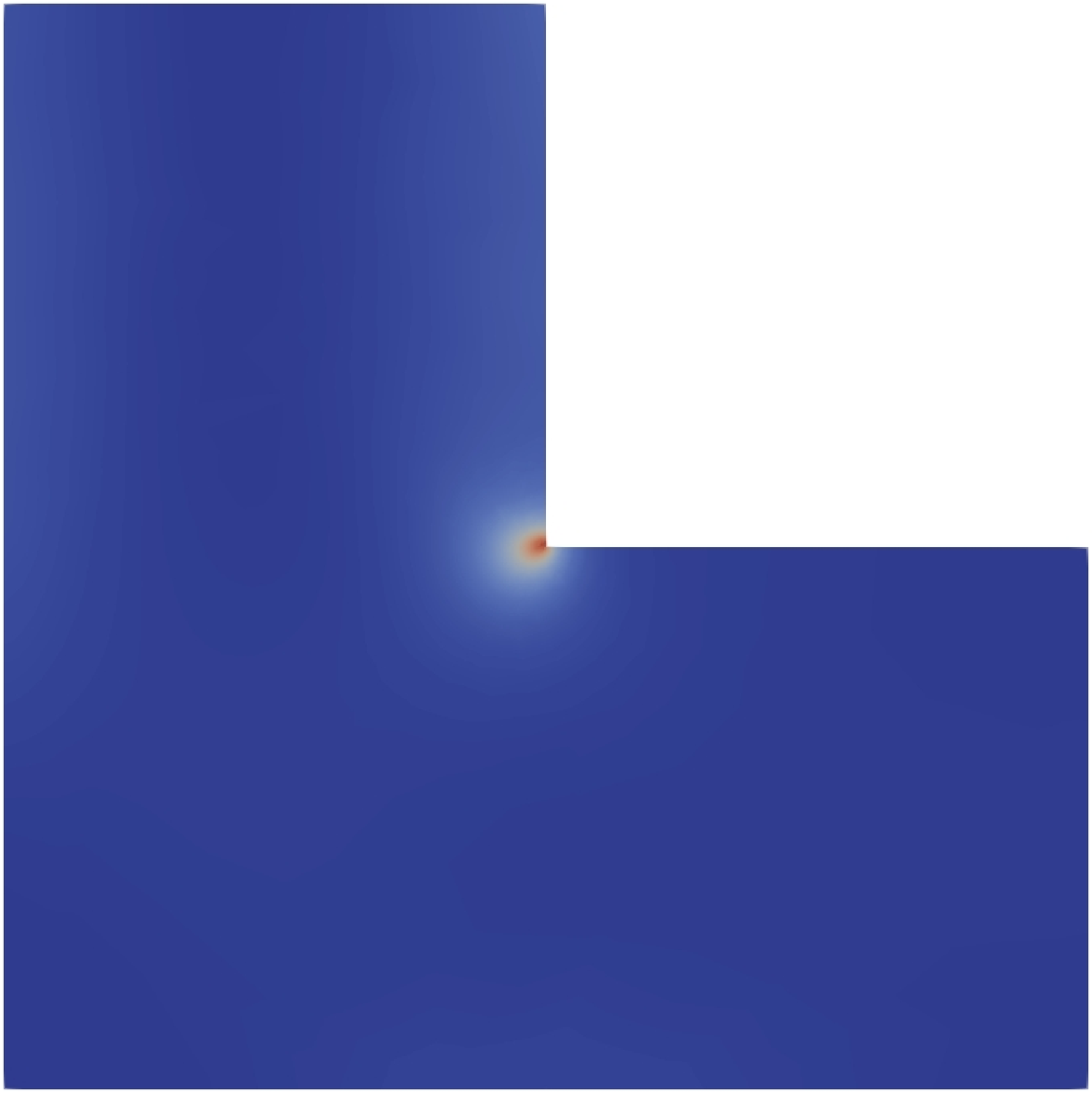}}
\subcaptionbox{$t^\rho_{138}=7.712$}
{\includegraphics[scale=0.14]{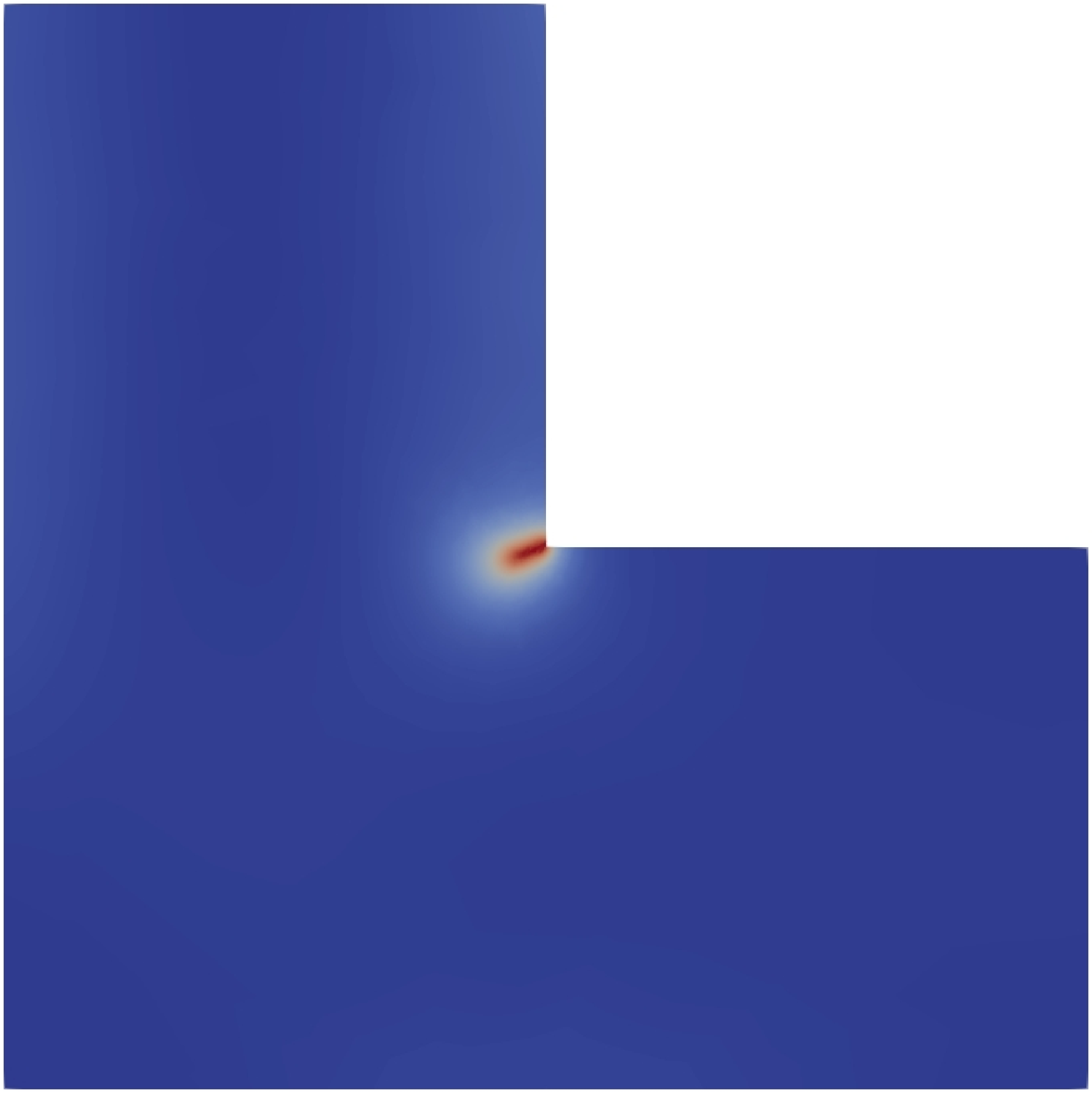}}
\subcaptionbox{$t^\rho_{156}=7.712$}
{\includegraphics[scale=0.14]{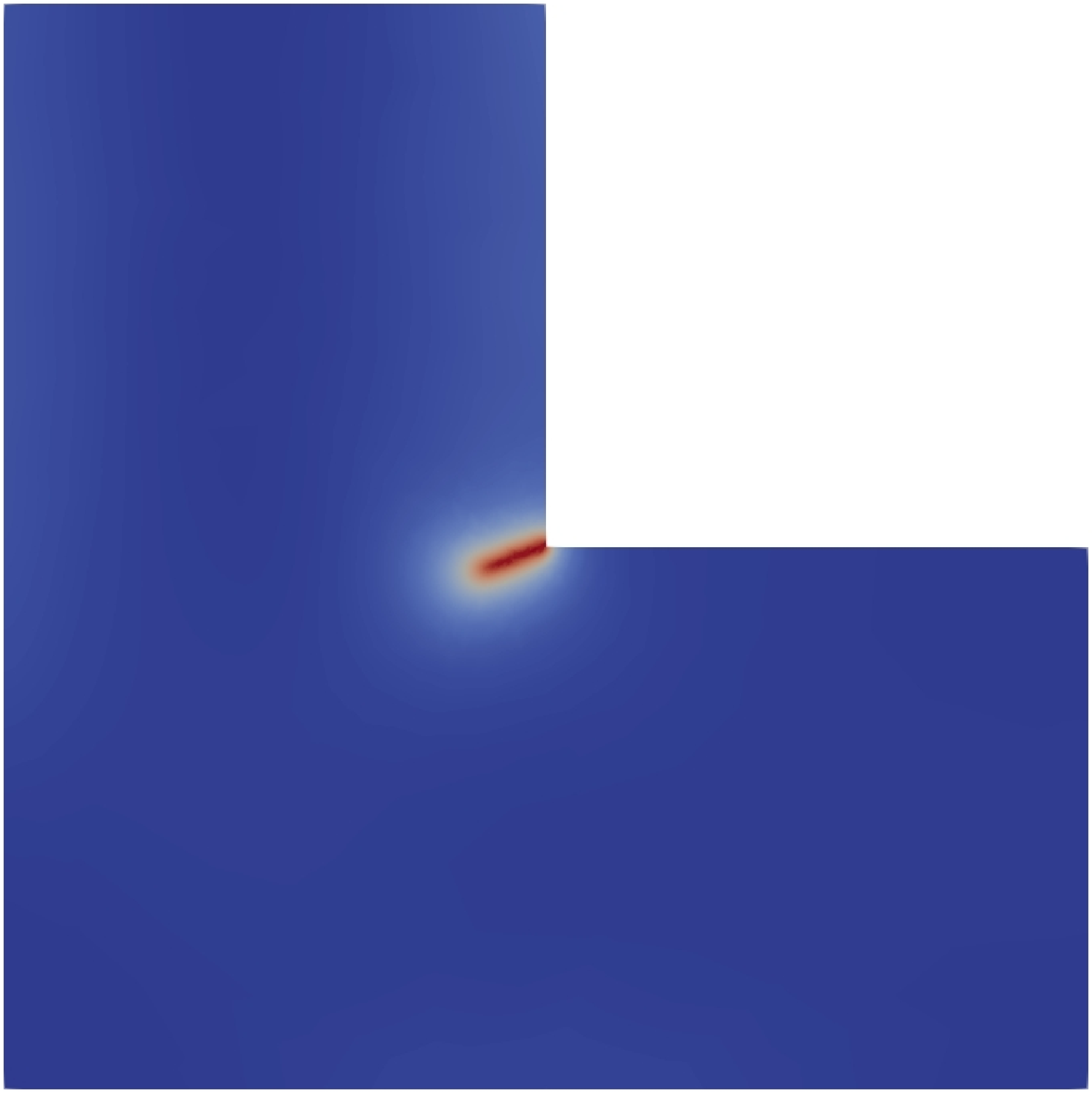}}
\subcaptionbox{$t^\rho_{174}=7.719$}
{\includegraphics[scale=0.14]{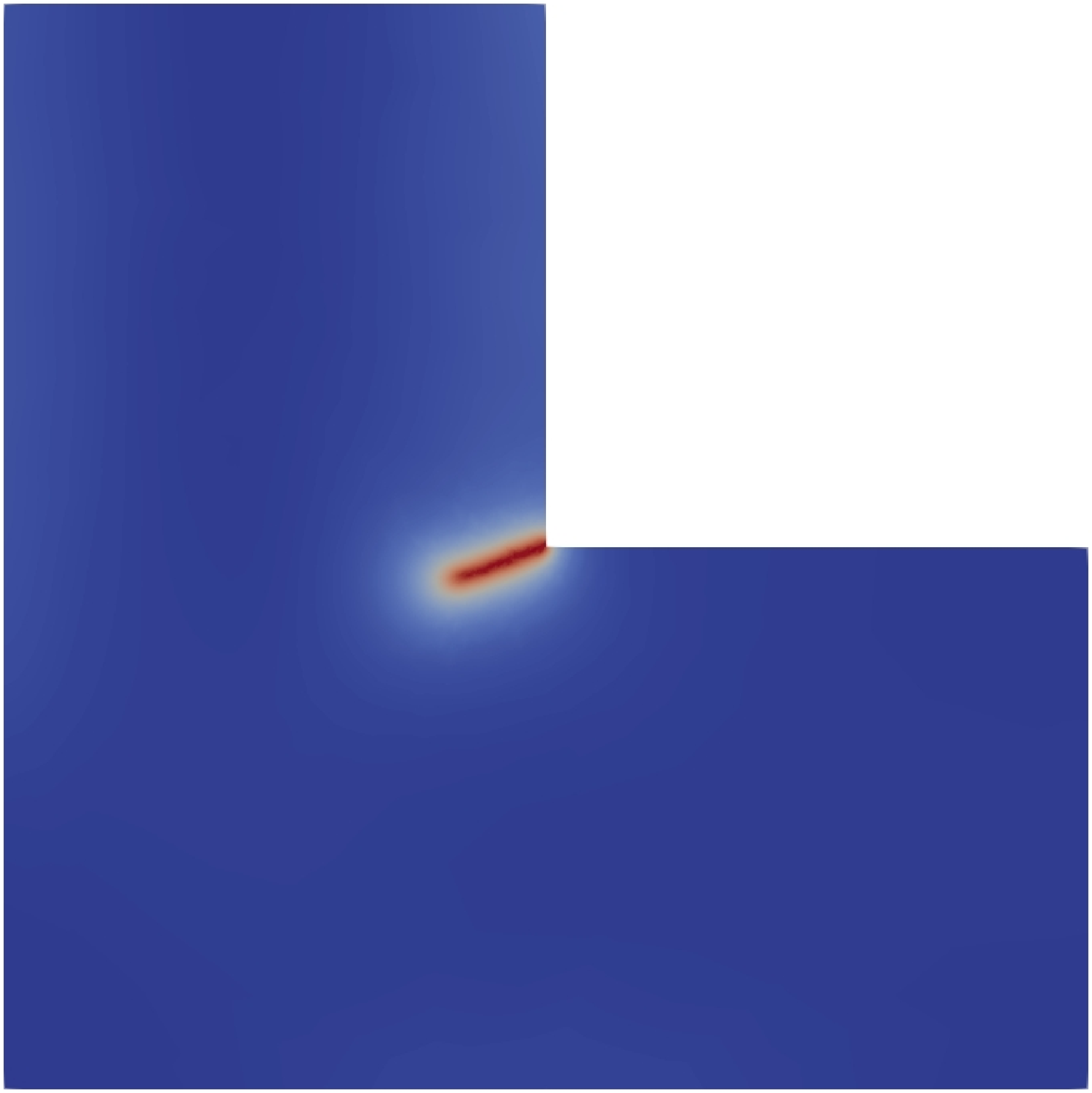}}
\subcaptionbox{$t^\rho_{196}=7.918$}
{\includegraphics[scale=0.14]{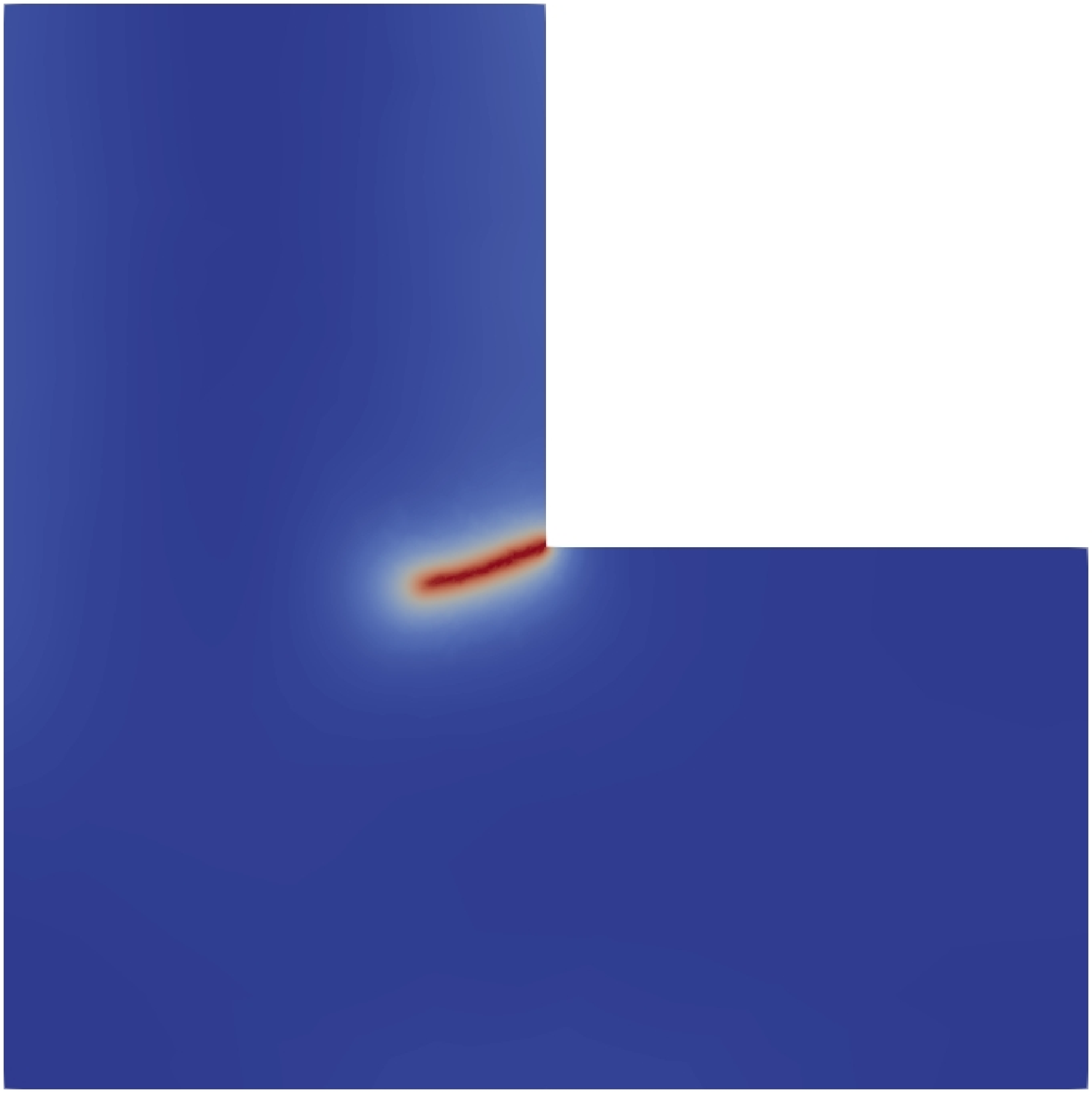}}
\subcaptionbox{$t^\rho_{221}=8.212$}
{\includegraphics[scale=0.14]{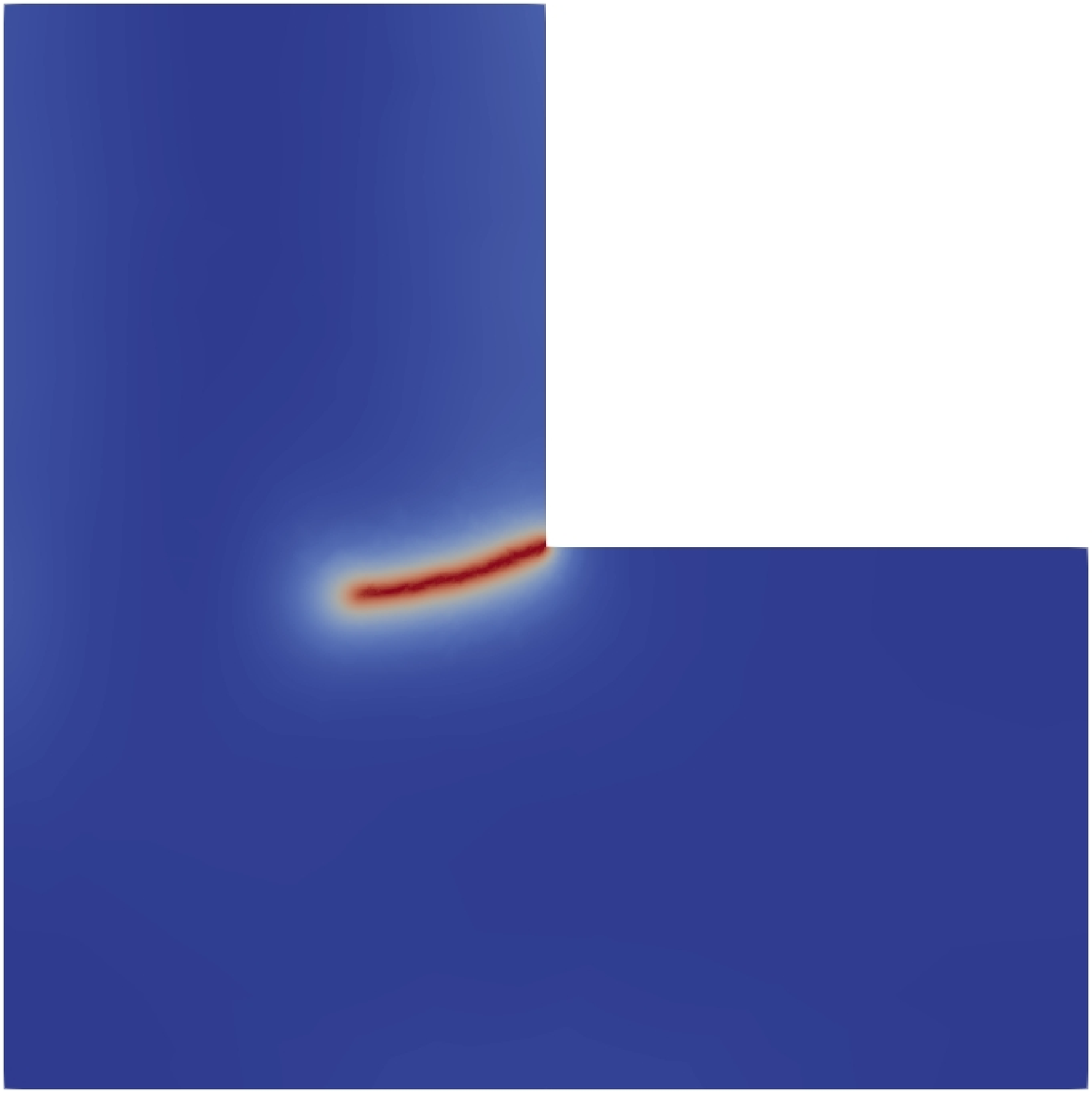}}
\subcaptionbox{$t^\rho_{247}=8.613$}
{\includegraphics[scale=0.14]{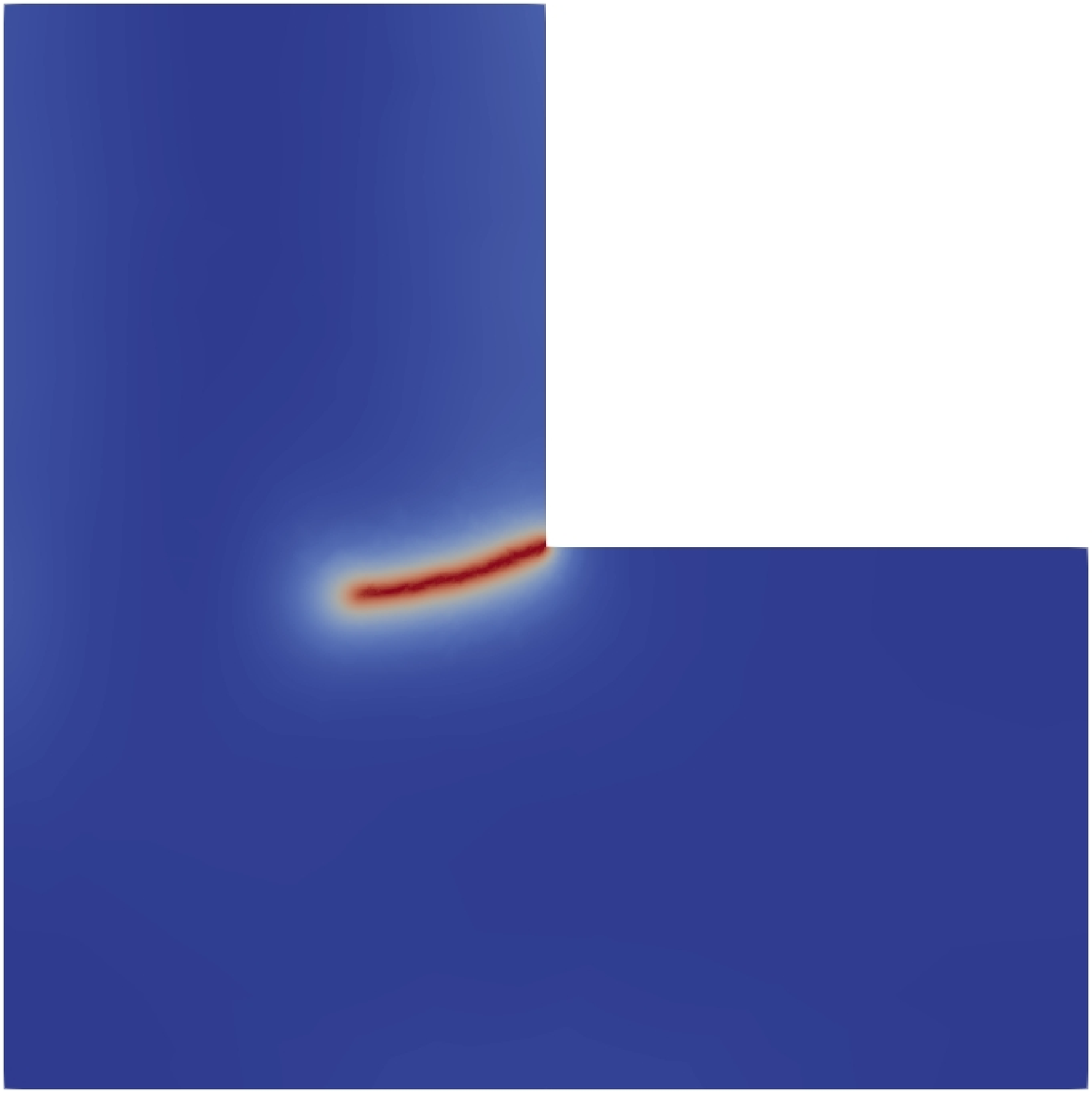}}
\caption{Numerical analysis of an L-shaped plate:  Distribution of the phase-field variable $z$ computed from the Efiendiev \& Mielke adaptive scheme combined with alternate minimization for $\rho=0.08658$ and $\alpha=4$. Brutal crack-growth occurs at $\bar u=0.228$ [mm].\label{fig:LShape_Plots}}
\end{figure}

According to the left diagram in Fig.~\ref{fig:Lshape_lpexperiment}, $\rho=0.08658$ and $\alpha=4$ correspond to brutal crack growth. The evolution of the phase-field parameter for this choice of model parameters is shown in Fig.~\ref{fig:LShape_Plots}.
It can be seen that brutal crack growth occurs at time $t=7.712$.

\begin{figure}[htbp]
	\centering
\begin{psfrags}%
\psfrag{s1}[tc][tc][1]{\color[rgb]{0,0,0}\setlength{\tabcolsep}{0pt}\begin{tabular}{c}$\bar u$ [mm]\end{tabular}}%
\psfrag{s2}[bc][bc][1]{\color[rgb]{0,0,0}\setlength{\tabcolsep}{0pt}\begin{tabular}{c}$F$ [N] \end{tabular}}%
\psfrag{4-normal}[cl][cl][0.5]{\color[rgb]{0,0,0}\setlength{\tabcolsep}{0pt}\begin{tabular}{c}$L^4\text{ norm},\,\rho=0.0866$\end{tabular}}%
\psfrag{4-small}[cl][cl][0.5]{\color[rgb]{0,0,0}\setlength{\tabcolsep}{0pt}\begin{tabular}{c}$L^4\text{ norm},\,\rho=0.0173$\end{tabular}}%
\psfrag{H1-normal}[cl][cl][0.5]{\color[rgb]{0,0,0}\setlength{\tabcolsep}{0pt}\begin{tabular}{c}$H^1$ norm, $\rho=0.0866$\end{tabular}}%
\psfrag{H1-large}[cl][cl][0.5]{\color[rgb]{0,0,0}\setlength{\tabcolsep}{0pt}\begin{tabular}{c}$H^1$ norm, $\rho=0.866$\end{tabular}}%
\includegraphics[scale=0.4]{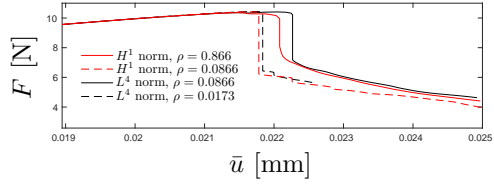}
\end{psfrags}
	\caption{Numerical analysis of an L-shaped plate: zoom-in of the force-displacement curves as computed by using the $H^1$ norm --- compared to results corresponding to the $L^\alpha$ norm.\label{fig:Lshape_H1Norm}}
\end{figure}

Finally, the computation is re-performed by using the $H^1$ norm. A zoom-in of the resulting force-displacement diagram can be seen in Fig.~\ref{fig:Lshape_H1Norm}.
Accordingly, also the $H^1$ norm is able to capture brutal crack growth --- depending on the choice of arc-length increment $\rho$.

It has been shown that all norms can capture brutal crack growth. However, this requires a sufficiently small arc-length increment which, in turn, depends on the chosen norm as well as on the analyzed structure.

	\section{Acknowledgment}

D. Knees and J. Mosler acknowledge the support in the framework of the DFG-Priority Programme SPP 2256  "Variational Methods for Predicting Complex Phenomena in Engineering Structures and Materials" funded by the Deutsche Forschungsgemeinschaft (DFG, German research foundation) within the project  "Rate-independent systems in solid mechanics: physical properties, mathematical analysis, efficient numerical algorithms" (Project-ID 441222077).
	
		\appendix
	\renewcommand{\thesubsection}{\Roman{subsection}}
	\section{Continuous dependence}
	
	\begin{lemma}[Lemma~2.5 from \cite{KRZ}]\label{lemmaA1}
		For $t\in\left[0,T\right]$ and $z\in\calZ$ let $u_{\text{min}}(t,z)$ denote the minimizer of $\calE(t,\cdot,z)$ over $\calU$. Moreover, let $p>2$ be given according to Theorem~\ref{thm:improved integrability} with $M=1$. Then 
		\begin{multline*}
			\exists C>0,\,\forall \tilde{p}\in\left[2,p\right),\,\forall t_1,t_2\in\left[0,T\right],\,\forall z_1,z_2\in\calZ_{\left[0,1\right]}:\quad\\
			\bnorm{u_{\text{min}}(t_2,z_2)-u_{\text{min}}(t_1,z_1)}_{\calU^{\tilde{p}}}
			\leq C\abs{t_2-t_1}+C\bnorm{z_2-z_1}_{L^{r_1}(\Omega)}
		\end{multline*}
		with $r_1=\frac{\tilde{p}p}{p-\tilde{p}}>\tilde{p}$. Observe that $r_1>\tilde{p}$.
	\end{lemma}
	\begin{lemma}[Variant of Lemma~A.2 from \cite{KN2017} including E\&M constraint]\label{lemmaA2EM}
		Let $\rho>0$ and $\beta\in\left[1,\infty\right]$ be arbitrary. Moreover, let $p>2$ be given according to Theorem~\ref{thm:improved integrability} with $M=1$.
		Let $\tilde{p}\in\left(2,p\right]$, $t\in\left[0,T\right]$, $u_1,u_2\in \calU^{\tilde{p}}$ and $z_0\in\calZ_{\left[0,1\right]}$. Furthermore, let
		\begin{align}\label{assume:lemmaA2}
			z_i&=\argmin\Set{\calE(t,u_i,z)+\calR(z-z_0)}{z\in\calZ,\,\norm{z-z_0}_{L^\beta(\Omega)}\leq\rho}\quad\text{for }i\in\{1,2\}.
		\end{align}
		Then there exists a constant $C>0$ such that for
		$r_2=\frac{\tilde{p}}{\tilde{p}-2}$ we have
		\begin{align*}
			\norm{z_1-z_2}_\calZ^2\leq C\bnorm{z_2-z_1}_{L^{r_2}(\Omega)}\bnorm{u_1-u_2}_{\calU^{\tilde{p}}}\Bigl(\bnorm{u_1}_{\calU^{\tilde{p}}}+\bnorm{u_2}_{\calU^{\tilde{p}}}\Bigr).
		\end{align*}
		Observe that $r_2>1$.
	\end{lemma}
	\begin{proof}
		We will denote with $K$ the set of functions $z\in\calZ$ that satisfy the Efendiev\&Mielke constraint as well as the irreversibility constraint, i.e.
		\begin{align*}
			K\coloneq \Set{z\in\calZ}{z\leq z_0\text{ a.e. in } \Omega,\, \norm{z-z_0}_{L^\beta(\Omega)}\leq \rho}.
		\end{align*}
		Let then $\chi_{_K}:\calZ\to\left[0,\infty\right]$ denote the indicator function of $K$, i.e.
		\begin{align*}
			\chi_{_K}(z)\coloneq \begin{cases}
				0,\quad z\in K\\
				\infty,\quad \text{otherwise.}
			\end{cases}
		\end{align*}
		Now we set for $t\in\left[0,T\right]$, $u\in\calU^{\tilde{p}}$ and $z\in\calZ$
		\begin{align*}
			\calF(t,u,z)\coloneq \calE(t,u,z)+\kappa\int_{\Omega}\left(z_0-z\right)\dd{x}
		\end{align*}
		such that we can reformulate minimality condition \eqref{assume:lemmaA2}, where
		\begin{align*}
			\calR(v)=\begin{cases}
				\kappa \int_{\Omega}\abs{v}\dd{x},\quad\text{if } v\leq 0 \text{ a.e. in }\Omega\\
				\infty,\quad\text{otherwise},
			\end{cases}
		\end{align*}
		as
		\begin{align*}
			z_i&=\argmin\Set{\calF(t,u_i,z)+\chi_{_K}(z)}{z\in\calZ}\quad\text{for }i\in\{1,2\}.
		\end{align*}
		This implies
		\begin{align*}
			0\in \D_z\calF(t,u_i,z_i)+\partial \chi_{_K}(z_i)\quad\text{for }i\in\{1,2\}
		\end{align*}
		and by definition of the subdifferential this is equivalent to
		\begin{align*}
			\forall w\in\calZ:\quad \chi_{_K}(w)\geq \chi_{_K}(z_i)+\langle-\D_z\calF(t,u_i,z_i),w-z_i\rangle \quad\text{for }i\in\{1,2\}
		\end{align*}
		respectively
		\begin{align*}
			\forall w\in K:\quad 0\geq \langle-\D_z\calF(t,u_i,z_i),w-z_i\rangle \quad\text{for }i\in\{1,2\}.
		\end{align*}
		By \eqref{assume:lemmaA2} it obviously holds $z_i\in K$ for $i\in\{1,2\}$. So if we test for $i=1$ with $w=z_2$ and for $i=2$ with $w=z_1$, we obtain
		\begin{align*}
			0\geq \langle-\D_z\calF(t,u_1,z_1),z_2-z_1\rangle\quad\text{and}\quad0\geq \langle-\D_z\calF(t,u_2,z_2),z_1-z_2\rangle.
		\end{align*}
		Thus summation leads to 
		\begin{align*}
			\langle \D_z\calF(t,u_2,z_2)-\D_z\calF(t,u_1,z_1),z_2-z_1\rangle\leq 0.
		\end{align*}
		Moreover, we have
		\begin{align*}
			\langle \D_z\calF(t,u_2,z_2)-\D_z&\calF(t,u_2,z_1),z_2-z_1\rangle \\
			&= \kappa\int_{\Omega}(z_2-z_1)^2+\abs{\nabla(z_2-z_1)}^2 \dd{x}+\underbrace{\int_{\Omega}(z_2-z_1)^2\bigl(\boldC\strain(u_2)\bigr)\colon\strain(u_2)\dd{x}}_{\geq 0}\\
			&\geq \kappa\norm{z_2-z_1}_{\calZ}^2.
		\end{align*}
		Since $z_0\in\calZ_{\left[0,1\right]}$, analogously to Lemma~\ref{lemma:be_AM} we get $z_1,z_2\in\calZ_{\left[0,1\right]}$. 
		Therefore, we obtain by 
		Hölder's inequality with $\frac{1}{r_2}+\frac{2}{\tilde{p}}=1$ that
		\begin{align*}
			\langle \D_z\calF(t,u_1,z_1)&-\D_z\calF(t,u_2,z_1),z_2-z_1\rangle\\
			&=\int_\Omega z_1\left(z_2-z_1\right) \Bigl(\bigl(\boldC\strain(u_1)\bigr)\colon\strain(u_1)- \bigl(\boldC\strain(u_2)\bigr)\colon\strain(u_2)\Bigr)\dd{x}\\
			&=\int_\Omega \underbrace{z_1}_{\in\calZ_{\left[0,1\right]}}\left(z_2-z_1\right) \bigl(\boldC\strain(u_1+u_2)\bigr)\colon\strain(u_1-u_2)\dd{x}\\
			&\leq C\bnorm{z_2-z_1}_{L^{r_2}(\Omega)}\bnorm{\strain(u_1-u_2)}_{L^{\tilde{p}}(\Omega)}\left(\bnorm{\strain(u_1)}_{L^{\tilde{p}}(\Omega)}+\bnorm{\strain(u_2)}_{L^{\tilde{p}}(\Omega)}\right)\\
			&\leq C\bnorm{z_2-z_1}_{L^{r_2}(\Omega)}\bnorm{u_1-u_2}_{\calU^{\tilde{p}}}\Bigl(\bnorm{u_1}_{\calU^{\tilde{p}}}+\bnorm{u_2}_{\calU^{\tilde{p}}}\Bigr).
		\end{align*}
		In summary this implies
		\begin{align*}
			\kappa&\norm{z_1-z_2}_\calZ^2\leq \langle \D_z\calF(t,u_2,z_2)-\D_z\calF(t,u_2,z_1),z_2-z_1\rangle\\
			&= \underbrace{\langle \D_z\calF(t,u_2,z_2)-\D_z\calF(t,u_1,z_1),z_2-z_1\rangle}_{\leq 0} +\langle \D_z\calF(t,u_1,z_1)-\D_z\calF(t,u_2,z_1),z_2-z_1\rangle\\
			&\leq  C\bnorm{z_2-z_1}_{L^{r_2}(\Omega)}\bnorm{u_1-u_2}_{\calU^{\tilde{p}}}\Bigl(\bnorm{u_1}_{\calU^{\tilde{p}}}+\bnorm{u_2}_{\calU^{\tilde{p}}}\Bigr).
		\end{align*}
	\end{proof}
		
	\section{Convex Analysis}\label{appendix:convexanalysis}
	Within this section we collect some convex analysis tools.
	\begin{definition}[Subdifferential]\label{defi:subdifferential}
		Let $\calQ$ be a normed vector space. Furthermore, let $f:\calQ\to\bbR_\infty$ be a convex function and let $u\in\calQ$. $f$ is called \textit{subdifferentiable} in $u$, if
		\begin{align}
			\partial f(u)\coloneq \Set{\xi\in\calQ^\ast}{f(v)\geq f(u)+\left\langle \xi, v-u\right\rangle \,\forall v\in\calQ}\neq\emptyset.
		\end{align}
		Every $\xi\in\partial f(u)$ is called a \textit{subgradient} of $f$ in $u$ and $\partial f(u)$ is called the \textit{subdifferential} of $f$ in $u$. Whenever it is not obvious which underlying space $\calQ$ is chosen to build the subdifferential $\partial f(u)$, we use the notation $\partial^\calQ f(u)$.
	\end{definition}
	For the next lemma we refer to \cite[Proposition~1.26; Proposition~1.8]{Phelps1989}.
	\begin{lemma}[Minimizers of convex functions]\label{lemma:minimizerofconvexfunctions}
		Let $\calQ$ be a convex set and let $f:\calQ\to\R_\infty$ be a convex function. Then
		\begin{align*}
			x^\ast=\argmin\Set{f(x)}{x\in\calQ}\hspace{7pt}\Leftrightarrow\hspace{7pt} f(x)\geq f(x^\ast)=f(x^\ast)+\langle 0,x-x^\ast\rangle\,\,\forall x\in\calQ\hspace{7pt}\Leftrightarrow\hspace{7pt} 0\in\partial f(x^\ast).
		\end{align*}
		If $f$ is additionally continuous and Gateaux differentiable in $x^\ast$ with Gateaux derivative $\D f(x^\ast)$, then it holds $\partial f(x^\ast)=\{\D f(x^\ast)\}$, thus in this case $x^\ast$ is a global minimizer of $f$ if and only if $x^\ast$ is a stationary point of $f$, i.e. $\D f(x^\ast)=0$.
	\end{lemma}
	A proof of the following subdifferential sumrule, also known as Moreau-Rockafellar theorem, can be found for instance in \cite[Theorem~9.5.4]{Attouch2014} or \cite[Prop~5.6, see also Def.~2.1 and Prop.~3.1 for the assumptions]{EkelandTemam1976}.
	\begin{lemma}[Subdifferential sum rule/Moreau-Rockafellar Theorem]\label{lemma:subdiffsumrule}
		Let $\calQ$ be a normed vector space. Furthermore, let $f_1,f_2:\calQ\to\bbR_\infty$ be convex and lower semicontinuous functions. Moreover, let there exist $x_0\in \dom(f_1)\cap \dom(f_2)$ such that $f_1$ is continuous in $x_0$. Here $\dom(f)\coloneq \Set{x\in\calQ}{f(x)<\infty}$ denotes the effective domain of a function $f$. Then we have
		\begin{align*}
			\forall x\in\calQ:\quad\partial(f_1+f_2)(x)=\partial f_1(x)+\partial f_2(x).
		\end{align*}
	\end{lemma}
	
	\begin{definition}[Convex conjugate]
		Let $\calQ$ be a normed vector space and let $f:\calQ\to\R_\infty$ be a proper function. The \textit{convex conjugate} or \textit{Legendre-Fenchel conjugate} of $f$ is the function \begin{align*}
			f^\ast:\calQ^\ast\to\R_\infty,\quad f^\ast(\xi)\coloneq \sup_{x\in\calQ}\left(\langle\xi,x\rangle-f(x)\right)
		\end{align*}
		Whenever it is not obvious which underlying space $\calQ$ is chosen to build the Legendre-Fenchel conjugate $f^\ast$, we use the notation $f^{\ast^\calQ}$.
	\end{definition}

	 A proof of the following lemma about Fenchel's inequality (also called Young-Fenchel-inequality) and Fenchel's identity (also called Fenchel extremality relation) can be found for instance in \cite[Proposition~9.5.1, Theorem~9.5.1]{Attouch2014} or \cite[3.3.1 and 4.2.1 Hilfssatz~1]{IT1979}.
	
	\begin{lemma}[Fenchel's inequality and identity]\label{fenchel}
		Let $\calQ$ be a real reflexive Banach space and let $\calR:\calQ\to\R_\infty$ be a proper, convex and lower semicontinuous function. Then the \textit{Fenchel inequality}
		\begin{align*}
			\forall v\in\calQ,\,\forall\xi\in\calQ^\ast:\quad \calR(v)+\calR^\ast(\xi)\geq \langle\xi,v\rangle
		\end{align*}
		holds and we have the following equivalent characterizations of the \textit{Fenchel identity}
		\begin{align*}
			\calR(v)+\calR^\ast(\xi)=\langle\xi,v\rangle\quad\Leftrightarrow\quad \xi\in\partial\calR(v)\quad\Leftrightarrow\quad v\in\partial\calR^\ast(\xi).
		\end{align*}
	\end{lemma}

		\begin{definition}[Positively homogeneous of degree one]\label{defi:poshomogeneousdegone}
		Let $\calQ$ be a real vector space.
		A function $\calR:\calQ\to\bbR_\infty$ is called \textit{positively homogeneous of degree one} or shorter \textit{positively 1-homogeneous}, if
		\begin{align}
			\calR(0)=0\quad\text{and}\quad\forall \lambda> 0,\,\forall u\in\calQ:\quad\calR(\lambda u)=\lambda \calR(u).
		\end{align}
	\end{definition}
	
	\begin{lemma}[Subdifferentials and convex conjugates of positively 1-homogeneous functions]\label{prop:subdifferentialpos1homogen}
		Let $\calQ$ be a real reflexive Banach space and let $\calR:\calQ\to\left[0,\infty\right]$ be a proper, convex and lower semicontinuous function, which is additionally positively homogeneous of degree one. Then the subdifferentials satisfy
		\begin{align}\label{prop:sdphdoSubdiffin0}
			\partial \calR (0) &=\Set{\xi\in\calQ^\ast}{\calR(v)\geq \left\langle \xi,v\right\rangle \,\forall v\in\calQ},\\\label{prop:sdphdoSubdiff}
			\forall u\in\calQ:\quad\partial \calR(u) &=\Set{\xi\in\partial \calR(0)}{\calR(u)= \left\langle \xi,u\right\rangle}
		\end{align}
		and $\partial\calR(0)$ is closed and convex with $0\in\partial\calR(0)$.
		
		Furthermore, the convex conjugate $\calR^\ast:\calQ^\ast\to\R_\infty$ equals the indicator function of the subdifferential $\partial\calR(0)$, i.e.
		\begin{align}\label{prop:convexconjugatepos1hom}
			\forall\xi\in\calZ^\ast:\quad\calR^\ast(\xi)=\chi_{_{\partial\calR(0)}}(\xi)=
			\begin{cases}
				0, &\text{if }\xi\in\partial\calR(0)\\
				\infty, &\text{otherwise.}
			\end{cases}
		\end{align}
	\end{lemma}
	\begin{proof}
		By Definition~\ref{defi:poshomogeneousdegone} we have $\calR(0)=0$, so \eqref{prop:sdphdoSubdiffin0} follows immediately from  Definition~\ref{defi:subdifferential}. The mapping $\calQ^\ast\ni\xi\mapsto \left\langle \xi,v\right\rangle$ is linear and continuous for every $v\in\calQ$. Therefore, the set $\partial\calR(0)$ is convex and closed. Since $\calR(v)\geq 0$ for every $v\in\calQ$, we have $0\in\partial\calR(0)$. 
		It remains to show \eqref{prop:sdphdoSubdiff}. Let $u\in\calQ$ be arbitrary.
		
		\glqq$\subset$\grqq: Let $\xi\in\partial\calR(u)$. 
		If we test with $v\coloneq 0$ we obtain from Definition~\ref{defi:subdifferential} and with using the positive 1-homogeneity of $\calR$
		\begin{align*}
			0=\calR(0)\geq \calR(u)+\left\langle \xi,-u\right\rangle\text{, i.e. }\left\langle \xi,u\right\rangle\geq \calR(u).
		\end{align*}
		Whereas if we test with $v\coloneq 3u$ we obtain
		\begin{align}\label{proof:sdphdo1}
			3\calR(u)=\calR(v)\geq \calR(u)+ \left\langle \xi,v-u\right\rangle =\calR(u)+\left\langle \xi,2u\right\rangle,\quad \text{i.e.}\quad \calR(u)\geq\left\langle \xi,u\right\rangle.
		\end{align}
		Thus $\calR(u)=\langle\xi,u\rangle$.
		Next, together with $\xi\in\partial\calR(u)$, we also obtain
		\begin{align*}
			\forall v\in\calQ:\quad \calR(v)\geq\calR(u)+\langle\xi,v-u\rangle=\langle \xi,v\rangle,
		\end{align*}
		and hence, by \eqref{prop:sdphdoSubdiffin0}, $\xi\in\partial\calR(0)$.
		
		\glqq$\supset$\grqq: Let $\xi\in\partial\calR(0)$ with $\calR(u)=\left\langle \xi,u\right\rangle$. Then for every $v\in\calQ$ we have
		\begin{align*}
			\calR(v)\geq \left\langle \xi,v\right\rangle =\left\langle \xi,v\right\rangle+\underbrace{\calR(u)-\left\langle \xi,u\right\rangle}_{=0}=\calR(u)+\left\langle \xi,v-u\right\rangle.
		\end{align*}
		For the proof of \eqref{prop:convexconjugatepos1hom} first notice that, since $\calR$ is positively $1$-homogeneous, we have
		\begin{align*}
			\forall\xi\in\calQ^\ast,\,\forall\lambda>0:\quad \calR^\ast(\xi)&=\sup_{x\in\calQ}\left(\langle\xi,x\rangle-\calR(x)\right)=\sup_{x\in\calQ}\left(\langle\xi,\lambda x\rangle-\calR(\lambda x)\right)\\
			&\hspace{-21pt}\overset{\calR\text{ pos. 1-hom.}}{=}\lambda \sup_{x\in\calQ}\left(\langle\xi,x\rangle-\calR(x)\right)=\lambda \calR^\ast(\xi).
		\end{align*}
		This implies that $\calR^\ast(\xi)\in\left\{0,\infty\right\}$ for all $\xi\in\calZ^\ast$. Secondly, by the equivalent characterizations of the Fenchel identity from Lemma~\ref{fenchel} it holds
		\begin{align*}
			\xi\in\partial\calR(0)\quad\Leftrightarrow\quad \underbrace{\calR(0)}_{=0}+\calR^\ast(\xi)=\langle \xi, 0\rangle=0\quad \Leftrightarrow\quad \calR^\ast(\xi)=0.
		\end{align*}
		Thus we obtain \eqref{prop:convexconjugatepos1hom}.
	\end{proof}

	For the inf-convolution formula we refer to \cite[Theorem~9.4.1]{Attouch2014} or \cite[p.~163-164, Satz~1]{IT1979}.
	
	\begin{lemma}[Inf-convolution formula]\label{infconvolution}
		Let $f_1,f_2:\calQ\to\bbR_\infty$ be as in Lemma~\ref{lemma:subdiffsumrule}. Then the following equality holds
		\begin{align*}
			\forall\xi\in\calQ^\ast:\quad\left(f_1+f_2\right)^\ast(\xi)=\inf_{\zeta\in\calQ^\ast}\left(f_1^\ast(\zeta)+f_2^\ast(\xi-\zeta)\right).
		\end{align*}
	\end{lemma}
	
	\begin{lemma}[Modification of Lemma~A.1 in \cite{Knees2018}]\label{lemmaA1knees2018}
		Let $\calZ$ and $\calV$ be the spaces from section~\ref{sec:spaces}. Further let $\rho>0$ be given and let $\calR:\calZ\to\R_\infty$ be a proper positively 1-homogeneous convex lower semicontinuous function.
		Define $\psi_\rho:\calZ\to\R_\infty$ as
		\begin{align*}
			\psi_\rho\coloneq \calR+\chi_\rho\quad\text{with indicator function}\quad \chi_\rho:\calZ\to\R_\infty\quad\text{given by}\quad\chi_\rho(v)=\begin{cases}
				0, & \text{ if } \norm{v}_\calV\leq \rho\\
				\infty, &\text{ otherwise}.
			\end{cases}
		\end{align*}
		Then its convex conjugate $\psi_\rho^{*^\calZ}:\calZ^\ast\to\R_\infty$ is given by
		\begin{align*}
			\forall \zeta\in\calZ^\ast:\quad\psi_\rho^{*^\calZ}(\zeta)=\rho\,\dist_{\calV^\ast}\bigl(\zeta,\partial^\calZ\calR(0)\bigr),
		\end{align*}
		where $\dist_{\calV^\ast}\big(\zeta,\partial^\calZ\calR(0)\big)=\infty$, if there exists no $\xi\in\partial^\calZ\calR(0)$ such that $\zeta-\xi\in\calV^*$.
	\end{lemma}
	\begin{proof}
		We follow \cite[Lemma~A.1]{Knees2018}. By the inf-convolution formula (see Lemma~\ref{infconvolution}; notice that $0\in\dom(\calR)\cap\dom(\chi_\rho)$ and $\chi_\rho$ is continuous in $0$)  we obtain
		\begin{align*}
			\forall\zeta\in\calZ^\ast:\quad\psi_\rho^{*^\calZ}(\zeta)=\inf_{\xi\in\calZ^*}\Big(\calR^{*^\calZ}(\xi)+\chi_\rho^{*^\calZ}(\zeta-\xi)\Big).
		\end{align*}
		First observe that for all $\xi\in\calZ^\ast$
		\begin{align*}
			\chi_\rho^{*^\calZ}(\xi)=\begin{cases}
				\rho\norm{\xi}_{\calV^*},&\text{ if } \xi\in\calV^*\\
				\infty,&\text{ otherwise}.
			\end{cases}
		\end{align*}
		This can be seen as follows: Assume that $\xi\in\calZ^\ast$ is chosen such that
		\begin{align*}
			\chi_\rho^{*^\calZ}(\xi)=\sup_{z\in\calZ}\big(\langle\xi,z\rangle-\chi_\rho(z)\big)=\sup_{\substack{z\in\calZ\\ \norm{z}_\calV\leq \rho}}\langle \xi, z\rangle<\infty,
		\end{align*}
		i.e. there exists $c>0$ such that for all $z\in\calZ$ with $\norm{z}_\calV\leq\rho$ it holds $\abs{\langle \xi,z\rangle}\leq c$. By testing with $z\coloneq \frac{\rho \,v}{\norm{v}_\calV}$ we obtain for any $v\in\calZ\setminus\left\{0\right\}$
		\begin{align*}
			\abs{\langle \xi,v\rangle}\leq \frac{c}{\rho}\norm{v}_\calV.
		\end{align*}
		Therefore, the density of $\calZ$ in $\calV$ allows us to extend $\xi\in\calZ^\ast$ in a unique way to an element of $\calV^\ast$. Thus we are actually in the case $\xi\in\calV^\ast$. Since also $\Set{z\in\calZ}{\norm{z}_\calV\leq\rho}$ lies dense in $\Set{v\in\calV}{\norm{v}_\calV\leq\rho}$, we finally obtain
		\begin{align*}
			\forall \xi\in\calV^\ast:\quad\chi_\rho^{*^\calZ}(\xi)=\sup_{\substack{v\in\calZ\\ \norm{\frac{v}{\rho}}_\calV\leq 1}}\rho\,\Big\langle\xi,\frac{v}{\rho}\Big\rangle= \sup_{\substack{v\in\calV\\ \norm{\frac{v}{\rho}}_\calV\leq 1}}\rho\,\Big\langle\xi,\frac{v}{\rho}\Big\rangle=\rho\,\norm{\xi}_{\calV^*}.
		\end{align*}
		Secondly, by \eqref{prop:convexconjugatepos1hom} we have
		\begin{align*}
			\forall \xi\in\calZ^\ast:\quad\calR^{*^\calZ}(\xi)=\chi_{_{\partial^\calZ\calR(0)}}(\xi)=
			\begin{cases}
				0, &\text{if } \xi\in\partial^\calZ\calR(0)\\
				\infty, &\text{ otherwise}.
			\end{cases}
		\end{align*}
		In summary we thus have for all $\zeta\in\calZ^\ast$
		\begin{align*}
			\psi_\rho^{*^\calZ}(\zeta)=\inf_{\xi\in\calZ^*}\left(\calR^{*\calZ}(\xi)+\chi_\rho^{*^\calZ}(\zeta-\xi)\right)=\inf_{\substack{\xi\in\partial^\calZ\calR(0)\\
					\zeta-\xi\in\calV^\ast}}\rho\,\norm{\zeta-\xi}_{\calV^*}=\rho\,\dist_{\calV^\ast}\bigl(\zeta,\partial^\calZ\calR(0)\bigr),
		\end{align*}
		where $\inf\emptyset=\infty$ respectively $\dist_{\calV^\ast}\big(\zeta,\partial^\calZ\calR(0)\big)=\infty$, if there exists no $\xi\in\partial^\calZ\calR(0)$ such that $\zeta-\xi\in\calV^*$.
	\end{proof}
	Within the proof of Lemma~\ref{lemmaA1knees2018} we showed that $\chi_\rho^{\ast^\calZ}(\xi)=\chi_\rho^{\ast^\calV}(\xi)$ for every $\xi\in\calV^\ast$. Analogously one can show the following identity concerning subdifferentials instead of convex conjugates of $\chi_\rho$, cf. \cite[Lemma~A.3.10]{Sievers2020}.
	\begin{lemma}[Lemma~A.3.10 from \cite{Sievers2020}]\label{lemmaA310}
		Let $\calZ$ and $\calV$ be the spaces from section~\ref{sec:spaces}.
		Further let $\rho>0$ be given and let $\psi_\rho:\calZ\to\R_\infty$ be as in Lemma~\ref{lemmaA1knees2018}. For all $v\in\calZ$ it holds 
		\begin{align*}
			\partial^\calZ\chi_\rho(v)=\partial^\calV\chi_\rho(v).
		\end{align*}
		This is to be understood as follows:
		\begin{align*}
			\xi\in\partial^\calZ\chi_\rho(v)\subset\calZ^\ast\quad &\Rightarrow\quad \xi \text{ can be uniquely extended to }\tilde{\xi}\in\calV^\ast\text{ and }\tilde{\xi}\in\partial^\calV\chi_\rho(v)
		\shortintertext{and}
			\xi\in\partial^\calV\chi_\rho(v)\subset\calV^\ast\quad &\Rightarrow\quad \xi\vert_\calZ\in\calZ^\ast\text{ and }\xi\vert_\calZ\in\partial^\calZ\chi_\rho(v).
		\end{align*}
	\end{lemma}
	\begin{proof}
		The inclusion $\partial^\calV\chi_\rho(v)\subset\partial^\calZ\chi_\rho(v)$ is a direct consequence from $\calZ\hookrightarrow\calV$. For the opposite inclusion let $\xi\in\partial^\calZ\chi_\rho(v)$ be arbitrary, then in particular $\chi_\rho(v)=0$ and $\norm{v}_\calV\leq\rho$. By the definition of $\chi_\rho$ and Definition~\ref{defi:subdifferential} of the subdifferential this yields
		\begin{align*}
			\forall z\in\calZ\text{ with } \norm{z}_\calV\leq \rho:&\quad \langle\xi,z\rangle\leq\langle\xi,v\rangle\leq \norm{\xi}_{\calZ^\ast}\norm{v}_\calZ.
		\end{align*}
		Thus by testing with $z=\frac{\rho\,w}{\norm{w}_\calV}$, we obtain for any $w\in\calZ\setminus\left\{0\right\}$
		\begin{align*}
			\langle\xi,w\rangle\leq\left(\frac{1}{\rho}\norm{\xi}_{\calZ^\ast}\norm{v}_\calZ\right)\norm{w}_\calV.
		\end{align*}
		Therefore, the density of $\calZ$ in $\calV$ allows us to extend $\xi\in\calZ^\ast$ in a unique way to an element $\tilde{\xi}\in\calV^\ast$. Since also $\Set{z\in\calZ}{\norm{z}_\calV\leq\rho}$ lies dense in $\Set{w\in\calV}{\norm{w}_\calV\leq\rho}$, we finally obtain
		\begin{align*}
			\forall w\in\calV\text{ with }\norm{w}_\calV\leq \rho:\quad 0\geq\chi_\rho(v)+\langle\tilde{\xi},z-v\rangle,\quad\text{i.e.} \quad \tilde{\xi}\in\partial^\calV\chi_\rho(v).
		\end{align*}
	\end{proof}
	
	\begin{lemma}[Attained finite distance]\label{lemma:attaineddist}
		Let $\calZ$ and $\calV$ be the spaces from section~\ref{sec:spaces} and let $\calR:\calZ\to\R_\infty$ be a proper, positively 1-homogeneous, convex and continuous function.
		Then for all $\zeta\in\calZ^\ast$ the distance
		$\dist_{\calV^\ast}\bigl(\zeta,\partial^\calZ\calR(0)\bigr)$
		is attained, if it is finite.
	\end{lemma}
	\begin{proof}
		Let $\zeta\in\calZ^\ast$ be chosen such that $\dist_{\calV^\ast}\bigl(\zeta,\partial^\calZ\calR(0)\bigr)<\infty$. Then we can find a sequence $\left(\sigma_n\right)_{n\in\N}\subset\partial^\calZ\calR(0)$ such that $\left(\sigma_n-\zeta\right)_{n\in\N}\subset\calV^\ast$ and
		\begin{align*}
			\lim_{n\to\infty}\norm{\sigma_n-\zeta}_{\calV^\ast}=\dist_{\calV^\ast}\bigl(\zeta,\partial^\calZ\calR(0)\bigr).
		\end{align*}
		In particular, since $\left(\norm{\sigma_n-\zeta}_{\calV^\ast}\right)_{n\in\N}$ is bounded and $\calV^\ast$ is a reflexive space, we find a weakly convergent subsequence $\left(\sigma_{n_k}-\zeta\right)_{k\in\N}$ in $\calV^\ast$, i.e.
		\begin{align*}
			\sigma_{n_k}-\zeta\rhupk \xi\quad\text{in }\calV^\ast.
		\end{align*}
		By $\calZ\hookrightarrow\calV$ respectively $\calV^\ast\hookrightarrow\calZ^\ast$ this especially implies
		\begin{align*}
			\sigma_{n_k}\rhupk \xi+\zeta\eqcolon\sigma\quad\text{in }\calZ^\ast.
		\end{align*}
		Since $\partial^\calZ\calR(0)$ is weakly closed in $\calZ^\ast$ (as closed and convex set in $\calZ^\ast$; cf. Lemma~\ref{prop:subdifferentialpos1homogen}), we infer that $\sigma\in \partial^\calZ\calR(0)$. Moreover, since $\norm{\cdot}_{\calV^\ast}$ is weakly lower semicontinuous on $\calV^\ast$, we obtain
		\begin{align*}
			\dist_{\calV^\ast}\bigl(\zeta,\partial^\calZ\calR(0)\bigr)=\liminf_{k\to\infty}\norm{\sigma_{n_k}-\zeta}_{\calV^\ast}\geq \norm{\xi}_{\calV^\ast}=\norm{\sigma-\zeta}_{\calV^\ast}.
		\end{align*}
		Thus by $\sigma\in \partial^\calZ\calR(0)$ and the definition of the distance, we indeed have equality here.
	\end{proof}
	%

	\section{Compactness}
	Within this section we collect some compactness theorems. A nice overview of these results is given in \cite{bochner}.
	The following Lemma (cf. \cite[Lemma~8]{Simon1987}), known as Ehrling's or Lions' Lemma, is on the one hand used for the absorption trick within the proof of Proposition~\ref{prop:crucial}. 
	
	\begin{lemma}[Ehrling/Lions]
		\label{lemma:ehrling}
		Let $\calZ$, $\calV$ and $\calX$ be three Banach spaces such that $\calZ$ is compactly embedded in $\calV$ and $\calV$ is continuously embedded in $\calX$, i.e. $\calZ\hookrightarrow\hookrightarrow\calV\hookrightarrow\calX$.
		Then it holds
		\begin{align*}
			\forall\varepsilon>0,\,\exists C(\varepsilon)>0,\,\forall z\in\calZ:\quad \norm{z}_\calV\leq\varepsilon\norm{z}_\calZ+C(\varepsilon)\norm{z}_\calX.
		\end{align*}
	\end{lemma}
    On the other hand, it is the first tool needed to prove the subsequent Aubin-Lions Lemma~\ref{lemma:aubin-lions}. The second tool is given by the following general version of the Arzelà-Ascoli theorem \cite{bochner}, which can be found for instance in \cite[Theorem~7.5.6]{Dieudonne1969}.

	\begin{theorem}[Arzelà-Ascoli]\label{thm:arzela-ascoli}
		Let $\calQ$ be a Banach space and let $K$ be a compact metric space. Further let $M\subset C(K,\calQ)$. Then $M$ is relatively compact in $C(K,\calQ)$ if and only if
		\begin{enumerate}[(a)]
			\item $M$ is equicontinuous, i.e.
			\begin{align*}
				\forall t\in K,\,\forall\varepsilon>0,\,\exists \text{ a neighborhood } U\text{ of } t,\,\forall s\in U,\,\forall f\in M:\quad \norm{f(s)-f(t)}_\calQ\leq\varepsilon,
			\end{align*}
			 and
			\item the set $M(t)\coloneq \Set{f(t)}{f\in M}$ is relatively compact in $\calQ$ for every $t\in K$.
		\end{enumerate}
		In particular, if $\calQ$ is finite dimensional, then $M$ is relatively compact, if and only if $M$ is equicontinuous and bounded.
	\end{theorem}

	The following version of the Aubin-Lions Lemma is due to \cite[Corollary~4]{Simon1987}. 

	\begin{lemma}[Aubin-Lions] \label{lemma:aubin-lions}
			Let $\calZ$, $\calV$ and $\calX$ be Banach spaces satisfying the assumptions of Ehrling's Lemma~\ref{lemma:ehrling}. Further let $p,q\in\left[1,\infty\right]$ and $S>0$. Consider the space
			\begin{align*}
				W^{1,p,q}\left(\left[0,S\right];\calZ,\calX\right)\coloneq \Set{z\in L^{p}\left(\left[0,S\right];\calZ\right)}{\frac{\dd{z}}{\dd{t}}\in L^{q}\left(\left[0,S\right];\calX\right)}
			\end{align*}
		endowed with the norm
		\begin{align*}
			\norm{z}_{	W^{1,p,q}\left(\left[0,S\right];\calZ,\calX\right)}\coloneq \norm{z}_{L^{p}\left(\left[0,S\right];\calZ\right)}+\norm{\frac{\dd{z}}{\dd{t}}}_{L^{q}\left(\left[0,S\right];\calX\right)},
		\end{align*}
		where $\frac{\dd{z}}{\dd{t}}$ is the time derivative of $z$ in the sense of vector valued distributions.
		\begin{enumerate}[(a)]
			\item If $p<\infty$, then the embedding of $W^{1,p,q}\left(\left[0,S\right];\calZ,\calX\right)$ into $L^{p}\left(\left[0,S\right];\calV\right)$ is compact.
			\item If $p=\infty$ and $q>1$, then the embedding of $W^{1,p,q}\left(\left[0,S\right];\calZ,\calX\right)$ into $C\left(\left[0,S\right];\calV\right)$ is compact.
		\end{enumerate}		
	\end{lemma}

	The Banach-Alaoglu theorem states that the closed unit ball in the dual space of a Banach space is compact in the weak$*$ topology \cite[Theorem~3.16]{Brezis1973}. Moreover, this unit ball is metrizable in the weak$*$ topology, if and only if the Banach space is separable \cite[Theorem~3.28]{Brezis2011}. The Bochner spaces $L^p(\left[0,S\right];\calQ)$ are separable for $\calQ$ a separable Banach space and $1\leq p<\infty$ \cite[Korollar~2.11]{bochner}. Finally, in metric spaces sequentially compactness and compactness with respect to the open cover definition are equivalent. This gives us the following statement concerning sequentially weak$*$ compactness of $L^p(\left[0,S\right];\calQ)$.
	\begin{theorem}[Sequentially weak$*$ compactness]\label{thm:weakastcompactness}
		Let $\calQ$ be a separable Banach space and let $1\leq p<\infty$ and $S>0$. Then the closed unit ball in $\left(L^p(\left[0,S\right];\calQ)\right)^\ast$
		\begin{align*}
			\overline{B}_1^{\left(L^p(\left[0,S\right];\calQ)\right)^\ast}\coloneq \Set{f\in \left(L^p(\left[0,S\right];\calQ)\right)^\ast}{\norm{f}_{\left(L^p(\left[0,S\right];\calQ)\right)^\ast}\leq 1}
		\end{align*}
		is sequentially weakly$*$ compact. Especially, for every bounded sequence $\left(f_n\right)_{n\in\N}\subset \left(L^p(\left[0,S\right];\calQ)\right)^\ast$ there exists $f\in \left(L^p(\left[0,S\right];\calQ)\right)^\ast$ and a subsequence $\left(f_{n_k}\right)_{k\in\N}\subset\left(f_n\right)_{n\in\N}$ such that we have for $k\to\infty$
		\begin{align*}
			f_{n_k}\rhupastk f\quad\text{in } \left(L^p(\left[0,S\right];\calQ)\right)^\ast.
		\end{align*}
	\end{theorem}
\noindent
	The value of the previous Theorem becomes apparent by the following characterization of the dual spaces $\left(L^p(\left[0,S\right],\calQ)\right)^\ast$.
	\begin{lemma}[Characterization of the dual spaces]\label{thm:characterizationofdualspaces}
		Let $\calQ$ be a reflexive and separable Banach space and let $1\leq p<\infty$. Then the mapping
		\begin{align*}
			T:L^{p'}\left(\left[0,S\right];\calQ^\ast\right)&\to \left(L^{p}\left(\left[0,S\right];\calQ\right)\right)^\ast,\\
			f&\mapsto T_f\text{ with }\langle T_f,v\rangle_{\left(L^{p}\left(\left[0,S\right];\calQ\right)\right)^\ast,L^{p}\left(\left[0,S\right];\calQ\right)}\coloneq \int_{0}^{S}\langle f(t), v(t)\rangle_{\calQ^\ast,\calQ}\dd{t}
		\end{align*}
		is an isometric isomorphism.
	\end{lemma}
	\begin{remark}\label{remark:dualspaces}		
		Given a Banach space $\calQ$ note that we have the equivalence
		\begin{align*}
			\calQ\text{ reflexive and separable}\quad\Leftrightarrow\quad \calQ^\ast\text{ reflexive and seperable}.
		\end{align*}
		Indeed, the spaces $\calU^{\tilde{p}}$, $\calV$ and $\calZ$ from section~\ref{sec:set-up and notation} are reflexive and separable, such that we can apply Theorem~\ref{thm:weakastcompactness} and use the characterization of the dual spaces from Lemma~\ref{thm:characterizationofdualspaces} within the proof and the statement of Theorem~\ref{thm:convergence}.
	\end{remark}

	\section{Lower semicontinuity properties}
	
	For the following two lemmas on lower semicontinuity properties we refer to \cite[Appendix~B]{Knees2018} and the references therein.
	\begin{lemma}\label{Knees18:lemmaB1}
		Let $\left(v_n\right)_{n\in\N}\subset L^\infty((0,S);\calV)$ and $v\in L^\infty((0,S);\calV)$ with $v_n\rhupastn v$ in $L^\infty((0,S);\calV)$. Further let $\left(\delta_n\right)_{n\in\N}\subset L^1((0,S);\left[0,\infty\right))$and $\delta\in L^1((0,S);\left[0,\infty\right))$ with $\liminf_{n\to\infty}\delta_n(s)\geq \delta(s)$ f.a.a. $s\in (0,S)$. Then
		\begin{align*}
			\liminf_{n\to\infty} \int_{0}^{S}\norm{v_n(s)}_\calV \delta_n(s)\dd{s}\geq \int_{0}^{S}\norm{v(s)}_\calV \delta(s)\dd{s}.
		\end{align*}
	\end{lemma}
	\begin{lemma}\label{Knees18:lemmaB2}
		For all $n\in\N$ let $\delta_n,\delta,\tau_n,\tau:(0,S)\to\left[0,\infty\right)$ be measurable functions satisfying $\liminf_{n\to\infty}\delta_n(s)\geq \delta(s)$ f.a.a. $s\in (0,S)$ and $\tau_n \rhupn \tau$ weakly in $L^1((0,S);\R)$. Then
		\begin{align*}
			\liminf_{n\to\infty}\int_{0}^{S}\tau_n(s)\delta_n(s)\dd{s}\geq \int_{0}^{S}\tau(s)\delta(s)\dd{s}.
		\end{align*}
	\end{lemma}

	In view of Lemma~\ref{lemma:be_AM} the following closedness property is useful.
	\begin{lemma}[Weak closedness of the relevant range for the phase-field variable]\label{lemma:closedness}
		The physically relevant range $\calZ_{\left[0,1\right]}$ for the phase-field variable
		is weakly closed in $\calZ$.
	\end{lemma}
	\begin{proof}
		Let $\left(z_n\right)_{n\in\N}\subset \calZ$ and $z\in\calZ$ satisfy
		\begin{align*}
			\forall n\in\N:\,z_n\in\left[0,1\right]\text{ a.e. in }\Omega \quad\text{and}\quad z_n\rhupn z\text{ in }\calZ.
		\end{align*}
		By the compact embedding $\calZ\hookrightarrow\hookrightarrow L^1(\Omega)$ we have $z_n\rton z$ in $L^1(\Omega)$. Now suppose there exists $\tilde{\Omega}\subset\Omega$ with $\big\vert\tilde{\Omega}\big\vert>0$ such that $z\notin\left[0,1\right]$ on $\tilde{\Omega}$. Then there exists $c>0$ such that
		\begin{align*}
			0<c<\int_{\tilde{\Omega}}\abs{z-z_n}\dd{x}\leq \int_{\Omega}\abs{z-z_n}\dd{x}=\norm{z-z_n}_{L^1(\Omega)}\rton 0,
		\end{align*}
		a contradiction.
	\end{proof}

	The following lower semicontinuity result is a variant of \cite[Lemma~3.2.18]{Sievers2020} and is needed to treat  the unbounded dissipation.
	\begin{lemma}[Lower semicontinuity]\label{lemma:lowersemicon}
		Let $\left(\xi_n\right)_{n\in\bbN}\subset \calZ^\ast$ and $\xi\in\calZ^\ast$ with $\liminf\limits_{n\to\infty}\left\langle \xi_n,-v\right\rangle\geq\langle \xi,-v\rangle$ for all $v\in\calZ$ with $\calR(v)<\infty$ and let the distance terms be uniformly bounded, i.e. $\dist_{\calV^\ast}(-\xi_n,\partial^\calZ\calR(0))\leq C$ with $C>0$ independent of $n\in\bbN$. Then the following weak lower semicontinuity result holds true
		\begin{align}
			\liminf_{n\to\infty}\dist_{\calV^\ast}(-\xi_n,\partial^\calZ\calR(0))\geq \dist_{\calV^\ast}(-\xi,\partial^\calZ\calR(0)).
		\end{align}
	\end{lemma}
	\begin{proof}
		At first we choose a subsequence $\left(n_k\right)_{k\in\bbN}$, such that
		\begin{align*}
			\lim_{k\to\infty} \dist_{\calV^\ast}\big(-\xi_{n_k},\partial^\calZ\calR(0)\big)=\liminf_{n\to\infty} \dist_{\calV^\ast}(-\xi_n,\partial^\calZ\calR(0)).
		\end{align*}
		By assumption the distance terms are finite, hence by Lemma~\ref{lemma:attaineddist} the distances are attained, i.e. for all $k\in\N$ there exists a $\sigma_k\in\partial^\calZ\calR(0)$, such that
		\begin{align*}
			\big\Vert\xi_{n_k}+\sigma_k\big\Vert_{\calV^\ast} =\dist_{\calV^\ast}\big(-\xi_{n_k},\partial^\calZ\calR(0)\big)\leq C,
		\end{align*}
		where $C$ is the uniform bound from the assumption.
		Hence, we can find a subsequence (not relabeled) $\left(\xi_{n_k}+\sigma_k\right)_{k\in\bbN}$ with
		\begin{align*}
			\zeta_k\coloneq \xi_{n_k}+\sigma_k\rhupk\zeta\quad\text{in }\calV^\ast.
		\end{align*}
		By assumption it is $\liminf\limits_{n\to\infty} \langle \xi_n,-v\rangle\geq \langle \xi,-v\rangle$, thus $\zeta_k\rhupk\zeta$ in $\calV^\ast$ implies for all $v\in\calZ$ with $\calR(v)<\infty$
		\begin{align*}
			\liminf_{k\to\infty} \langle \sigma_k,v\rangle&=\liminf_{k\to\infty}\langle \zeta_k-\xi_{n_k},v\rangle=\underbrace{\lim_{k\to\infty}\langle \zeta_k,v\rangle}_{=\langle \zeta,v\rangle}
			+\underbrace{\liminf_{k\to\infty}\langle \xi_{n_k} ,-v\rangle}_{\geq \langle \xi,-v\rangle}\\
			&\geq \langle \zeta,v\rangle+\langle\xi,-v\rangle=\langle\zeta-\xi,v\rangle.
		\end{align*}
		Together with $\sigma_k\in\partial^\calZ\calR(0)$, i.e. $\calR(v)\geq\langle \sigma_k,v\rangle$ for all $v\in\calZ$, it follows by taking the limes inferior
		\begin{align*}
			\calR(v)\geq\langle \zeta-\xi,v\rangle\quad\text{for all } v\in\calZ,\quad\text{i.e. } \sigma\coloneq \zeta-\xi\in\partial^\calZ\calR(0) \text{  by Lemma~\ref{prop:subdifferentialpos1homogen}}.
		\end{align*}
		Therefore, we have shown that
		\begin{align*}
			\xi_{n_k}+\sigma_k\rhupk\xi+\sigma\quad\text{in }\calV^\ast\text{ with } \sigma\in\partial^\calZ\calR(0).
		\end{align*}
		Now by $\sigma\in\partial^\calZ\calR(0)$ and the weak lower semicontinuity of $\norm{\cdot}_{\calV^\ast}$ it follows
		\begin{align*}
			\dist_{\calV^\ast}(-\xi,\partial^\calZ\calR(0))&\leq \norm{\xi+\sigma}_{\calV^\ast} =\norm{\zeta}_{\calV^\ast}\leq \liminf_{k\to\infty}\big\Vert\xi_{n_k}+\sigma_k\big\Vert_{\calV^\ast}\\
			&=\liminf_{k\to\infty}\dist_{\calV^\ast}\big(-\xi_{n_k},\partial^\calZ\calR(0)\big)=\liminf_{n\to\infty}\dist_{\calV^\ast}(-\xi_n,\partial^\calZ\calR(0)).
		\end{align*}
	\end{proof}

	\printbibliography
	
\end{document}